\documentclass[11pt]{article}
\usepackage{amssymb}
\usepackage{amscd}
\usepackage{amsmath}
\input amssym.def
\input amssym  

\newtheorem{theorem}{Theorem}[section]
\newtheorem{corollary}[theorem]{Corollary}
\newtheorem{lemma}[theorem]{Lemma}

\newtheorem{conjecture}[theorem]{Conjecture}

\newtheorem{remark}[theorem]{Remark}

\usepackage{amssymb} 
\title{Presentations of the principal
subspaces of the higher-level standard
$\widehat{\mathfrak{sl}(3)}$-modules} \author{ Christopher Sadowski}
\date{}
\begin{document}
\maketitle

\renewcommand{\theequation}{\thesection.\arabic{equation}}
\renewcommand{\thetheorem}{\thesection.\arabic{theorem}}
\setcounter{equation}{0} \setcounter{theorem}{0}
\setcounter{section}{0}

\begin{abstract}
Using the theory of vertex operator algebras and intertwining
operators, we obtain presentations for the principal subspaces of all
the standard $\widehat{\goth{sl}(3)}$-modules. Certain of these
presentations had been conjectured and used in work of Calinescu to
construct exact sequences leading to the graded dimensions of certain
principal subspaces. We prove the conjecture in its full generality
for all standard $\widehat{\goth{sl}(3)}$-modules. We then provide a conjecture
for the case of $\widehat{\goth{sl}(n+1)}$.
\end{abstract}      

\section{Introduction}

 There is a long-standing connection between the theories of vertex
 operators and vertex operator algebras (\cite{B}, \cite{FLM},
 \cite{LL}, etc.) and affine Lie algebras (cf. \cite{K}) on the one
 hand, and Rogers-Ramanujan-type combinatorial identities
 (cf. \cite{A}) on the other hand (\cite{LM}, \cite{LW1}--\cite{LW4},
 \cite{LP1}--\cite{LP2}, and many other references). We begin by sketching a
 brief history of some of these connections and the importance of
 giving a priori proofs of presentations of certain algebraic
 structures.
 
Many difference-two type partition conditions have been interpreted
and obtained by the study of certain natural substructures of standard
(i.e., integrable highest weight) modules for affine Lie algebras. In
particular, in \cite{FS1}--\cite{FS2}, Feigin and Stoyanovsky
introduced the notion of ``principal subspace" of a standard module
for an affine Lie algebra, and in the case of $A_1^{(1)}(=\widehat{\goth{sl}(2)})$ and
$A_2^{(1)}(=\widehat{\goth{sl}(3)})$ obtained, under certain assumptions (in particular,
presentations for these principal subspaces in terms of generators and
relations) the graded dimensions (``characters") of the principal
subspaces of the ``vacuum" standard modules. Interestingly enough,
these graded dimensions were related to the Rogers-Ramanujan partition
identities, and more generally, the Gordon-Andrews identities, but in
a different setting than the original vertex-algebraic interpretation
of these identities in \cite{LW2}--\cite{LW4}. A more general case was
considered by Georgiev in \cite{G1}, where combinatorial bases were
constructed for the principal subspaces associated to certain standard
$A_n^{(1)}$-modules. Using these bases, Georgiev obtained the graded
dimensions of these principal subspaces.

Then, in \cite{CLM1}--\cite{CLM2}, the authors addressed the problem
of vertex-algebraically interpreting the classical Rogers-Ramanujan recursion and, more
generally, the Rogers-Selberg recursions (cf. \cite{A}) by using
intertwining operators among modules for vertex operator algebras to
construct exact sequences leading to these recursions. In particular,
the solutions of these recursions gave the graded dimensions of the
principal subspaces of the standard $A_1^{(1)}$-modules. In
\cite{CalLM1}--\cite{CalLM2}, the authors gave an a priori proof,
again using intertwining operators, of
the presentations assumed in \cite{FS1} and \cite{CLM1}--\cite{CLM2}
(these presentations had been known to be true, but it was important
to seek an a priori proof; see below). In \cite{CalLM3} these results
were extended to the 
principal subspaces of the standard level $1$ modules for the $ADE$-type
untwisted affine Lie algebras. The desired presentations
were proved, and exact sequences were obtained leading to recursions
and the graded dimensions of the principal subspaces of the level $1$
standard modules.

More specifically, principal subspaces of standard modules can be
viewed as the quotients of the universal enveloping algebra of a
certain nilpotent subalgebra of an affine Lie algebra by the left
ideals, annihilating the highest weight vectors,
generated by certain natural elements (the generators and
relations mentioned earlier). In \cite{CLM1}--\cite{CLM2}, the authors
assumed (as was done in \cite{FS1}--\cite{FS2}) that these natural
elements generated the annhilating ideals --- the presentations
for the principal subspaces of standard $A_1^{(1)}$-modules.
These presentations followed from the knowledge of
combinatorial bases constructed in \cite{LP2} (which was also invoked
in \cite{FS1}), but an a priori proof of these presentations was
lacking at the time.  Using intertwining operators, the authors
constructed exact sequences among the principal subspaces to obtain
the Rogers-Ramanujan recursion and the Rogers-Selberg recursions,
whose solutions gave the graded dimensions of the principal
subspaces. A key part of the proof of exactness of these sequences
followed directly from knowledge of the presentations of the principal
subspaces. The important problem of proving these presentations a
priori (without knowledge of combinatorial bases) was later solved in
\cite{CalLM1}--\cite{CalLM2}, and this work was further generalized in
\cite{CalLM3}. Knowledge of these presentations (which look relatively simple but
are highly nontrivial to prove) has been enough to obtain graded dimensions
of principal subspaces, without knowledge of combinatorial bases,
 which are significantly more complicated. Moreover,
in this theory, one wants to construct combinatorial bases as a
consequence of the theory and not be required to use such bases in
order to justify steps, notably, presentations of principal subspaces,
needed in the proof of their construction.

In \cite{C3}, Calinescu conjectured presentations for principal
subspaces of certain higher level standard $A_2^{(1)}$-modules and
gave exact sequences which led to the graded dimensions of the
principal subspaces considered. In \cite{C4}, she obtained
presentations for the principal subspaces of the level $1$ standard
$A_n^{(1)}$-modules and constructed exact sequences among these
modules, using them to obtain the graded dimensions of these principal
subspaces. A new method for proving the presentations in the $A_n^{(1)}$ level $1$ case
 was given in
\cite{CalLM3}, where, as we mentioned above, such presentations were proved a priori for the
principal subspaces of the standard modules of the untwisted affine Lie
algebras of $ADE$--type. These presentations were used to obtain exact
sequences, leading to recursions and the graded dimensions of the
principal subspaces.

The present work is another step forward in the same spirit. We
exploit intertwining operators among vertex operator algebra modules to solve
the problem of giving an a priori proof of presentations for the
principal subspaces of all the standard modules for $A_2^{(1)}(=\widehat{\goth{sl}(3)})$, including
those assumed conjecturally and used in \cite{C3}. The methods used in the proof of these
presentations are similar to those in
\cite{CalLM1}--\cite{CalLM3}, in that certain minimal counterexamples
are postulated and shown not to exist. However, in the general case, we needed to
introduce certain new ideas to prove our presentations. We then proceed to
formulate the presentations for principal subspace of all the standard modules
for $A_n^{(1)}$ as a conjecture.

In \cite{P} and \cite{MiP}, the authors considered the
integral lattice case (a case more general than the Lie algebra root and weight
lattice cases considered in \cite{CLM1}-\cite{CLM2},
\cite{CalLM1}-\cite{CalLM3}) and constructed combinatorial bases for
principal subspaces in this context. The problem of finding
combinatorial bases for principal subspaces of certain standard
modules in the $B_{2}^{(1)}$ case has been solved in \cite{Bu}. In a similar vein, the
problem of finding combinatorial bases for suitably defined principal
subspaces in the quantum $A_{n}^{(1)}$ case has been solved in
\cite{Ko}. Work in similar directions concerning combinatorial bases
for certain ``commutative" variants of principal subspaces has also
been pursued in \cite{Pr},  \cite{J1}--\cite{J3},
\cite{T1}--\cite{T4}, \cite{Ba}, \cite{JPr}. A different variant of principal subspace was
considered in \cite{AKS} and \cite{FFJMM}. In \cite{AKS}, the authors cite
well-known presentations for standard modules, and use these to provide (without proof) a 
set of defining relations for each principal subspace.
In 
\cite{FFJMM}, in which the authors consider
$A_2^{(1)}$, they do indeed prove that certain relations form a set 
of defining relations for their variant of principal subspace.
 In the case of the vacuum
modules, the principal subspaces in \cite{FFJMM} are essentially identical to the principal subspaces
considered in the present work, and their defining
relations indeed agree with those in the present work.
For the non-vacuum modules, the principal subspaces considered in the present work
can be viewed as proper substructures of those considered in \cite{FFJMM},
 and correspondingly, the defining relations
we obtain are different. Our method for proving the completeness of our defining
relations is completely different from the method in \cite{FFJMM}.

We now give a brief overview of our main result.  Given a complex
semisimple Lie algebra $\goth{g}$, a fixed Cartan subalgebra
$\goth{h}$, a fixed set of positive roots $\Delta_+$, and a root
vector $x_\alpha$ for each $\alpha \in \Delta_+$, consider the
subalgebra $\goth{n} = \coprod_{\alpha \in \Delta_+}
\mathbb{C}x_\alpha \subset \goth{g}$ spanned by the positive root
vectors.  The affinization $\bar{\goth{n}} = \goth{n} \otimes
\mathbb{C}[t,t^{-1}]$ of $\goth{n}$ (without any central extension) is
a subalgebra of the affine Lie algebra $\widehat{\goth{g}}$ (with
central extension). Let $L(\Lambda)$ be the standard module of
$\widehat{\goth{g}}$ with highest weight $\Lambda$ and level $k$, a
positive integer, and let $v_\Lambda \in L(\Lambda)$ be a highest
weight vector. The principal subspace of $L(\Lambda)$ is defined by
\begin{eqnarray}
W(\Lambda) = U(\bar{\goth{n}}) \cdot v_\Lambda,
\end{eqnarray}
where $U(\cdot)$ is the universal enveloping algebra.  Consider the
natural surjective map
\begin{eqnarray}
f_\Lambda : U(\bar{\goth{n}}) \rightarrow W(\Lambda).
\end{eqnarray}

We now take $\goth{g} = \goth{sl}(3)$.
We precisely determine Ker$f_\Lambda$ in terms of certain natural left
ideals of $U(\bar{\goth{n}})$. Specifically, in terms of the fundamental weights
of $A_2^{(1)}$, which we label $\Lambda_0, \Lambda_1,$ and
$\Lambda_2$, we may express $\Lambda$ as
\begin{eqnarray*}
\Lambda = k_0 \Lambda_0 + k_1 \Lambda_1 + k_2\Lambda_2,
\end{eqnarray*}
for some nonnegative integers $k_0, k_1,$ and $k_2$.  We define an
ideal $I_{k \Lambda_0}$ in terms of left ideals generated by the
coefficients of certain vertex operators associated with singular vectors
in a natural way. This left ideal is then used to define a larger left
ideal
\begin{eqnarray*}
I_\Lambda = I_{k\Lambda_0} + 
U(\bar{\goth{n}})x_{\alpha_1}(-1)^{k_0+k_2+1} +
U(\bar{\goth{n}})x_{\alpha_2}(-1)^{k_0+k_1+1} +
U(\bar{\goth{n}})x_{\alpha_1 + \alpha_2}(-1)^{k_0+1},
\end{eqnarray*}
where we use $x(n)$ to denote the action of $x \otimes t^n \in 
\widehat{\goth{g}}$ for $x \in \goth{g}$ and $n \in \mathbb{Z}$.
We then proceed to show that
\begin{eqnarray*}
\mathrm{Ker}f_{ \Lambda}= I_\Lambda.
\end{eqnarray*}
The proof of this result is similar in structure to the proof of the
presentations in \cite{CalLM2}. Considering all dominant integral
weights together, we choose minimal counterexamples (certain elements
in $\mathrm{Ker}f_{ \Lambda} \setminus I_\Lambda$) and show that a
contradiction is reached for each $\Lambda$. Certain maps used in
\cite{CalLM3} are also generalized and used in the proof, but these
ideas do not extend to the most general case. We develop a method 
for reaching the desired contradictions for each $\Lambda$ which
``rebuilds'' the minimal counterexample to show that it is in fact
an element of $I_\Lambda$. This ``rebuilding'' technique can also be used
to show all of the presentations proved in the works \cite{CalLM1}-\cite{CalLM3}
in the type $A$ case with suitable modifications (see remarks at the end
of Section 4).

In \cite{C3}, certain of these presentations were conjectured and used
to construct exact sequences among principal subspaces. Using these
exact sequences, Calinescu obtained the previously unknown graded
dimensions for principal subspaces whose highest weights are of the
form $k_1\Lambda_1 + k_2\Lambda_2$, where $k_1, k_2$ are positive
integers. The problem of constructing exact sequences for more general
highest weights is still unsolved.

The present work is organized as follows: In Section 2 we recall the vertex-algebraic
construction of the standard $A_n^{(1)}$-modules. We also recall
various results on intertwining operators (\cite{FHL}, \cite{DL}),
which will be used throughout, and we specialize to the case $n=2$. In
Section 3 we recall the definition of principal subspace from
\cite{FS1}--\cite{FS2} and define certain natural left ideals
$I_\Lambda$ for each dominant integral weight $\Lambda$ of
$A_2^{(1)}$. In Section 4 we prove our main result, showing that
Ker$f_\Lambda = I_{\Lambda}$. In Section 5 we formulate, as a conjecture,
the corresponding presentations 
of the principal subspaces of the standard $A_n^{(1)}$-modulues and relate
this conjecture to the presentations proved in this work.

{\bf Acknowledgement} We are very grateful to James Lepowsky for many helpful
discussions and suggestions for improvements to this work. 
We would also like to thank Corina Calinescu for
discussions which led to some of the ideas which are found here.

\section{Preliminaries}
\setcounter{equation}{0}

We begin by recalling certain vertex-algebraic constructions for the
untwisted affine Lie algebra $\widehat{\goth{sl}(n+1)}$, $n$ a
positive integer. We shall be working in the setting of \cite{FLM} and
\cite{LL}.

Fix a Cartan
subalgebra $\goth{h}$ of $\goth{sl}(n+1)$. Also fix a set of roots
$\Delta$, a set of simple roots $\{\alpha_1,\dots ,\alpha_n\}$, and a
set of positive roots $\Delta_+$. Let $\langle \cdot, \cdot \rangle$
denote the Killing form, rescaled so that $\langle \alpha, \alpha
\rangle =2$ for each $\alpha \in \Delta$. Using this form, we
identify $\goth{h}$ with $\goth{h}^*$. Let $\lambda_1, \dots,
\lambda_n \in \goth{h} \simeq \goth{h}^*$ denote the fundamental
weights of $\goth{sl}(n+1)$. Recall that $\langle \lambda_i,\alpha_j
\rangle = \delta_{ij}$ for each $i,j=1,\dots, n$.  Denote by $Q=
\sum_{i=1}^n \mathbb{Z}\alpha_i$ and $P=\sum_{i=1}^n
\mathbb{Z}\lambda_i$ the root lattice and weight lattice of
$\goth{sl}(n+1)$, respectively.

For each root $\alpha \in \Delta$, we have a root vector $x_{\alpha} \in \goth{sl}(n+1)$
(recall that $[h,x_\alpha] = \langle \alpha,h \rangle x_\alpha$ for each $h \in \goth{h}$).
We define
\begin{eqnarray*} 
\goth{n}= \sum_{\alpha \in \Delta_+} \mathbb{C} x_\alpha,
\end{eqnarray*}
a nilpotent subalgebra of $\goth{sl}(n+1)$.

We have the corresponding untwisted affine Lie algebra given by 
\begin{eqnarray*}
\widehat{{\goth{sl}(n+1)}}= {\goth{sl}(n+1)} \otimes \mathbb{C}[t, t^{-1}] \oplus \mathbb{C}c,
\end{eqnarray*}
where $c$ is a non-zero central element and 
\begin{eqnarray*}
[ x \otimes t^m, y \otimes t^p ] = [x, y] \otimes t^{m+p} + m\langle x, y \rangle \delta _{m+p, 0} c
\end{eqnarray*}
for any $x, y \in {\goth{sl}(n+1)}$ and $m, p \in \mathbb{Z}$.  If we
adjoin the degree operator $d$, where $$[d, x \otimes t^m]=mx\otimes
t^m$$ $$[d,c]=0,$$ we obtain the affine Kac-Moody Lie algebra
$\widetilde{{\goth{sl}(n+1)}}=\widehat{{\goth{sl}(n+1)}} \oplus
\mathbb{C}d$ (cf. \cite{K}). We define two important subalgebras of
$\widehat{\goth{sl}(n+1)}$:
\begin{eqnarray*}
\widehat{\goth{h}} = \goth{h} \otimes \mathbb{C}[t, t^{-1}] \oplus
\mathbb{C}c
\end{eqnarray*}
and the Heisenberg subalgebra
\begin{eqnarray*}
  \widehat {\goth{h}}_{\mathbb{Z}} = \coprod _{m \in \mathbb{ Z}
    \setminus \{0\}} \goth{h} \otimes t^m \oplus \mathbb{C}c
\end{eqnarray*}
(in the notation of \cite{FLM}, \cite{LL}).
We extend our form $\langle \cdot,\cdot \rangle$ to $\goth{h} \oplus \mathbb{C}c \oplus \mathbb{C}d $ by defining 
\begin{eqnarray*}\langle c, c\rangle&=&0\\
 \langle d, d \rangle &=&0 \\
 \langle c, d \rangle &=&1.
\end{eqnarray*}
Using this form,
we may identify $\goth{h} \oplus \mathbb{C}c \oplus \mathbb{C}d$ with
$(\goth{h} \oplus \mathbb{C}c \oplus \mathbb{C}d)^{*}$.  The simple
roots of $\widehat{{\goth{sl}(n+1)}}$ are $\alpha _{0}$, $\alpha
_{1},\dots , \alpha _{n}$ and the fundamental weights of
$\widehat{{\goth{sl}(n+1)}}$ are $\Lambda_{0}, \Lambda _{1},\dots,
\Lambda _{n} $, given by
\begin{eqnarray*}
\alpha_0= c-(\alpha_1+\alpha_2 + \dots + \alpha_n)
\end{eqnarray*}
and
\begin{eqnarray*}
 \Lambda_0= d, \; \Lambda_i=\Lambda_0+\lambda_i
\end{eqnarray*}
for each $i=1, \dots, n$.

An $\widehat{{\goth{sl}(n+1)}}$-module $V$ is said to have level $k
\in \mathbb{C}$ if the central element $c$ acts as multiplication by
$k$ (i.e. $c \cdot v = kv$ for all $v \in V$). Any standard
(i.e. irreducible integrable highest weight) module $L(\Lambda)$ with
$\Lambda \in ( \goth{h} \oplus \mathbb{C}c \oplus \mathbb{C}d)^{*}$
has nonnegative integral level, given by $\langle \Lambda, c \rangle$
(cf. \cite{K}). Let $L(\Lambda _{0})$, $ L( \Lambda _{1}), \dots,
L(\Lambda_{n})$ denote the standard $ \widehat
{{\goth{sl}(n+1)}}$-modules of level $1$ with $v_{\Lambda_0}$,
$v_{\Lambda_1}, \dots, v_{\Lambda_n}$ as highest weight vectors,
respectively.

Continuing to work in the setting of \cite{FLM} and \cite{LL},
we now recall the lattice vertex operator construction of the level
$1$ standard modules for $\widehat{\goth{sl}(n+1)}$. We use $U(\cdot)$ to denote the universal enveloping algebra.
 The induced module
\begin{eqnarray*} 
M(1)= U( \widehat{\goth{h}}) \otimes _{U(\goth{h} \otimes \mathbb{C}[t] \oplus \mathbb{C}c)}\mathbb{C}
\end{eqnarray*}
has a natural $\widehat{\goth{h}}$-module structure, where
$\goth{h}\otimes \mathbb{C}[t]$ acts trivially and $c$ acts as
identity on the one-dimensional module $ \mathbb{C}.$ Let
$s=2(n+1)^2$. We fix a primitive $s^\mathrm{th}$ root of unity
$\nu_s$, and a central extension $\widehat{P}$ of the weight lattice
$P$ by the finite cyclic group $\langle \kappa\rangle= \langle \kappa
\; | \: {\kappa}^s=1 \rangle$ of order $s$,
\begin{eqnarray*}
1 \rightarrow \langle \kappa \rangle \rightarrow \widehat{P}
\bar{\longrightarrow} P \rightarrow 1
\end{eqnarray*}
with associated commutator map $c_0: P \times P \longrightarrow
\mathbb{Z}/s\mathbb{Z}$, defined by
$aba^{-1}b^{-1}={\kappa}^{c_{0}(\bar{a},\bar{b})}$ for $a,b \in
\widehat{P}$. Let $c:P \times P \longrightarrow \mathbb{C}^{\times}$
denote the alternating $\mathbb{Z}$-bilinear map defined by
$c(\lambda, \mu)=\nu_s^{c_0(\lambda, \mu)}$ for $\lambda, \mu \in
P$. We require that
\begin{eqnarray*} \label{com}
c(\alpha, \beta)=(-1)^{\langle \alpha, \beta \rangle} \; \; \;
\mbox{for} \; \; \; \alpha, \beta \in Q.
\end{eqnarray*}
Such a central extension $\widehat{P}$ of $P$ does indeed exist (see Remark 6.4.12 in \cite{LL}).

We define the faithful character $\chi: \langle \kappa \rangle
\longrightarrow \mathbb{C}^{\times}$ by $\chi(\kappa)= \nu_s$. Let
$\mathbb{C}_{\chi}$ be the one dimensional $\langle \kappa
\rangle$-module, where the action of $\kappa$ is given by $ \kappa
\cdot 1=\nu_s$, and form the induced $\widehat{P}$-module
\begin{eqnarray*}
\mathbb{C} \{ P \} = \mathbb{C}[\widehat{P}]\otimes_{\mathbb{C}[\langle \kappa \rangle]}\mathbb{C}_{\chi}.
\end{eqnarray*}
For any subset $E \subset P$, we define $\widehat{E} = \{ a \in P |
\bar{a} \in E\}$, and we form $\mathbb{C} \{ E \}$ in the obvious way.
Then, the space
\begin{eqnarray*}
V_Q=M(1) \otimes \mathbb{C}\{ Q \}
\end{eqnarray*}
carries a natural vertex operator algebra structure, with
$1$ as vacuum vector, and the space 
\begin{eqnarray*}
V_P=M(1) \otimes \mathbb{C}\{ P \}
\end{eqnarray*}
is naturally a $V_Q$-module. 
We now recall some important details of this construction (cf. \cite{LL}).  

Choose a section 
\begin{eqnarray} \label{section}
e: P & \longrightarrow & \widehat{P} \\ \nonumber
\alpha & \mapsto & e_{\alpha}, \nonumber
\end{eqnarray}
(i.e. a map which satisfies $\ \bar{ } \circ e = 1$) such that $e_0 =
1$. Let $\epsilon_0: P \times P \longrightarrow
\mathbb{Z}/s\mathbb{Z}$ the corresponding $2$-cocycle, defined by the
condition $e_{\alpha}e_{\beta}={\kappa}^{\epsilon_0(\alpha, \beta)}
e_{\alpha+\beta}$ for $\alpha, \beta \in P $ and define the map
$\epsilon:P\times P \longrightarrow \mathbb{C}^{\times}$ by
$\epsilon(\alpha,\beta)=\nu_s^{\epsilon_0(\alpha, \beta)}$ For any
$\alpha,\beta \in P$ we have
\begin{eqnarray} \label{e-c}
\epsilon(\alpha, \beta)  / \epsilon(\beta, \alpha)= c(\alpha, \beta)
\end{eqnarray}
and 
\begin{eqnarray} \label{epsilon-zero}
\epsilon(\alpha, 0)=\epsilon(0, \alpha)=1.
\end{eqnarray}
We use this choice of section (\ref{section}) identify
$\mathbb{C}\{P \}$ and the group algebra $\mathbb{C}[P]$. In
particular, we have a vector space isomorphism given by
\begin{eqnarray} \label{identification}
\mathbb{C}[P] & \longrightarrow & \mathbb{C}\{ P \} \\ 
e^{\alpha} & \mapsto & \iota(e_{\alpha}) \nonumber
\end{eqnarray}
for $\alpha \in P$, where, for $a \in \widehat{P}$, we set $\iota(a)=a
\otimes 1 \in \mathbb{C} \{ P \}$.  By restriction, we also have the
identification $\mathbb{C}[Q] \simeq \mathbb{C}\{ Q \}$.  There is a
natural action $\widehat{P}$ on $\mathbb{C}[P]$ given by
\begin{eqnarray*}
e_{\alpha} \cdot e^{\beta}&=&
{\epsilon(\alpha, \beta)}e^{\alpha + \beta},\\ 
\kappa \cdot e^{\beta}&=&
\nu_{s} e^{\beta}
\end{eqnarray*}
for $\alpha, \beta \in P.$ As operators on $\mathbb{C}[P] \simeq \mathbb{C}\{P \}$ we have
\begin{eqnarray} \label{e-multiplication}
e_{\alpha}e_{\beta}=\epsilon(\alpha, \beta)e_{\alpha+\beta}.
\end{eqnarray}
We make the identifications 
$$ V_P = M(1) \otimes \mathbb{C}[P],$$
$$V_Q= M(1) \otimes \mathbb{C}[Q]$$ 
and we set
$$ V_Qe^{ \lambda_{i}}= M(1)\otimes \mathbb{C}[Q]e^{ \lambda_{i}}, \;
\; \; i=1,\dots, n$$ Given a Lie algebra element $a \otimes t^m \in
\widehat{\goth{sl}(n+1)}$, where $a \in \goth{sl}(n+1), m \in
\mathbb{Z}$, we will denote its action on an
$\widehat{\goth{sl}(n+1)}$-module using the notation $a(m)$. In particular, for $h \in
\goth{h}$ and $m \in \mathbb{Z}$, we have the operators $h(m)$ on
$V_P$:
$$h(0)(v\otimes \iota(e_\alpha)) = \langle h,\alpha \rangle (v\otimes \iota(e_\alpha)) $$ 
$$h(m) (v\otimes \iota(e_\alpha)) =( h(m) v\otimes \iota(e_\alpha)).$$
For a formal variable $x$ and $\lambda \in P$, we define the
operator $x^\lambda$ by
$$
x^{\lambda} (v \otimes \iota(e_{\mu}))=x^{\langle \lambda, \mu \rangle} (v \otimes \iota(e_{\mu}))
$$ 
for $v \in M(1)$ and $\mu \in P$. For each $\lambda \in P$, we define the vertex operators
\begin{equation} \label{vertex}
Y(\iota(e_{\lambda}), x)=E^{-}(-\lambda, x)E^{+}(-\lambda,x)e_{\lambda}x^{\lambda},
\end{equation}
where
\begin{equation*}
E^{\pm}(-\lambda, x)=\mbox{exp} \left ( \sum_{\pm n >0}
\frac{-\lambda(n)}{n} x^{-n} \right ) \in( \mbox{End} \; V_P) [[x,
x^{-1}]]
\end{equation*}

Using the identification (\ref {identification}) we write
$Y(e^{\lambda}, x)$ for $Y(\iota(e_{\lambda}), x)$.  In particular,
for any root $\alpha \in \Delta$ we have the operators $x_{\alpha}(m)$
defined by
\begin{equation} 
Y(e^{\alpha}, x)= \sum_{m \in \mathbb{Z}} x_{\alpha}(m) x^{-m-1}.
\end{equation}
It is easy to see that 
\begin{equation}
x^{\lambda}e_{\mu}=x^{\langle \lambda, \mu \rangle}e_{\mu}x^{\lambda}
\end{equation}
and 
\begin{equation} \label{lambda-e}
\lambda(m) e_{\mu}=e_{\mu} \lambda(m)
\end{equation}
for all $\lambda, \mu \in P$ and $m \in \mathbb{Z}$. 
Using
(\ref{e-c}), (\ref{e-multiplication}) and
(\ref{vertex})-(\ref{lambda-e}) we obtain, for $\alpha \in \Delta$, $\mu \in P$,
\begin{equation}\label{almost-comm}
x_{\alpha}(m)e_{\mu}=c(\alpha, \mu) e_{\mu}x_{\alpha}(m+\langle
\alpha, \mu \rangle).
\end{equation} 

Along with the action of $\widehat{\goth{h}}$, the operators
$x_\alpha(m)$, $m \in \mathbb{Z}$, give $V_P$ a
$\widehat{\goth{sl}(n+1)}$-module structure. In particular, we have
that$$V_P = V_Q\oplus V_Qe^{\lambda_1} \oplus \dots \oplus
V_Qe^{\lambda_n}$$ and that $V_Q$, $V_Qe^{\lambda_1},\dots
,V_Qe^{\lambda_n}$ are the level $1$ basic representations of
$\widehat{{\goth{sl}(n+1)}}$ with highest weights $\Lambda_0$,
$\Lambda_1, \dots, \Lambda_n$ and highest weight vectors
$v_{\Lambda_0}=1 \otimes 1$, $v_{\Lambda_1}= 1 \otimes
e^{\lambda_1},\dots ,v_{\Lambda_n}=1 \otimes e^{\lambda_n}$,
respectively.  We make the identifications
$$L(\Lambda_0) = V_Q$$
$$L(\Lambda_i) = V_Qe^{\lambda_i}$$ for each $i=1, \dots, n$.
Moreover, taking 
$$\omega = \frac{1}{2}\sum_{i=1}^{n}
u^{(i)}(-1)^2 v_{\Lambda_0}$$ 
to be the standard conformal vector, where
$\{u^{(1)}, \dots, u^{(n)}\}$ is an orthonormal basis of $\goth{h}$,
the operators $L(m)$ defined by
\begin{equation}\label{Yomega}
Y(\omega,x)=\sum_{m \in \mathbb{Z}} L(m)x^{-m-2}
\end{equation}
provide a representation of the Virasoro algebra of central charge
$n$.  The vertex operators (\ref{vertex}) and (\ref{Yomega}) give
$L(\Lambda_0)$ the structure of a vertex operator algebra whose
irreducible modules are precisely $L(\Lambda_0)$,$L(\Lambda_1), \dots,
L(\Lambda_n)$.  We shall write
\begin{equation} \label{vectors}
v_{\Lambda_0}=\textbf{1},\  \ v_{\Lambda_1}=e^{\lambda_1}, \dots, v_{\Lambda_n}=e^{\lambda_n}.
\end{equation}

As in \cite{G1}, \cite{CLM1}--\cite{CLM2}, \cite{C3}--\cite{C4}, and \cite{CalLM1}--\cite{CalLM3},
we need certain intertwining
operators among standard modules.  We recall some facts from
\cite{FHL} and \cite{DL} about intertwining operators and, in
particular, the intertwining operators between $L(\Lambda_0)$,
$L(\Lambda_1),\dots,L(\Lambda_n)$.

Given modules $W_1$, $W_2$ and $W_3$ for
the vertex operator algebra $V$, an intertwining operator
of type
$$ \left( \begin{array} {c} W_3 \\ \begin{array}{cc} W_1 & W_2
\end{array} \end{array} \right) $$ is a linear map 
\begin{eqnarray}
\nonumber {\cal Y}(\cdotp, x) : W_1 & \longrightarrow & \mbox{Hom}
(W_2,W_3)\{x\} \nonumber \\ w &\mapsto & {\cal Y}(w,x)=\sum_{n\in
\mathbb{Q}}w_nx^{-n-1} \nonumber 
\end{eqnarray} such that all the axioms of vertex operator algebra which make sense hold (see \cite{FHL}). The main axiom is the
Jacobi identity:

   \begin{eqnarray*} \label{Jacobi-intertwining}
   \lefteqn{x_0^{-1}\delta \left ( \frac{x_1-x_2}{x_0} \right ) Y(u,
   x_1){\cal Y}(w_{(1)}, x_2)w_{(2)}} \\
   &&\hspace{2em}-x_0^{-1}\delta \left ( \frac{x_2-x_1}{x_0} \right
   ){\cal Y}(w_{(1)}, x_2)Y(u, x_1)w_{(2)} \\
   &&\displaystyle{=x_2^{-1}\delta \left ( \frac{x_1-x_0}{x_2} \right
   ) {\cal Y}(Y(u, x_0)w_{(1)}, x_2)w_{(2)}} \\ \end{eqnarray*} for
   $u \in V$, $w_{(1)} \in W_1$ and $w_{(2)} \in W_2$.

Define the operators $e^{i \pi \lambda}$ and $c(\cdotp, \lambda)$ on $V_P$ by:
 $$
e^{i \pi \lambda}(v \otimes e^{\beta})=e^{i \pi \langle \lambda ,\beta \rangle}v \otimes e^{\beta},
 $$
 $$
 c(\cdotp, \lambda)(v \otimes e^{\beta})=c(\beta, \lambda)v \otimes e^{\beta},
 $$
for $v \in M(1)$ and $ \beta , \lambda \in P$. We have that 
\begin{eqnarray} \label{IntwOp}{ \cal Y}( \cdot , x): L(\Lambda_r) &
 \longrightarrow & \mbox{Hom}(L(\Lambda_s), L( \Lambda_p))\{ x \} \\ w
 & \mapsto & {\cal Y}(w, x)=Y(w, x)e^{i\pi\lambda_r}c(\cdot,
 \lambda_r) \nonumber \end{eqnarray}
defines an intertwining operator of type 
\begin{equation} 
\label{type} \left( \begin{array} {c} L(
 \Lambda_p) \\ \begin{array}{cc} L( \Lambda_r) & L( \Lambda_s)
 \end{array} \end{array} \right)  
\end{equation}
if and only if $p \equiv r+s \;
\mbox{mod}\; (n+1)$ (cf. \cite{DL}).

If we take $u=e^{\alpha}$ and
$w_1=e^{\lambda_r}$ (for $r=1,\dots,n$) in the Jacobi identity
(\ref{Jacobi-intertwining}) and apply $\mbox{Res}_{x_0}$ (the formal residue operator,
giving us the coefficient of $x_0^{-1}$), we
have
\begin{equation} \label{CalYcommute}
[Y(e^{\alpha}, x_1), {\cal Y}(e^{\lambda_r}, x_2)]=0,
\end{equation}
whenever $\alpha \in \Delta_+$, which means that each coefficient of
the series ${\cal Y}(e^{\lambda_r}, x)$ commutes with the action of
$x_\alpha(m)$ for positive roots $\alpha$.

Given such an intertwining operator, we define a map
$${\cal Y}_c(e^{\lambda_r}, x): L(\Lambda_s) \longrightarrow L(\Lambda_p)$$ by
$${\cal Y}_c(e^{\lambda_r}, x) = \mathrm{Res}_x x^{-1-\langle \lambda_r , \lambda_s\rangle}{\cal Y}(e^{\lambda_r}, x)$$
and by (\ref{CalYcommute}) we have 
\begin{equation} \label{calYcommute}
[Y(e^{\alpha}, x_1), {\cal Y}_c(e^{\lambda_r}, x_2)]=0,
\end{equation}
which implies
\begin{equation} 
[x_\alpha(m), {\cal Y}_c(e^{\lambda_r}, x_2)]=0
\end{equation}
for each $m \in \mathbb{Z}$.

From now on through the end of Section 4, we work in the case where $n=2$.
 The finite-dimensional
simple Lie algebra $\goth{sl}(3)$ has a standard basis
$$\{ h_{\alpha_1}, \; h_{\alpha_2}, \; x_{\pm \alpha _{1}},
\; x_{\pm \alpha _{2}}, \; x_{\pm (\alpha_{1} + \alpha_{2})} \};$$
we do not need to normalize the root vectors.  We fix the
Cartan subalgebra 
$$
\goth{h}=\mathbb{C} h_{\alpha_1} \oplus
\mathbb{C}h_{\alpha_2}
$$
of ${\goth{sl}(3)}$.  Under
our identification of $\goth{h}$ with $\goth{h}^{*}$, we have
$$
\alpha_1 = h_{\alpha_1} \ \mathrm{and} \ \alpha_2 = h_{\alpha_2}.
$$ 
We also have the fundamental weights $\lambda_1, \lambda_2 \in
\goth{h}^{*}$ of $\goth{sl}(3)$, given by the
condition $\langle \lambda_i, \alpha_j \rangle = \delta_{i,j}$ for $i,
j =1,2$. In particular, we have
\begin{eqnarray*}
\lambda_1= \frac{2}{3} \alpha_1+ \frac{1}{3} \alpha_2 \; \; \mbox{and}
\; \; \lambda_2= \frac{1}{3} \alpha_1+\frac{2}{3} \alpha_2
\end{eqnarray*}
and
\begin{eqnarray*}
\alpha_1= 2\lambda_1 - \lambda_2 \; \; \mbox{and}
\; \; \alpha_2= -\lambda_1+2\lambda_2.
\end{eqnarray*}

The level $1$ standard modules of $\widehat{\goth{sl}(3)}$ are
$L(\Lambda_0)$, $L(\Lambda_1)$, and $L(\Lambda_2)$. Given
the intertwining operators (\ref{IntwOp}), we have that
\begin{eqnarray}
{\cal Y}_c(e^{\lambda_i}, x) v_{\Lambda_0} &=& r_1 v_{\Lambda_i}\\
{\cal Y}_c(e^{\lambda_i}, x) v_{\Lambda_i} &=& r_2 x_{\alpha_i}(-1)\cdot v_{\Lambda_j} = r_2'e_{\lambda_i} \cdot v_{\Lambda_i}\\
{\cal Y}_c(e^{\lambda_i}, x) v_{\Lambda_j} &=& r_3 x_{\alpha_1 + \alpha_2}(-1)\cdot v_{\Lambda_0} = r_3'e_{\lambda_i} \cdot v_{\Lambda_j}
\end{eqnarray}
for $i,j =1,2$, $i \neq j$ and some constants $r_1,r_2,r_3,r_2',r_3'
\in \mathbb{C}^{\times}$.

For any level $k$ standard $\widehat{\goth{sl}(3)}$-module
$L(\Lambda)$, its highest weight $\Lambda $ is of the form
$$\Lambda = k_0\Lambda_0 + k_1 \Lambda_1 + k_2 \Lambda_2$$ for some
nonnegative integers $k_0,k_1,k_2$ satisfying $k_0 + k_1 + k_2 =
k$. We now give a realization of these modules. Consider the space
 \begin{equation}
V_P^{\otimes k}= \underbrace{V_P \otimes \cdots \otimes V_P}_{k \; \;\mbox{times}},
\end{equation} 
and let 
\begin{equation}
v_{i_1, \dots, i_k}= v_{\Lambda_{i_1}} \otimes \cdots \otimes
v_{\Lambda_{i_k}} \in V_P^{\otimes k},
\end{equation}
where exactly $k_0$ indices are equal to $0$, $k_1$ indices are equal
to $1$ and $k_2$ indices are equal to $2$.  Then, we have that
$v_{i_1, \dots, i_k}$ is a highest weight vector for
$\widehat{\goth{sl}(3)}$, and
\begin{equation}
L(\Lambda) \simeq U(\widehat{{\goth{sl}(3)}}) \cdot v_{i_1, \dots,
i_k} \subset V_P ^{\otimes k}
\end{equation}
(cf. \cite{K}). Here, the action of $\widehat{{\goth{sl}(3)}}$ on
$V_P^{\otimes k}$ is given by the usual diagonal action of a Lie
algebra on a tensor product of modules:
\begin{equation} \label{comultiplication}
a \cdot v= \Delta(a) v= (a\otimes1 \otimes \cdots \otimes 1+ \cdots +
1 \otimes \cdots \otimes 1 \otimes a) v
\end{equation}
for $a \in \widehat{\goth{sl}(3)}, v \in V_P ^{\otimes k}$ and is
extended to $U(\widehat{{\goth{sl}(3)}})$ in the usual way. We also
have a natural vertex operator algebra structure, given by:

\begin{theorem}(\cite{FZ}, \cite{DL}, \cite{Li1}; cf. \cite{LL}) 
The standard module $L(k \Lambda_0)$ has a natural vertex operator
algebra structure. The level $k$ standard
$\widehat{{\goth{sl}(3)}}$-modules provide a complete list of
irreducible $L(k \Lambda_0)$-modules. 
\end{theorem}

The operators $L(0)$ defined in (\ref{Yomega}) provide each
$L(\Lambda)$ of level $k$ with a grading, which we refer to as the
{\em weight} grading:
\begin{equation} \label{WeightGrading}
L(\Lambda) = \coprod_{s \in \mathbb{Z}} L(\Lambda)_{(s+h_{\Lambda})}
\end{equation}
where $$h_{\Lambda} = \frac{\langle \Lambda, \Lambda + \alpha_1 +
 \alpha_2 \rangle }{2(k+3)}.$$ In particular, we have the grading
\begin{equation}
L(k\Lambda_0) = \coprod_{s \in \mathbb{Z}} L(\Lambda)_{(s)}.
\end{equation}
We denote the weight of an element $a \cdot v_{\Lambda} \in W(\Lambda)$ by 
$\mathrm{wt}(a\cdot v_{\Lambda})$. We will also write 
$$\mathrm{wt} (x_{\alpha}(m)) = -m,$$
where we view $x_\alpha(m)$ both as an operator and as an element of 
$U(\bar{\goth{n}})$.

We also have two different {\em charge} gradings on each $L(\Lambda)$
of level $k$, given by the eigenvalues of the operators $\lambda_1(0)$
and $\lambda_2(0)$:
\begin{equation} \label{ChargeGrading}
L(\Lambda) = \coprod_{r_i \in \mathbb{Z}} L(\Lambda)_{[r_i+\langle\lambda_i,\Lambda\rangle]}
\end{equation}
for each $i=1,2$. We call these the $\lambda_1$-{\em charge} and
$\lambda_2$-{\em charge} gradings, respectively. An element of
$L(\Lambda)$ with $\lambda_1$-charge $n_1$ and $\lambda_2$-charge
$n_2$ is said to have {\em total charge} $n_1+n_2$. The gradings
(\ref{WeightGrading}) and (\ref{ChargeGrading}) are compatible, and we
have that
\begin{equation} \label{TripleGrading}
L(\Lambda) = \coprod_{r_1,r_2,s \in \mathbb{Z}} L(\Lambda)_{r_1+\langle\lambda_1,\Lambda\rangle, r_2+\langle\lambda_2\Lambda\rangle;s+h_{\Lambda}}.
\end{equation}

\section{Principal subspaces, morphisms, and ideals}
\setcounter{equation}{0}

In this section, we recall the notion of principal subspace of an
$\widehat{\goth{sl}(3)}$ standard module given in
\cite{FS1}--\cite{FS2}. We then introduce some useful morphisms and
ideals, and prove various properties about these objects.

Consider the $\widehat{\goth{sl}(3)}$-subalgebras 
\begin{equation} \label{n}
\bar{\goth{n}}= \goth{n} \otimes \mathbb{C}[t, t^{-1}] ,
\end{equation}
\begin{equation}
\bar{\goth{n}}_{+}=\goth{n} \otimes \mathbb{C}[t],
\end{equation}
\begin{equation}
\bar{\goth{n}}_{-}=\goth{n} \otimes t^{-1}\mathbb{C}[t^{-1}]
\end{equation}
and the universal enveloping algebra $U(\bar{\goth{n}})$, which has the decomposition
\begin{equation}\label{decomp}
U(\bar{\goth{n}}) = U(\bar{\goth{n}}_{-}) \oplus
U(\bar{\goth{n}})\bar{\goth{n}}_{+}.
\end{equation}
For any weight $\Lambda = k_0\Lambda_0 +
k_1\Lambda_1+k_2\Lambda_2$, consider the standard module $L(\Lambda)$
with highest weight vector $v_\Lambda$. We define the {\it principal
subspace} of $L(\Lambda)$ by $$W(\Lambda) = U(\bar{\goth{n}}) \cdot
v_\Lambda \subset L(\Lambda).$$ Standard arguments show that
$W(k\Lambda_0)$ is a vertex subalgebra of $L(k\Lambda_0)$ and that
each $W(\Lambda)\subset L(\Lambda)$, $L(\Lambda)$ of level $k$, is a module for this vertex
algebra. Moreover, each $W(\Lambda)$ is preserved by $L(0)$ and by the
action of $\goth{h}$.  The gradings (\ref{WeightGrading}),
(\ref{ChargeGrading}), and (\ref{TripleGrading}) give us
gradings on $W(\Lambda)$, which are again compatible:
\begin{equation} \label{WWeightGrading}
W(\Lambda) = \coprod_{s \in \mathbb{Z}} W(\Lambda)_{(s+h_{\Lambda})},
\end{equation}
\begin{equation} \label{WChargeGrading}
W(\Lambda) = \coprod_{r_i \in \mathbb{Z}} W(\Lambda)_{[r_i+\langle\lambda_i,\Lambda\rangle]}
\end{equation}
for each $i=1,2$, and 
\begin{equation} \label{WTripleGrading}
W(\Lambda) = \coprod_{r_1,r_2,s \in \mathbb{Z}} W(\Lambda)_{r_1+\langle\lambda_1,\Lambda\rangle, r_2+\langle\lambda_2\Lambda\rangle;s+h_{\Lambda}}
\end{equation}
For each such $\Lambda$, we have a surjective map
\begin{eqnarray} \label{surj1}
F_{\Lambda}: U(\widehat{\goth{g}}) & \longrightarrow &
L(\Lambda) \\ a &\mapsto& a \cdot v_{\Lambda} \nonumber
\end{eqnarray}
and its surjective restriction $f_{\Lambda}$:
\begin{eqnarray} \label{surj2}
f_{\Lambda}: U(\bar{\goth{n}}) & \longrightarrow &
W(\Lambda)\\ a & \mapsto & a \cdot v_{\Lambda}. \nonumber
\end{eqnarray}
We will give a precise description of the kernels
$\mbox{Ker}f_{\Lambda}$ for every each $\Lambda =k_0\Lambda_0 +
k_1\Lambda_1 + k_2\Lambda_2$ and thus a presentation of the principal
subspaces $W(\Lambda)$ for $\widehat{\goth{sl}(3)}$.

As in \cite{CalLM1}--\cite{CalLM3}, define the formal sums
\begin{equation}  \label{R}
R_{t}^i= \sum_{m_1+...+m_{k+1}=-t}x_{\alpha_i}(m_1)\cdots x_{\alpha_i}(m_{k+1}) \ \ t
\in \mathbb{Z}, \ \ i=1,2;
\end{equation}
notice that $R_{t}^i$ acts in a well-defined way on any highest weight
$\widehat{\goth{g}}$-module, since, applied to any element of the
module, the formal sum becomes finite. We also define the usual
truncation of $R_t^i$ as follows:
\begin{equation} \label{r}
R_{-1, t}^i= \sum_{m_1+\dots +m_n=-t, \ m_i \le -1}x_{\alpha_i}(m_1)x_{\alpha_i}(m_2)\cdots x_{\alpha_i}(m_{k+1})
\end{equation}
for $t \in \mathbb{Z}, t \ge k+1$, and $i=1, 2$. Notice that the sums
$R_t^i$ are not actual elemnts of $U(\bar{\goth{n}})$, but that their truncations 
$R_{-1,t}^i$ are. Let $J$ be the left
ideal of $U(\bar{\goth{n}})$ generated by the elements $R_{-1, t}^i$
for $t \geq k+1$ and $i=1,2$:
\begin{equation}
J = \sum_{i=1}^2 \sum_{t \geq k+1} U(\bar{\goth{n}}) R_{-1,t}^i.
\end{equation}
As in \cite{CalLM1}--\cite{CalLM3}, we define certain left ideals of $U(\bar{\goth{n}})$ by:
\begin{equation*} 
I_{k\Lambda_{0}}= J+U(\bar{\goth{n}})\bar{\goth{n}}_{+}
\end{equation*}
and
\begin{eqnarray*} 
\lefteqn{I_{k_0\Lambda_0 + k_1\Lambda_1 + k_2\Lambda_2}}\\ &&
=I_{k\Lambda_0} + U(\bar{\goth{n}})x_{\alpha_1}(-1)^{k_0 + k_2 + 1} +
U(\bar{\goth{n}})x_{\alpha_2}(-1)^{k_0 + k_1 + 1} +
U(\bar{\goth{n}})x_{\alpha_1+\alpha_2}(-1)^{k_0 + 1}
\end{eqnarray*}
These ideals are generated by certain natural elements of $U(\bar{\goth{n}})$ which
annihilate the highest weight vectors in $W(k_0\Lambda_0 + k_1 \Lambda_1 + k_2\Lambda_2)$,
as we will see later.

\begin{remark} \em
Note that, as in  \cite{CalLM1}--\cite{CalLM3}, we have that 
$$ L(0)I_\Lambda \subset I_\Lambda$$ 
and
$$ L(0) \mathrm{Ker}f_\Lambda \subset \mathrm{Ker}f_\Lambda$$
for each $\Lambda = 
k_0 \Lambda_0 + k_1 \Lambda_1 + k_2 \Lambda_2$ with $k_0, k_1, k_2 \in \mathbb{N}$
We also have that $R_{-1,t}^1$ and $R_{-1,t}^2$ have weight $t$:
$$L(0)R_{-1,t}^1 = tR_{-1,t}^1$$
$$L(0)R_{-1,t}^2 = tR_{-1,t}^2.$$
In particular, the spaces $I_\Lambda$ and $\mathrm{Ker}f_{\Lambda}$ are $L(0)$-stable and
are graded by weight.  Moreover, 
$R_{-1,t}^1$ and $R_{-1,t}^2$ each have total charge $k+1$, and the spaces $I_\Lambda$ and $\mathrm{Ker}f_{\Lambda}$
are graded by charge. The charge and weight gradings are compatible. 
\em
\end{remark}

We recall the translation maps used in
\cite{CalLM1}--\cite{CalLM3}. For each $\lambda \in P$ and character $\nu
: Q \longrightarrow \mathbb{C}^\times$, we define a map
$\tau_{\lambda,\nu}$ on $\bar{\goth{n}}$ by
\begin{equation*}
\tau_{\lambda, \nu}(x_{\alpha}(m))=\nu(\alpha) x_{\alpha}(m-\langle
\lambda, \alpha \rangle)
\end{equation*}
for $\alpha \in \Delta_{+}$ and $m \in \mathbb{Z}$.  It is easy to see
that $\tau_{\lambda,\nu}$ is an automorphism of $\bar{\goth{n}}$.  In
the special case when $\nu$ is trivial (i.e., $\nu=1$), we set
\begin{equation*}
\tau_{\lambda}=\tau_{\lambda,1}.
\end{equation*}
The map $\tau_{\lambda,\nu}$ extends canonically to an automorphism of
$U(\bar{\goth{n}})$, also denoted by $\tau_{\lambda,\nu}$, given by
\begin{equation}\label{def-lambda_i} 
\tau_{\lambda,\nu}(x_{\beta_1}(m_1) \cdots x_{\beta_r}(m_r))=
\nu(\beta_1+\cdots + \beta_r)
x_{\beta_1}(m_1-\langle \lambda, \beta_1 \rangle ) \cdots
x_{\beta_r}(m_r-\langle \lambda, \beta_r \rangle)
\end{equation} 
for $\beta_1, \dots, \beta_r \in \Delta_{+}$ and $m_1, \dots, m_r \in
\mathbb{Z}$.  Notice that if $\lambda = \lambda_i$ for $i=1,2$, we have that
$$\mathrm{wt}(\tau_{-\lambda}(a)) \le \mathrm{wt}(a)$$
for each $a \in U(\bar{\goth{n}})$. We will use this fact frequently without mention.

We proceed to further generalize these maps. Define the
injective maps
\begin{eqnarray}
\tau_{\lambda_1, \nu}^{k_0 \Lambda_0 + k_1 \Lambda_1 + k_2 \Lambda_2}:
U(\bar{\goth{n}}) & \longrightarrow & U(\bar{\goth{n}})\\ a & \mapsto
& \tau_{\lambda_1,\nu}(a)x_{\alpha_1}(-1)^{k_1}x_{\alpha_1 +
  \alpha_2}(-1)^{k_2}. \nonumber
\end{eqnarray}
and
\begin{eqnarray}
\tau_{\lambda_2, \nu}^{k_0 \Lambda_0 + k_1 \Lambda_1 + k_2 \Lambda_2}:
U(\bar{\goth{n}}) & \longrightarrow & U(\bar{\goth{n}})\\ a & \mapsto
& \tau_{\lambda_2,\nu}(a)x_{\alpha_2}(-1)^{k_2}x_{\alpha_1 +
  \alpha_2}(-1)^{k_1}. \nonumber
\end{eqnarray}
Let $\omega_i = \alpha_i- \lambda_i \in P$ for $i=1,2$. Generalizing
the idea of \cite{CalLM3}, we define, for each character $\nu: Q \rightarrow \mathbb{C}^\times$,
injective linear maps
\begin{eqnarray}
\sigma_{\omega_1, \nu}^{k_1 \Lambda_1 + k_2 \Lambda_2}: U(\bar{\goth{n}}) & \longrightarrow &
U(\bar{\goth{n}})\\ a & \mapsto & \tau_{\omega_1,\nu}(a)x_{\alpha_1}(-1)^{k_1}. \nonumber
\end{eqnarray}
and 
\begin{eqnarray}
\sigma_{\omega_2, \nu}^{k_1 \Lambda_1 + k_2 \Lambda_2}: U(\bar{\goth{n}}) & \longrightarrow &
U(\bar{\goth{n}})\\ a & \mapsto & \tau_{\omega_2,\nu}(a)x_{\alpha_2}(-1)^{k_2}. \nonumber
\end{eqnarray}

The following facts about $U(\bar{\goth{n}})$ will be useful:

\begin{lemma} \label{rootcommute}
Given $r,k \in \mathbb{N}$ and root vectors $x_\alpha, x_\beta \in
\goth{sl}(3)$ with $\alpha, \beta, \alpha+\beta \in \Delta_+$ and
$[x_\alpha,x_\beta] = C_{\alpha,\beta}x_{\alpha+\beta}$ for some
constant $C_{\alpha, \beta} \in \mathbb{C}^{\times}$, we have
\begin{eqnarray}
\lefteqn{x_\beta(m_1) \dots x_\beta(m_r)x_\alpha(-1)^k}  \nonumber \\ &=& 
\sum_{p=0}^kx_\alpha(-1)^{k-p}
\mathop{\sum_{ j_1, \cdots , j_p=1}^r}_{  j_1 < \cdots < j_p }
C_{j_1,\dots, j_p}x_\beta(m_1)\cdots x_{\alpha+\beta}(m_{j_1}-1) \dots
x_{\alpha+\beta}(m_{j_p}-1) \cdots x_\beta(m_r) \nonumber
\end{eqnarray}
for some constants $C_{j_1,\dots j_p} \in \mathbb{C}$. The constants
$C_{j_1,\dots j_p}$ are understood to be $0$ when $p>r$.
\end{lemma}
{\em Proof:} We induct on $k \in \mathbb{N}$. For $k=1$ we have:
\begin{eqnarray}
\lefteqn{x_\beta(m_1) \cdots x_\beta(m_r)x_\alpha(-1)} \nonumber
\\ &=& x_\alpha(-1) x_\beta(m_1) \cdots x_\beta(m_r) + \sum_{j=1}^r
C_{\beta,\alpha} x_\beta(m_1) \cdots x_{\alpha+ \beta}(m_j-1) \cdots
x_\beta(m_r)
\end{eqnarray}
 and so our claim is true for $k=1$. Assume that our claim is true for
 some $k \ge 1$.  Then, we have:
\begin{eqnarray*}
\lefteqn{x_\beta(m_1) \cdots x_\beta(m_r)x_\alpha(-1)^{k+1}}  \nonumber \\ &=& 
\sum_{p=0}^kx_\alpha(-1)^{k-p}
\mathop{\sum_{ j_1, \dots , j_p=1}^r}_{  j_1 < \dots < j_p }
\bigg( C_{j_1,\dots,j_p}x_\beta(m_1)\cdots x_{\alpha+\beta}(m_{j_1}-1)\cdots \\
&&\hspace{2in}  \cdots x_{\alpha+\beta}(m_{j_p}-1) \cdots x_\beta(m_r) x_\alpha(-1)\bigg)
\\ &=&
\sum_{p=0}^kx_\alpha(-1)^{k-p+1}
\mathop{\sum_{ j_1, \dots , j_p=1}^r}_{  j_1 < \dots < j_p }
C_{j_1,\dots,j_p}x_\beta(m_1)\cdots x_{\alpha+\beta}(m_{j_1}-1) \cdots
x_{\alpha+\beta}(m_{j_p}-1) \cdots x_\beta(m_r) \nonumber \\ &+&
\sum_{p=0}^kx_\alpha(-1)^{k-p}
\mathop{\sum_{s \neq j_q, s=1}}^r_{  q=1, \dots ,p} 
\mathop{\sum_{ j_1, \dots , j_p=1}^r}_{  j_1 < \dots < j_p}
\bigg( C_{j_1,\dots,j_p}C_{\beta,\alpha} x_\beta(m_1)\cdots x_{\alpha+\beta}(m_{j_1}-1)
\dots  \\ &&
\indent \hspace{2in} \cdots x_{\alpha+\beta}(m_s-1) \cdots
x_{\alpha+\beta}(m_{j_p}-1) \cdots x_\beta(m_r) \bigg) \\ &=&
\sum_{p=0}^{k+1}x_\alpha(-1)^{k+1-p}
\mathop{\sum_{ j_1, \dots , j_p=1}^r}_{  j_1 < \dots < j_p }
C_{j_1,\dots ,j_p}'x_\beta(m_1)\cdots x_{\alpha+\beta}(m_{j_1}-1) \cdots
x_{\alpha+\beta}(m_{j_p}-1) \cdots x_\beta(m_r) 
\end{eqnarray*} 
for some constants $C_{j_1,\dots,j_p}' \in \mathbb{C}$, concluding our proof.

\begin{corollary}  \label{Rcommute}
For $0 \le m\le k$ and simple roots $\alpha_i,\alpha_j \in \Delta_+$ 
such that $\alpha_i + \alpha_j \in \Delta_+$, we have
\begin{eqnarray*}
\lefteqn{R_{-1,t}^ix_{\alpha_j}(-1)^m}  \\ &=&
x_{\alpha_j}(-1)^mR_{-1,t}^i +
r_1 x_{\alpha_j}(-1)^{m-1}[R_{-1,t+1}^i,x_{\alpha_j}(0)] + \dots\\ && +
r_m[\dots[R_{-1,t+m}^i,x_{\alpha_j}(0)],\dots ,x_{\alpha_j}(0)] +
bx_{\alpha_i + \alpha_j}(-1) + c 
\end{eqnarray*}
for some $r_1 \dots r_{m} \in \mathbb{C}, b\in U(\bar{\goth{n}}),$ and
$c \in U(\bar{\goth{n}})\bar{\goth{n}}_+$. In particular, we have that
$$ R_{-1,t}^ix_{\alpha_j}(-1)^m \in I_{k\Lambda_0} + U(\bar{\goth{n}})x_{\alpha_i + \alpha_j}(-1). $$
Moreover, if $a \in I_{k\Lambda_0}$ then 
$$ax_{\alpha_j}(-1)^m \in I_{k\Lambda_0} + U(\bar{\goth{n}})x_{\alpha_i + \alpha_j}(-1).$$
\end{corollary}

The next two lemmas show that the maps $\tau_{\lambda_i, \nu}^{k_0\Lambda_0+k_1\Lambda_1+k_2\Lambda_2}$
and $\sigma_{\omega_i, \nu}^{k_1\Lambda_1+k_2\Lambda_2}$, $i=1,2$, allow us to move between the left ideals
we have defined.

\begin{lemma} \label{taulemma}
For every character $\nu$, we have that 
$$ \tau_{\lambda_1, \nu}^{k_0 \Lambda_0 + k_1 \Lambda_1 + k_2
  \Lambda_2} ( I_{k_0\Lambda_0 + k_1 \Lambda_1 + k_2 \Lambda_2})
\subset I_{k_2 \Lambda_0 + k_0 \Lambda_1 + k_1 \Lambda_2} $$
and
$$ \tau_{\lambda_2, \nu}^{k_0 \Lambda_0 + k_1 \Lambda_1 + k_2
  \Lambda_2} ( I_{k_0\Lambda_0 + k_1 \Lambda_1 + k_2 \Lambda_2})
\subset I_{k_1 \Lambda_0 + k_2 \Lambda_1 + k_0 \Lambda_2}. $$
\end{lemma}
{\em Proof:} We prove the claim for $\tau_{\lambda_1, \nu}^{k_0
  \Lambda_0 + k_1 \Lambda_1 + k_2 \Lambda_2} $. The claim for
$\tau_{\lambda_2, \nu}^{k_0 \Lambda_0 + k_1 \Lambda_1 + k_2 \Lambda_2}
$ follows similarly.  Since $I_{k_0\Lambda_0 + k_1 \Lambda_1 + k_2
  \Lambda_2}$ is a homogeneous ideal, it suffices to prove our claim
for $\nu = 1$.

We have that 
\begin{eqnarray*}
\lefteqn{\tau_{\lambda_1}^{k_0 \Lambda_0 + k_1 \Lambda_1 +
    k_2\Lambda_2} (R_{-1,t}^1)} \\ &=& \tau_{\lambda_1}^{k_0 \Lambda_0
  + k_1 \Lambda_1 + k_2 \Lambda_2}
\bigg(\sum_{m_1+...+m_{k+1}=-t,\  m_i \le -1}x_{\alpha_1}(m_1)\cdots
x_{\alpha_1}(m_{k+1})\bigg) \\ &=& \sum_{m_1+\dots
  +m_{k+1}=-t, \ m_i \le -1}\tau_{\lambda_1}\bigg( x_{\alpha_1}(m_1)\cdots
x_{\alpha_1}(m_{k+1})\bigg)x_{\alpha_1}(-1)^{k_1}x_{\alpha_1+\alpha_2}(-1)^{k_2}\\ &=&
\sum_{m_1+\dots+m_{k+1}=-t,\  m_i \le -1}x_{\alpha_1}(m_1-1)\cdots
x_{\alpha_1}(m_{k+1}-1)x_{\alpha_1}(-1)^{k_1}x_{\alpha_1+\alpha_2}(-1)^{k_2}\\ &=&
x_{\alpha_1}(-1)^{k_1}x_{\alpha_1+\alpha_2}(-1)^{k_2}R_{-1,t+(k+1)}^1
+ ax_{\alpha_1}(-1)^{k_1+1}x_{\alpha_1+\alpha_2}(-1)^{k_2}\\ &=&
x_{\alpha_1}(-1)^{k_1}x_{\alpha_1+\alpha_2}(-1)^{k_2}R_{-1,t+(k+1)}^1
+ b[x_{\alpha_2}(0), \dots
  [x_{\alpha_2}(0),x_{\alpha_1}(-1)^{k_1+k_2+1}] \dots ]
\end{eqnarray*}
for some $a,b \in U(\bar{\goth{n}})$. Clearly 
\begin{eqnarray*}
\lefteqn{x_{\alpha_1}(-1)^{k_1}x_{\alpha_1+\alpha_2}(-1)^{k_2}(-1)R_{-1,t+(k+1)}^1}\\ &&+
s[x_{\alpha_2}(0), \dots
  [x_{\alpha_2}(0),x_{\alpha_1}(-1)^{k_1+k_2+1}] \dots ] \in I_{k_2
  \Lambda_0 + k_0 \Lambda_1 + k_1 \Lambda_2}
\end{eqnarray*}
and so $$\tau_{\lambda_1}^{k_0 \Lambda_0 + k_1 \Lambda_1 + k_2 \Lambda_2}
(R_{-1,t}^1) \in I_{k_2
  \Lambda_0 + k_0 \Lambda_1 + k_1 \Lambda_2}.$$\\
We also have
\begin{eqnarray*}
\lefteqn{\tau_{\lambda_1}^{k_0 \Lambda_0 + k_1 \Lambda_1 + k_2
    \Lambda_2} (R_{-1,t}^2)}\\ &=& \tau_{\lambda_1}^{k_0 \Lambda_0 +
  k_1 \Lambda_1 + k_2
  \Lambda_2}\bigg(\sum_{m_1+...+m_{k+1}=-t, \ m_i \le -1}x_{\alpha_2}(m_1)...x_{\alpha_2}(m_{k+1})\bigg)
\\ &=& \sum_{m_1+\dots+m_{k+1}=-t, \ m_i \le -1}\tau_{\lambda_1}\bigg(
x_{\alpha_2}(m_1)\cdots
x_{\alpha_2}(m_{k+1})\bigg)x_{\alpha_1}(-1)^{k_1}x_{\alpha_1+\alpha_2}(-1)^{k_2}\\ &=&
\sum_{m_1+\dots +m_{k+1}=-t, \ m_i \le -1}x_{\alpha_2}(m_1)\cdots
x_{\alpha_2}(m_{k+1})x_{\alpha_1}(-1)^{k_1}x_{\alpha_1+\alpha_2}(-1)^{k_2}\\ &=&
R_{-1,t}^2x_{\alpha_1}(-1)^{k_1}x_{\alpha_1+\alpha_2}(-1)^{k_2}\\ &=&
a + bx_{\alpha_1 + \alpha_2}(-1)^{k_2+1}
\end{eqnarray*}
for some $a \in I_{k\Lambda_0}$ and $b \in U(\bar{\goth{n}})$, with
the last equality following from Corollary \ref{Rcommute}. So we have
that $$\tau_{\lambda_1, \nu}^{k_0 \Lambda_0 + k_1 \Lambda_1 + k_2
  \Lambda_2} (R_{-1,t}^2) \in I_{k_2\Lambda_0 +
  k_0\Lambda_1+k_1\Lambda_2}.$$  Since $J$ is the left ideal of
$U(\bar{\goth{n}})$ generated by $R_{-1,t}^1$ and $R_{-1,t}^2$, we
have that
\begin{eqnarray*}
\tau_{\lambda_1}^{k_0 \Lambda_0 + k_1 \Lambda_1 + k_2 \Lambda_2}
(J) \subset I_{k_2\Lambda_0 + k_0\Lambda_1+k_1\Lambda_2}.
\end{eqnarray*}

We now show that $\tau_{\lambda_1}^{k_0 \Lambda_0 + k_1 \Lambda_1 +
  k_2 \Lambda_2}(U(\bar{\goth{n}})\bar{\goth{n}}_+) \subset
I_{k_2\Lambda_0 + k_0\Lambda_1+k_1\Lambda_2}.$ If $m\in \mathbb{N}$,
we have have that
\begin{eqnarray*}
\tau_{\lambda_1}^{k_0 \Lambda_0 + k_1 \Lambda_1 + k_2
  \Lambda_2}(x_{\alpha_1}(m)) &=&
x_{\alpha_1}(m-1)x_{\alpha_1}(-1)^{k_1}x_{\alpha_1+\alpha_2}(-1)^{k_2}\\
&=& x_{\alpha_1}(-1)^{k_1}x_{\alpha_1+\alpha_2}(-1)^{k_2}
x_{\alpha_1}(m-1) \in U(\bar{\goth{n}})\bar{\goth{n}}_+
\end{eqnarray*}
if $m>0$ and
\begin{eqnarray*}
\tau_{\lambda_1}^{k_0 \Lambda_0 + k_1 \Lambda_1 + k_2
  \Lambda_2}(x_{\alpha_1}(m)) &=&
x_{\alpha_1}(-1)^{k_1+1}x_{\alpha_1+\alpha_2}(-1)^{k_2}\\ &=& r
[x_{\alpha_2}(0),
  \dots,[x_{\alpha_2}(0),x_{\alpha_1}(-1)^{k_1+k_2+1}]\dots]\\ &\in &
I_{k_2\Lambda_0 + k_0 \Lambda_1 + k_1 \Lambda_2}
\end{eqnarray*}
for some $r \in \mathbb{C}$ if $m=0$. We also have
\begin{eqnarray*}
\tau_{\lambda_1}^{k_0 \Lambda_0 + k_1 \Lambda_1 + k_2
  \Lambda_2}(x_{\alpha_2}(m)) &=&
x_{\alpha_2}(m)x_{\alpha_1}(-1)^{k_1}x_{\alpha_1+\alpha_2}(-1)^{k_2}\\
&=&
x_{\alpha_1}(-1)^{k_1}x_{\alpha_1+\alpha_2}(-1)^{k_2}x_{\alpha_2}(m)\\
&& + rx_{\alpha_1}(-1)^{k_1-1}x_{\alpha_1 + \alpha_2}(-1)^{k_2}
x_{\alpha_1+\alpha_2}(m-1)
\end{eqnarray*}
for some $r \in \mathbb{C}$ and so $$\tau_{\lambda_1}^{k_0 \Lambda_0 +
  k_1 \Lambda_1 + k_2 \Lambda_2}(x_{\alpha_2}(m)) \in I_{k_2\Lambda_0
  + k_0\Lambda_1+k_1\Lambda_2}.$$ Finally, we have that, for $m \ge
0$,
\begin{eqnarray*}
\tau_{\lambda_1}^{k_0 \Lambda_0 + k_1 \Lambda_1 + k_2
  \Lambda_2}(x_{\alpha_1+\alpha_2}(m)) &=& x_{\alpha_1 +
    \alpha_2}(m-1)x_{\alpha_1}(-1)^{k_1}x_{\alpha_1+\alpha_2}(-1)^{k_2}\\
&=&x_{\alpha_1}(-1)^{k_1}x_{\alpha_1+\alpha_2}(-1)^{k_2}x_{\alpha_1+\alpha_2}(m-1)\\
&\in&  I_{k_2\Lambda_0 +
  k_0\Lambda_1+k_1\Lambda_2}.
\end{eqnarray*}
Since $U(\bar{\goth{n}})\bar{\goth{n}}_+$ is a left ideal of
$U(\bar{\goth{n}})$, we have that $$\tau_{\lambda_1, \nu}^{k_0
  \Lambda_0 + k_1 \Lambda_1 + k_2
  \Lambda_2}(U(\bar{\goth{n}})\bar{\goth{n}}_+) \subset
I_{k_2\Lambda_0 + k_0\Lambda_1+k_1\Lambda_2}$$ and so we
have $$\tau_{\lambda_1, \nu}^{k_0 \Lambda_0 + k_1 \Lambda_1 + k_2
  \Lambda_2}(I_{k\Lambda_0}) \subset I_{k_2\Lambda_0 +
  k_0\Lambda_1+k_1\Lambda_2}.$$
We now check the remaining terms.  We have
\begin{eqnarray*}
\lefteqn{\tau_{\lambda_1}^{k_0 \Lambda_0 + k_1 \Lambda_1 + k_2
  \Lambda_2}(x_{\alpha_1}(-1)^{k_0+k_2+1})}\\ &=&
x_{\alpha_1}(-2)^{k_0+k_2+1}x_{\alpha_1}(-1)^{k_1}x_{\alpha_1+\alpha_2}(-1)^{k_2}\\
&=& c x_{\alpha_1 + \alpha_2}(-1)^{k_2}R_{-1,2(k_0+k_2+1)+k_1}^1 +
a_1x_{\alpha_1}(-1)^{k_1+1}x_{\alpha_1 + \alpha_2}(-1)^{k_2} \\
&=& c x_{\alpha_1 + \alpha_2}(-1)^{k_2}R_{-1,2(k_0+k_2+1)+k_1}^1 +
a_2[x_{\alpha_2}(0), \dots [x_{\alpha_2}(0),x_{\alpha_1}(-1)^{k_1+k_2+1}]\\
&\in& I_{k_2\Lambda_0 + k_0\Lambda_1 + k_1\Lambda_2}
\end{eqnarray*}
for some $c,a_1,a_2 \in U(\bar{\goth{n}})$. So, since
$U(\bar{\goth{n}})x_{\alpha_1}(-1)^{k_0+k_2+1}$ is the left ideal of
$U(\bar{\goth{n}})$ generated by $x_{\alpha_1}(-1)^{k_0+k_2+1}$, we
have that $$\tau_{\lambda_1, \nu}^{k_0 \Lambda_0 + k_1 \Lambda_1 + k_2
  \Lambda_2}(U(\bar{\goth{n}})x_\alpha(-1)^{k_0+k_2+1}) \subset
I_{k_2\Lambda_0 + k_0 \Lambda_1 + k_1\Lambda_2}.$$ By Lemma \ref{rootcommute} we have
\begin{eqnarray*}
\lefteqn{\tau_{\lambda_1}^{k_0 \Lambda_0 + k_1 \Lambda_1 + k_2
  \Lambda_2}(x_{\alpha_2}(-1)^{k_0+k_1+1})}\\ &=&
x_{\alpha_2}(-1)^{k_0+k_1+1}x_{\alpha_1}(-1)^{k_1}x_{\alpha_1+\alpha_2}(-1)^{k_2}\\
&=&
x_{\alpha_1}(-1)^{k_1}x_{\alpha_2}(-1)^{k_0+k_1+1}x_{\alpha_1+\alpha_2}(-1)^{k_2}
\\
&&+ r_1 x_{\alpha_1}(-1)^{k_1-1}x_{\alpha_1 +
  \alpha_2}(-2)x_{\alpha_2}(-1)^{k_0+k_1}x_{\alpha_1+\alpha_2}(-1)^{k_2}\\
&&+ \dots + r_{k_1}x_{\alpha_1 +
  \alpha_2}(-2)^{k_1}x_{\alpha_2}(-1)^{k_0+1}x_{\alpha_1+\alpha_2}(-1)^{k_2}\\
&=&
r_0'x_{\alpha_1}(-1)^{k_1}[x_{\alpha_1}(0),\dots[x_{\alpha_1}(0),R_{-1,k+1}^2]
\dots ]\\
&& +r_1'x_{\alpha_1}(-1)^{k_1-1}[x_{\alpha_1}(0),\dots[x_{\alpha_1}(0),R_{-1,k+2}^2]\dots
]\\
&& + \dots
r_{k_1}'[x_{\alpha_1}(0),\dots[x_{\alpha_1}(0),R_{-1,2k_1+k_0+k_2+1}]\dots]
+ ax_{\alpha_1+\alpha_2}(-1)^{k_2+1}\\
&\in& I_{k_2\Lambda_0 + k_0\Lambda_1  k_1\Lambda_2}
\end{eqnarray*}
for some $a \in U(\bar{\goth{n}})$ and
$r_0',r_1,r_1',\dots,r_{k_1},r_{k_1}' \in \mathbb{C}$. So, since
$U(\bar{\goth{n}})x_{\alpha_2}(-1)^{k_0+k_1+1}$ is the left ideal of
$U(\bar{\goth{n}})$ generated by $x_{\alpha_2}(-1)^{k_0+k_1+1}$, we
have that $$\tau_{\lambda_1}^{k_0 \Lambda_0 + k_1 \Lambda_1 + k_2
  \Lambda_2}(U(\bar{\goth{n}})x_{\alpha_2}(-1)^{k_0+k_1+1}) \subset
I_{k_2\Lambda_0 + k_0 \Lambda_1 + k_1\Lambda_2}.$$ Finally, we have
that
\begin{eqnarray*}
\lefteqn{\tau_{\lambda_1}^{k_0 \Lambda_0 + k_1 \Lambda_1 + k_2
  \Lambda_2}(x_{\alpha_1 +\alpha_2}(-1)^{k_0+1})}\\ &=&
x_{\alpha_1 + \alpha_2}(-2)^{k_0+1}x_{\alpha_1}(-1)^{k_1}x_{\alpha_1+\alpha_2}(-1)^{k_2}\\
&=& r [x_{\alpha_2}(0), \dots
[x_{\alpha_2}(0),R_{-1,2k_0+2+k_1+k_2}^1],\dots] + ax_{\alpha_1+\alpha_2}(-1)^{k_2+1}
\end{eqnarray*}
for some $a \in U(\bar{\goth{n}})$ and some constant $r \in
\mathbb{C}$. So, since
$U(\bar{\goth{n}})x_{\alpha_1+\alpha_2}(-1)^{k_0+1}$ is the left ideal
of $U(\bar{\goth{n}})$ generated by
$x_{\alpha_1+\alpha_2}(-1)^{k_0+1}$, we have
that $$\tau_{\lambda_1}^{k_0 \Lambda_0 + k_1 \Lambda_1 + k_2
  \Lambda_2}(U(\bar{\goth{n}})x_{\alpha_1 + \alpha_2}(-1)^{k_01})
\subset I_{k_2\Lambda_0 + k_0 \Lambda_1 + k_1\Lambda_2}.$$ This
concludes our proof.

\begin{lemma}\label{sigmalemma}
For every character $\nu$, we have that
$$ \sigma_{\omega_1, \nu}^{k_1 \Lambda_1 + k_2
  \Lambda_2} ( I_{k_1 \Lambda_1 + k_2 \Lambda_2})
\subset I_{k_1 \Lambda_0 + k_2 \Lambda_1} $$
and
$$ \sigma_{\omega_2, \nu}^{k_1 \Lambda_1 + k_2
  \Lambda_2} ( I_{k_1 \Lambda_1 + k_2 \Lambda_2})
\subset I_{k_2 \Lambda_0 + k_1 \Lambda_2}. $$
\end{lemma}
{\em Proof:} We prove the claim for $ \sigma_{\omega_1, \nu}^{k_1 \Lambda_1 + k_2
  \Lambda_2}  $. The claim for $ \sigma_{\omega_2, \nu}^{k_1 \Lambda_1 + k_2
  \Lambda_2}  $ follows similarly.  Since $I_{k_1
  \Lambda_1 + k_2 \Lambda_2}$ is a homogeneous ideal, it suffices to
prove our claim for $\nu = 1$.  We have that 
\begin{eqnarray*}
\lefteqn{\sigma_{\omega_1}^{k_1 \Lambda_1 + k_2 \Lambda_2}
  (R_{-1,t}^1)}\\ &=& \sigma_{\omega_1}^{k_1 \Lambda_1 + k_2
  \Lambda_2}\bigg(\sum_{m_1+\dots +m_{k+1}=-t, \ m_i \le -1}x_{\alpha_1}(m_1)\cdots
x_{\alpha_1}(m_{k+1})\bigg) \\ &=& \sum_{m_1+\dots
  +m_{k+1}=-t, \ m_i \le -1}\sigma_{\omega_1}\bigg(x_{\alpha_1}(m_1)\cdots
x_{\alpha_1}(m_{k+1})\bigg)x_{\alpha_1}(-1)^{k_1}\\ &=& \sum_{m_1+\dots
  +m_{k+1}=-t, m_i \le -1}x_{\alpha_1}(m_1-1)\cdots
x_{\alpha_1}(m_{k+1}-1)x_{\alpha_1}(-1)^{k_1}\\ &=& R_{-1,t+(k+1)}^1
+ ax_{\alpha_1}(-1)^{k_1+1}
\end{eqnarray*}
for some $a \in U(\bar{\goth{n}})$ and so $\sigma_{\omega_1}^{k_1 \Lambda_1 + k_2 \Lambda_2}
(R_{-1,t}^1) \in I_{k_1\Lambda_0 + k_2\Lambda_1}$.
We also have, by Lemma \ref{rootcommute}, that
\begin{eqnarray*}
\lefteqn{\sigma_{\omega_1}^{k_1 \Lambda_1 + k_2 \Lambda_2}
  (R_{-1,t}^2)}\\ &=& \sigma_{\omega_1}^{k_1 \Lambda_1 + k_2
  \Lambda_2}\bigg(\sum_{m_1+\dots +m_{k+1}=-t, \ m_i \le -1}x_{\alpha_2}(m_1)\cdots
x_{\alpha_2}(m_{k+1})\bigg) \\ &=& \sum_{m_1+\dots
  +m_{k+1}=-t, \ m_i \le -1}\sigma_{\omega_1}\bigg(x_{\alpha_2}(m_1)\cdots
x_{\alpha_2}(m_{k+1})\bigg)x_{\alpha_1}(-1)^{k_1}\\ &=& \sum_{m_1+\dots
  +m_{k+1}=-t, \ m_i \le -1}x_{\alpha_2}(m_1+1) \cdots
x_{\alpha_2}(m_{k+1}+1)x_{\alpha_1}(-1)^{k_1}\\ &=&\sum_{m_1+\dots
  +m_{k+1}=-t, \ m_i \le -1} \sum_{p=0}^{k_1}x_{\alpha_1}(-1)^{k_1-p} \mathop{\sum_{
    j_1, \dots , j_p=1}^{k+1}}_{ j_1 < \dots < j_p }
\bigg(C_{j_1,\dots,j_p} x_{\alpha_2}(m_1+1)\cdots
x_{\alpha_1+\alpha_2}(m_{j_1})\cdots\\ &&\hspace{2.7in}\cdots
x_{\alpha_1+\alpha_2}(m_{j_p})\cdots
x_{\alpha_2}(m_{k+1}+1)\bigg)\\ &=&
\sum_{p=0}^{k_1}x_{\alpha_1}(-1)^{k_1-p} [\dots
  [R_{-1,t-(k+1-p)}^2,x_{\alpha_1}(0)],\dots,x_{\alpha_1}(0)] +
bx_{\alpha_2}(0)
\end{eqnarray*}
for some $b \in U(\bar{\goth{n}})$ and constants $C_{j_1,\dots,j_p}
\in \mathbb{C}$. Since $J$ is the left ideal of $U(\bar{\goth{n}})$
generated by $R_{-1,t}^1$ and $R_{-1,t}^2$, we have that
$$\sigma_{\omega_1, \nu}^{k_1 \Lambda_1 + k_2 \Lambda_2}(J) \subset
I_{k_1\Lambda_0 + k_2\Lambda_1}.$$ We now show that
$\sigma_{\omega_1}^{k_1 \Lambda_1 + k_2
  \Lambda_2}(U(\bar{\goth{n}})\bar{\goth{n}}_+) \subset
I_{k_1\Lambda_0 + k_2\Lambda_1}$.  We have
\begin{eqnarray*}
\sigma_{\omega_1}^{k_1 \Lambda_1 + k_2 \Lambda_2}(x_{\alpha_1}(m))=
x_{\alpha_1}(m-1)x_{\alpha_1}(-1)^{k_1} \in
U(\bar{\goth{n}})\bar{\goth{n}}_+ +
U(\bar{\goth{n}})x_{\alpha_1}(-1)^{k_1+1}
\end{eqnarray*}
\begin{eqnarray*}
\sigma_{\omega_1}^{k_1 \Lambda_1 + k_2
  \Lambda_2}(x_{\alpha_2}(m)) &=&
x_{\alpha_2}(m+1)x_{\alpha_1}(-1)^{k_1}\\
&=& cx_{\alpha_1}(-1)^{k_1-1}x_{\alpha_1 + \alpha_2}(m) +
x_{\alpha_1}(-1)^{k_1}x_{\alpha_2}(m+1)\\
&\in& U(\bar{\goth{n}})\bar{\goth{n}}_+
\end{eqnarray*}
and
\begin{eqnarray*}
\sigma_{\omega_1}^{k_1 \Lambda_1 + k_2
  \Lambda_2}(x_{\alpha_1 + \alpha_2}(m)) &=& x_{\alpha_1 +
  \alpha_2}(m)x_{\alpha_1}(-1)^{k_1}\\
&=& x_{\alpha_1}(-1)^{k_1}x_{\alpha_1 +
  \alpha_2}(m)\\
&\in& U(\bar{\goth{n}})\bar{\goth{n}}_+
\end{eqnarray*}
for $m \ge 0$. Since $U(\bar{\goth{n}})\bar{\goth{n}}_+$ is the left
ideal of $U(\bar{\goth{n}})$ generated by $\bar{\goth{n}}_+$, we have
that $$\sigma_{\omega_1}^{k_1 \Lambda_1 + k_2
  \Lambda_2}(U(\bar{\goth{n}})\bar{\goth{n}}_+) \subset
I_{k_1\Lambda_0 + k_2\Lambda_1}$$ and so we have $$\sigma_{\omega_1,
  \nu}^{ k_1 \Lambda_1 + k_2 \Lambda_2}(I_{k\Lambda_0}) \subset
I_{k_1\Lambda_0 + k_2\Lambda_1}.$$ We now check the remaining terms. We
have
\begin{eqnarray*}
\sigma_{\omega_1}^{k_1 \Lambda_1 + k_2
  \Lambda_2}(x_{\alpha_1}(-1)^{k_2+1}) &=&
x_{\alpha_1}(-2)^{k_2+1}x_{\alpha_1}(-1)^{k_1}\\
&=& rR_{-1,2k_2+2+k_1}^1 + ax_{\alpha_1}(-1)^{k_1 + 1}
\end{eqnarray*}
for some $a \in U(\bar{\goth{n}})$ and $r \in \mathbb{C}$, and
so $$\sigma_{\omega_1}^{k_1 \Lambda_1 + k_2
  \Lambda_2}(x_{\alpha_1}(-1)^{k_2+1}) \in I_{k_1\Lambda_0 +
  k_2\Lambda_1}.$$ We also have, by Lemma \ref{rootcommute},
\begin{eqnarray*}
\sigma_{\omega_1}^{k_1 \Lambda_1 + k_2
  \Lambda_2}(x_{\alpha_2}(-1)^{k_1+1}) &=&
x_{\alpha_2}(0)^{k_1+1}x_{\alpha_1}(-1)^{k_1}\\
&=& x_{\alpha_1}(-1)^{k_1}x_{\alpha_2}(0)^{k_1+1} + r_1
x_{\alpha_1}(-1)^{k_1-1}x_{\alpha_1 +
  \alpha_2}(-1)x_{\alpha_2}(0)^{k_1}\\
&& + \dots + r_{k_1}x_{\alpha_1 + \alpha_2}(-1)^{k_1}x_{\alpha_2}(0)
\end{eqnarray*}
for some constants $r_1,\dots , r_{k_1} \in \mathbb{C}$, and
so $$\sigma_{\omega_1}^{k_1 \Lambda_1 + k_2
  \Lambda_2}(x_{\alpha_2}(-1)^{k_1+1}) \in
U(\bar{\goth{n}})\bar{\goth{n}}_+\subset I_{k_1\Lambda_0 +
  k_2\Lambda_1}.$$ Finally,
\begin{eqnarray*}
\sigma_{\omega_1}^{k_1 \Lambda_1 + k_2
  \Lambda_2}(x_{\alpha_1 + \alpha_2}(-1)) &=& x_{\alpha_1 +
  \alpha_2}(-1)x_{\alpha_1}(-1)^{k_1}\\
&=& r[x_{\alpha_2}(0),x_{\alpha_1}(-1)^{k_1+1}]
\end{eqnarray*}
for some constant $r \in \mathbb{C}$. So we have
that $$\sigma_{\omega_1}^{k_1 \Lambda_1 + k_2 \Lambda_2}(x_{\alpha_1 +
  \alpha_2}(-1)^{k_1+1}) \in U(\bar{\goth{n}})x_{\alpha_1}(-1)^{k_1+1}
\subset I_{k_1\Lambda_0 + k_2\Lambda_1}.$$ This concludes our proof.

\begin{remark} \em
Lemmas \ref{taulemma} and \ref{sigmalemma} here directly generalize
Lemma 3.1 and Lemma 3.2 in \cite{CalLM3}, respectively. Lemma
\ref{sigmalemma} in this paper does not have an analogue for
$I_{k_0\Lambda_0 + k_1 \Lambda_1 + k_2 \Lambda_2}$, and will be the
main reason our proof of the presentations needs ideas other than
those found in \cite{CalLM1}-\cite{CalLM3}. \em
\end{remark}

\begin{remark} \em
Note that $\tau_{\lambda_i,\nu}^{k\Lambda_0} = \tau_{\lambda_i,\nu}$, so that,
as in \cite{CalLM1}-\cite{CalLM3}, we have $$ 
\tau_{\lambda_i,\nu}(I_{k\Lambda_0}) \subset I_{k\Lambda_i}.$$\em
\end{remark} \

For any $\lambda \in P$ we have the linear isomorphism 
\begin{eqnarray*}
e_{\lambda}: V_P \longrightarrow V_P.
\end{eqnarray*}
In particular, for $i,j=1,2$ with $i+j=3$ we have
\begin{eqnarray*}
e_{\lambda_i} \cdot v_{\lambda_0}= v_{\lambda_i}
\end{eqnarray*}
\begin{eqnarray*}
e_{\lambda_i} \cdot v_{\lambda_i}=\epsilon(\lambda_i,\lambda_i)x_{\alpha_i}(-1) \cdot v_{\lambda_j}
\end{eqnarray*}
\begin{eqnarray*}
e_{\lambda_i} \cdot v_{\lambda_j}=\epsilon(\lambda_i,\lambda_j)x_{\alpha_1+\alpha_2}(-1)\cdot v_{\lambda_0}
\end{eqnarray*}
Since 
\begin{eqnarray*}
e_{\lambda_i}x_{\alpha}(m)={c(\alpha, -\lambda_i)}
x_{\alpha}(m-\langle \alpha, \lambda_i \rangle)e_{\lambda_i} \; \; \;
\mbox{for} \; \; \alpha \in \Delta_{+} \; \; \mbox{and} \; \; m
\in \mathbb{Z}
\end{eqnarray*}
we have that 
\begin{equation}
e_{\lambda_i}(a \cdot v_{\lambda_0})=\tau_{\lambda_i,c_{-\lambda_i}}(a) \cdot
v_{\lambda_i}, \; \; \; a \in U(\bar{\goth{n}}).
\end{equation}
For any $\lambda \in P$, we define linear isomorphisms on $V_P^{\otimes k}$ by
\begin{eqnarray*}
e_{\lambda}^{\otimes k} = \underbrace{ e_{\lambda_1} \otimes \dots
  \otimes e_{\lambda_1}}_{k \; \;\mbox{times}}: V_P^{\otimes k}
\longrightarrow V_P^{\otimes k}.
\end{eqnarray*}
In particular, we have
\begin{eqnarray*}
\lefteqn{e_{\lambda_1}^{\otimes k}(\underbrace{v_{\Lambda_0} \otimes
    \dots \otimes v_{\Lambda_0}}_{k_0 \; \; \mbox{times}} \otimes
  \underbrace{v_{\Lambda_1} \otimes \dots \otimes v_{\Lambda_1}}_{k_1
    \; \; \mbox{times}} \otimes \underbrace{v_{\Lambda_2} \otimes
    \dots \otimes v_{\Lambda_2}}_{k_2 \; \; \mbox{times}})}\\ &=&
\epsilon(\lambda_1,\lambda_1)^{k_1}\epsilon(\lambda_1,\lambda_2)^{k_2}\frac{1}{k_1!}\frac{1}{k_2!}x_{\alpha_1}(-1)^{k_1}x_{\alpha_1
  + \alpha_2}(-1)^{k_2}\\ &&\hspace{1in} \cdot
(\underbrace{v_{\Lambda_1} \otimes \dots \otimes v_{\Lambda_1}}_{k_0
  \; \; \mbox{times}} \otimes \underbrace{v_{\Lambda_2} \otimes \dots
  \otimes v_{\Lambda_2}}_{k_1 \; \; \mbox{times}} \otimes
\underbrace{v_{\Lambda_0} \otimes \dots \otimes v_{\Lambda_0}}_{k_2 \;
  \; \mbox{times}}).\\
\end{eqnarray*}
and
\begin{eqnarray*}
\lefteqn{e_{\lambda_2}^{\otimes k}(\underbrace{v_{\Lambda_0} \otimes
    \dots \otimes v_{\Lambda_0}}_{k_0 \; \; \mbox{times}} \otimes
  \underbrace{v_{\Lambda_1} \otimes \dots \otimes v_{\Lambda_1}}_{k_1
    \; \; \mbox{times}} \otimes \underbrace{v_{\Lambda_2} \otimes
    \dots \otimes v_{\Lambda_2}}_{k_2 \; \; \mbox{times}})}\\ &=&
\epsilon(\lambda_2,\lambda_1)^{k_1}\epsilon(\lambda_2,\lambda_2)^{k_2}\frac{1}{k_1!}\frac{1}{k_2!}x_{\alpha_2}(-1)^{k_2}x_{\alpha_1
  + \alpha_2}(-1)^{k_1}\\ &&\hspace{1in} \cdot
(\underbrace{v_{\Lambda_2} \otimes \dots \otimes v_{\Lambda_2}}_{k_0
  \; \; \mbox{times}} \otimes \underbrace{v_{\Lambda_0} \otimes \dots
  \otimes v_{\Lambda_0}}_{k_1 \; \; \mbox{times}} \otimes
\underbrace{v_{\Lambda_1} \otimes \dots \otimes v_{\Lambda_1}}_{k_2 \;
  \; \mbox{times}}).\\
\end{eqnarray*}
This, along with the fact that 
\begin{eqnarray*}
e_{\lambda_i}^{\otimes k}x_{\alpha}(m)={c(\alpha, -\lambda_i)}
x_{\alpha}(m-\langle \alpha, \lambda_i \rangle)e_{\lambda_i}^{\otimes k} \; \; \;
\mbox{for} \; \; \alpha \in \Delta_{+} \; \;, i = 1,2, \mbox{and} \; \; m
\in \mathbb{Z}
\end{eqnarray*}
gives us
\begin{eqnarray*}
\lefteqn{e_{\lambda_1}^{\otimes k}(a \cdot (\underbrace{v_{\Lambda_0}
    \otimes \dots \otimes v_{\Lambda_0}}_{k_0 \; \; \mbox{times}}
  \otimes \underbrace{v_{\Lambda_1} \otimes \dots \otimes
    v_{\Lambda_1}}_{k_1 \; \; \mbox{times}} \otimes
  \underbrace{v_{\Lambda_2} \otimes \dots \otimes v_{\Lambda_2}}_{k_2
    \; \; \mbox{times}}))}\\ &=&
\epsilon(\lambda_1,\lambda_1)^{k_1}\epsilon(\lambda_1,\lambda_2)^{k_2}\frac{1}{k_1!}\frac{1}{k_2!}\tau_{\lambda_1,c_{-\lambda_1}}^{k_0\Lambda_0
  + k_1\Lambda_1 + k_2\Lambda_2}(a) \\ &&\hspace{1in} \cdot
(\underbrace{v_{\Lambda_1} \otimes \dots \otimes v_{\Lambda_1}}_{k_0
  \; \; \mbox{times}} \otimes \underbrace{v_{\Lambda_2} \otimes \dots
  \otimes v_{\Lambda_2}}_{k_1 \; \; \mbox{times}} \otimes
\underbrace{v_{\Lambda_0} \otimes \dots \otimes v_{\Lambda_0}}_{k_2 \;
  \; \mbox{times}}).\\
\end{eqnarray*}
and
\begin{eqnarray*}
\lefteqn{e_{\lambda_2}^{\otimes k}(a \cdot(\underbrace{v_{\Lambda_0}
    \otimes \dots \otimes v_{\Lambda_0}}_{k_0 \; \; \mbox{times}}
  \otimes \underbrace{v_{\Lambda_1} \otimes \dots \otimes
    v_{\Lambda_1}}_{k_1 \; \; \mbox{times}} \otimes
  \underbrace{v_{\Lambda_2} \otimes \dots \otimes v_{\Lambda_2}}_{k_2
    \; \; \mbox{times}}))}\\ &=&
\epsilon(\lambda_1,\lambda_2)^{k_1}\epsilon(\lambda_2,\lambda_2)^{k_2}\frac{1}{k_1!}\frac{1}{k_2!}\tau_{\lambda_2,c_{-\lambda_2}}^{k_0
  \Lambda_0 + k_1\Lambda_1 + k_2 \Lambda_2}(a) \\ &&\hspace{1in}
\cdot (\underbrace{v_{\Lambda_2} \otimes \dots \otimes
  v_{\Lambda_2}}_{k_0 \; \; \mbox{times}} \otimes
\underbrace{v_{\Lambda_0} \otimes \dots \otimes v_{\Lambda_0}}_{k_1 \;
  \; \mbox{times}} \otimes \underbrace{v_{\Lambda_1} \otimes \dots
  \otimes v_{\Lambda_1}}_{k_2 \; \; \mbox{times}}).\\
\end{eqnarray*}
In particular, the above gives us:
\begin{theorem}
For any $k_0,k_1,k_2 \in \mathbb{N}$ such that $k=k_0 + k_1 + k_2$, we
have injective maps
\begin{eqnarray*}
e_{\lambda_1}^{\otimes k}:W(k_0\Lambda_0 + k_1\Lambda_1 +
k_2\Lambda_2) \longrightarrow W(k_2\Lambda_0 + k_0\Lambda_1 +
k_2\Lambda_2)
\end{eqnarray*}
and
\begin{eqnarray*}
e_{\lambda_2}^{\otimes k}:W(k_0\Lambda_0 + k_1\Lambda_1 +
k_2\Lambda_2) \longrightarrow W(k_1\Lambda_0 + k_2\Lambda_1 +
k_0\Lambda_2).
\end{eqnarray*}
\end{theorem}

\begin{remark} \em Notice that the maps $e_{\lambda_i}^{\otimes k}$, $i=1,2$, cyclically permute the
  weights. \em
\end{remark} 

\begin{remark} \em
The maps $e_{\lambda_1}^{\otimes k}$ and $e_{\lambda_2}^{\otimes k}$
 were the motivation for the definitions of
$\tau_{\lambda_1,\nu}^{k_0\Lambda_0+k_1\Lambda_1+k_2\Lambda_2}$ and
$\tau_{\lambda_2,\nu}^{k_0\Lambda_0+k_1\Lambda_1+k_2\Lambda_2}$,
respectively. \em
\end{remark}

Recall that $\omega_i = \alpha_i - \lambda_i$ for $i=1,2$.
For $i,j=1,2$ with $i+j=3$ we have maps
\begin{eqnarray*}
e_{\omega_i} \cdot v_{\lambda_i}=\epsilon(\omega_i,\lambda_i)x_{\alpha_i}(-1)\cdot v_{\lambda_0}
\end{eqnarray*}
\begin{eqnarray*}
e_{\omega_i} \cdot v_{\lambda_j}=\epsilon(\omega_i,\lambda_j)v_{\lambda_i}
\end{eqnarray*}
As operators, we have that
\begin{eqnarray*}
e_{\omega_i}x_{\alpha}(m)={c(\alpha, -\omega_i)}
x_{\alpha}(m-\langle \alpha, \omega_i \rangle)e_{\omega_i} \; \; \;
\mbox{for} \; \; \alpha \in \Delta_{+} \; \; \mbox{and} \; \; m
\in \mathbb{Z}
\end{eqnarray*} 
For any such $\omega_i \in P$, we define linear isomorphisms on $V_P^{\otimes k}$ by
\begin{eqnarray*}
e_{\omega_i}^{\otimes k} = \underbrace{ e_{\omega_i} \otimes \dots
  \otimes e_{\omega_i}}_{k \; \;\mbox{times}}: V_P^{\otimes k}
\longrightarrow V_P^{\otimes k}.
\end{eqnarray*}
In particular, we have
\begin{eqnarray*}
\lefteqn{e_{\omega_1}^{\otimes k}(\underbrace{v_{\Lambda_1} \otimes
    \dots \otimes v_{\Lambda_1}}_{k_1 \; \; \mbox{times}} \otimes
  \underbrace{v_{\Lambda_2} \otimes \dots \otimes v_{\Lambda_2}}_{k_2
    \; \; \mbox{times}})}\\ &=&
\epsilon(\omega_1,\lambda_1)^{k_1}\epsilon(\omega_1,\lambda_2)^{k_2}\frac{1}{k_1!}x_{\alpha_1}(-1)^{k_1}\cdot (\underbrace{v_{\Lambda_0}
  \otimes \dots \otimes v_{\Lambda_0}}_{k_1 \; \; \mbox{times}}
\otimes \underbrace{v_{\Lambda_1} \otimes \dots \otimes
  v_{\Lambda_1}}_{k_2 \; \; \mbox{times}}).\\
\end{eqnarray*}
and
\begin{eqnarray*}
\lefteqn{e_{\omega_2}^{\otimes k}(\underbrace{v_{\Lambda_1} \otimes
    \dots \otimes v_{\Lambda_1}}_{k_1 \; \; \mbox{times}} \otimes
  \underbrace{v_{\Lambda_2} \otimes \dots \otimes v_{\Lambda_2}}_{k_2
    \; \; \mbox{times}})}\\ &=&
\epsilon(\omega_2,\lambda_1)^{k_1}\epsilon(\omega_2,\lambda_2)^{k_2}\frac{1}{k_2!}x_{\alpha_2}(-1)^{k_2}\cdot (\underbrace{v_{\Lambda_2}
  \otimes \dots \otimes v_{\Lambda_2}}_{k_1 \; \; \mbox{times}}
\otimes \underbrace{v_{\Lambda_0} \otimes \dots \otimes
  v_{\Lambda_0}}_{k_2 \; \; \mbox{times}}).\\
\end{eqnarray*}
This, along with the fact that 
\[
e_{\omega_i}^{\otimes k}x_{\alpha}(m)={c(\alpha, -\omega_i)}
x_{\alpha}(m-\langle \alpha, \omega_i \rangle)e_{\omega_i}^{\otimes k} \; \; \;
\mbox{for} \; \; \alpha \in \Delta_{+} \; \;, i = 1,2, \mbox{and} \; \; m
\in \mathbb{Z}
\] 
give us
\begin{eqnarray*}
\lefteqn{e_{\omega_1}^{\otimes k}(a \cdot (\underbrace{v_{\Lambda_1} \otimes \dots \otimes
  v_{\Lambda_1}}_{k_1 \; \; \mbox{times}} \otimes \underbrace{v_{\Lambda_2} \otimes \dots \otimes
  v_{\Lambda_2}}_{k_2 \; \; \mbox{times}}))}\\ 
&=& \epsilon(\omega_1,\lambda_1)^{k_1}\epsilon(\omega_1,\lambda_2)^{k_2}\frac{1}{k_1!}\sigma_{\omega_1,c_{-\omega_1}}^{k_1\Lambda_1 +
  k_2\Lambda_2}(a) \cdot (\underbrace{v_{\Lambda_0} \otimes \dots \otimes
  v_{\Lambda_0}}_{k_1 \; \; \mbox{times}} \otimes \underbrace{v_{\Lambda_1} \otimes \dots \otimes
  v_{\Lambda_1}}_{k_2 \; \; \mbox{times}}).\\
\end{eqnarray*}
and
\begin{eqnarray*}
\lefteqn{e_{\omega_2}^{\otimes k}(a \cdot (\underbrace{v_{\Lambda_1} \otimes \dots \otimes
  v_{\Lambda_1}}_{k_1 \; \; \mbox{times}} \otimes \underbrace{v_{\Lambda_2} \otimes \dots \otimes
  v_{\Lambda_2}}_{k_2 \; \; \mbox{times}}))}\\ 
&=& \epsilon(\omega_2,\lambda_1)^{k_1}\epsilon(\omega_2,\lambda_2)^{k_2}\frac{1}{k_2!}\sigma_{\omega_2,c_{-\omega_2}}^{k_1\Lambda_1 +
  k_2\Lambda_2}(a) \cdot (\underbrace{v_{\Lambda_2} \otimes \dots \otimes
  v_{\Lambda_2}}_{k_1 \; \; \mbox{times}} \otimes \underbrace{v_{\Lambda_0} \otimes \dots \otimes
  v_{\Lambda_0}}_{k_2 \; \; \mbox{times}}).\\
\end{eqnarray*}
In particular, the above gives us:
\begin{theorem}
For any $k_1,k_2 \in \mathbb{N}$ such that $k=k_1 + k_2$, we
have injective maps
\begin{eqnarray*}
e_{\omega_1}^{\otimes k}:W(k_1\Lambda_1 +
k_2\Lambda_2) \longrightarrow W(k_1\Lambda_0 + k_2\Lambda_1)
\end{eqnarray*}
and
\begin{eqnarray*}
e_{\omega_2}^{\otimes k}:W(k_1\Lambda_1 +
k_2\Lambda_2) \longrightarrow W(k_1\Lambda_1 + k_2\Lambda_0).
\end{eqnarray*}
\end{theorem}

\begin{remark} \em
The maps $e_{\omega_1}^{\otimes k}$ and $e_{\omega_2}^{\otimes k}$
were the motivation for the definitions of
$\sigma_{\omega_1,\nu}^{k_1\Lambda_1+k_2\Lambda_2}$ and
$\sigma_{\omega_2,\nu}^{k_1\Lambda_1+k_2\Lambda_2}$, respectively. \em
\end{remark}

\section{Proof of the main result}
\setcounter{equation}{0}

We are now ready to state and prove the main result.

\begin{theorem}
Let $k \in \mathbb{N}$. For every $k_0,k_1,k_2 \in \mathbb{N}$ such
that $k_0 + k_1 + k_2 = k$ and weight $\Lambda = k_0 \Lambda_0 +
k_1\Lambda_1 + k_2 \Lambda_2$, we have that 
\begin{eqnarray*}
\mathrm{Ker}f_\Lambda = I_\Lambda.
\end{eqnarray*}
\end{theorem}

{\em Proof:} The fact that $I_\Lambda \subset \mbox{Ker} f_\Lambda$ is
clear. Indeed, the $(k+1)$-power of each vertex operator
$Y(e^{\alpha_j},x)$, $j=1,\dots , n$, is equal to $0$ on $V_P^{\otimes k}$, and so we
have $Y(e^{\alpha_j},x)^{k+1}=0$ on each $W(\Lambda)$ of level $k$. In
particular, we have
$$Y(e^{\alpha_j},x)^{k+1} = \sum_{t \in \mathbb{Z}}R_{t}^j
x^{t-(k+1)}$$ which implies
$$\mathrm{Res}_x x^{-t+k}Y(e^{\alpha_j},x)^{k+1}\cdot v_{\Lambda}
=R_{-1,t}^j\cdot v_{\Lambda} =0,$$ and so we have that $J \subset
\mathrm{Ker}f_\Lambda$. The fact that $U(\bar{\goth{n}})\bar{\goth{n}}_+ \subset
\mathrm{Ker}f_\Lambda$ is clear. Finally, the fact that
$$x_{\alpha_1}(-1)^{k_0+k_2+1}\cdot v_\Lambda = 0$$
$$x_{\alpha_2}(-1)^{k_0+k_1+1}\cdot v_\Lambda = 0$$
$$x_{\alpha_1+\alpha_2}(-1)^{k_0+1}\cdot v_\Lambda = 0$$
follow from the fact that, in the level $1$ case, we have
$$x_{\alpha_1}(-1)\cdot v_{\Lambda_1} = 0$$
$$x_{\alpha_2}(-1)\cdot v_{\Lambda_2} = 0$$
and 
$$Y(e^\alpha,x)^2 = 0$$ on $V_P$ for each $\alpha \in \Delta_+$. Thus,
we have that $I_\Lambda \subset \mathrm{Ker}f_\Lambda$.

It remains to show that $\mbox{Ker}f_\Lambda \subset I_\Lambda$.
We proceed by contradiction.
Consider the set of elements 
\begin{eqnarray}\label{mcset}
\{ a \in U(\bar{\goth{n}}) | a \in \mbox{Ker}f_\Lambda \  \mbox{and} \ a
\notin I_\Lambda \; \;  \mbox{for some} \; \;
 \Lambda = k_0 \Lambda_0 +
k_1\Lambda_1 + k_2 \Lambda_2\}.
\end{eqnarray}
We may and do assume that homogeneous elements of (\ref{mcset}) have
positive weight (otherwise, if $a$ is such an element with non-positive weight, we have 
$a \in U(\bar{\goth{n}})\bar{\goth{n}}+
\subset I_{k\Lambda_0}$, by our decomposition (\ref{decomp})).  Among
all elements in (\ref{mcset}), we look at those of lowest total
charge.  Among all elements of lowest total charge in (\ref{mcset}),
we choose a nonzero element of lowest weight. We call this element $a$.

First, we show that $\Lambda \neq k\Lambda_1$ (that is,
 we show that $a \notin (\mbox{Ker}f_{k\Lambda_1}) \setminus I_{k\Lambda_1}$). 
Indeed, suppose that
$\Lambda = k\Lambda_1$. Then
\begin{eqnarray*}
a \cdot (v_{\Lambda_1} \otimes \dots \otimes v_{\Lambda_1}) = 0
\end{eqnarray*}
We have two cases to consider: when the
$\lambda_1$-charge of $a$ is nonzero, and when  the
$\lambda_1$-charge of $a$ is zero. If the $\lambda_1$-charge of $a$ is
nonzero, then we have:
\begin{eqnarray*} 
a \cdot (v_{\Lambda_1} \otimes \dots \otimes v_{\Lambda_1}) =
e_{\lambda_1}^{\otimes k}(\tau_{\lambda_1,c_{-\lambda_1}}^{-1}(a)\cdot (v_{\lambda_0} \otimes
  \dots \otimes v_{\lambda_0})) = 0
\end{eqnarray*}
and so, by the injectivity of $e_{\lambda_1}^{\otimes k}$, we have 
\begin{eqnarray*}
\tau_{\lambda_1,c_{-\lambda_1}}^{-1}(a)\cdot (v_{\lambda_0} \otimes
  \dots \otimes v_{\lambda_0}) = 0.
\end{eqnarray*}
Now, we have $\mathrm{wt}(\tau_{\lambda_1,c_{-\lambda_1}}^{-1}(a)) < \mathrm{wt}(a)$,
and so by assumption on $a$, we have that
$$\tau_{\lambda_1,c_{-\lambda_1}}^{-1}(a) \in I_{k\Lambda_0}.$$ But
then, by Lemma \ref{taulemma}, we have that
$\tau_{\lambda_1,c_{-\lambda_1}}(\tau_{\lambda_1,c_{-\lambda_1}}^{-1}(a))=a
\in I_{k\Lambda_1}$, a contradiction. So $a$ cannot have positive
$\lambda_1$-charge. Suppose that a has $\lambda_1$-charge equal to
$0$, and so the $\lambda_2$-charge of $a$ is positive. In this case,
we have that $\tau_{\lambda_1,c_{-\lambda_1}}^{-1}(a)$ is a nonzero constant
multiple of $a$. As before, we have
\begin{eqnarray*} 
a \cdot (v_{\Lambda_1} \otimes \dots \otimes v_{\Lambda_1}) =
e_{\lambda_1}^{\otimes k}(\tau_{\lambda_1,c_{-\lambda_1}}^{-1}(a)\cdot (v_{\lambda_0} \otimes
  \dots \otimes v_{\lambda_0}))=0
\end{eqnarray*}
which implies 
\begin{eqnarray*}
e_{\lambda_1}^{\otimes k}(a\cdot (v_{\lambda_0} \otimes
  \dots \otimes v_{\lambda_0}))=0.
\end{eqnarray*}
By the injectivity of $e_{\lambda_1}^{\otimes k}$, we have that
$$a\cdot (v_{\lambda_0} \otimes \dots \otimes v_{\lambda_0})=0.$$ Applying
the map $\mathcal{Y}_c(e^{\lambda_2},x)^{\otimes k}$, we have that
$$ a\cdot (v_{\lambda_2} \otimes\dots \otimes v_{\lambda_2})=0.$$ This gives
that 
\begin{eqnarray*}
a\cdot (v_{\lambda_2} \otimes\dots \otimes v_{\lambda_2})=e_{\lambda_2}^{\otimes k}(\tau_{\lambda_2,c_{-\lambda_2}}^{-1}(a)\cdot (v_{\lambda_0} \otimes
  \dots \otimes v_{\lambda_0}))=0.
\end{eqnarray*}
By injectivity of $e_{\lambda_2}^{\otimes k}$, we have that
$$ \tau_{\lambda_2,c_{-\lambda_2}}^{-1}(a)\cdot (v_{\lambda_0} \otimes
  \dots \otimes v_{\lambda_0})=0$$
Since $\mathrm{wt}(\tau_{\lambda_2,c_{-\lambda_2}}^{-1}(a)) < \mathrm{wt}(a)$, we have
that $\tau_{\lambda_2,c_{-\lambda_2}}^{-1}(a) \in I_{k\Lambda_0}$, and
so
$$\tau_{\lambda_2,c_{-\lambda_2}}(\tau_{\lambda_2,c_{-\lambda_2}}^{-1}(a))
= a
\in I_{k\Lambda_2} = I_{k\Lambda_0} + U(\bar{\goth{n}})x_{\alpha_2}(-1).$$ So
we may write 
\begin{eqnarray*}
a = b_1+c_1x_{\alpha_2}(-1)
\end{eqnarray*}
for some $b_1 \in I_{k\Lambda_0}$ and $c_1 \in
U(\bar{\goth{n}})$. Since $a \notin I_{k\Lambda_1}$, we have that $c_1
\neq 0$ (otherwise, $a = b_1 \in I_{k\Lambda_0} \subset
I_{k\Lambda_1}$, which is a contradiction). So we have
$$ c_1x_{\alpha_2}(-1)\cdot (v_{\lambda_0} \otimes
  \dots \otimes v_{\lambda_0}) = (a-b_1)\cdot (v_{\lambda_0} \otimes
  \dots \otimes v_{\lambda_0}) = 0 $$
Applying
$\mathcal{Y}_c(e^\lambda_1,x)\otimes\mathcal{Y}_c(e^{\lambda_2},x)^{\otimes
  k-1}$, we have 
$$ c_1x_{\alpha_2}(-1)\cdot (v_{\lambda_1} \otimes v_{\lambda_2} \otimes
  \dots \otimes v_{\lambda_2})= 0 $$
which implies
$$ c_1\cdot (e_{\lambda_2}v_{\lambda_2} \otimes v_{\lambda_2} \otimes
  \dots \otimes v_{\lambda_2}) = e_{\lambda_2}^{\otimes
    k}(\tau_{\lambda_2,c_{-\lambda_2}}^{-1}(c_1) \cdot (v_{\lambda_2} \otimes
  v_{\lambda_0} \otimes \dots \otimes v_{\lambda_0})) =0.$$
By the injectivity of $e_{\lambda_2}^{\otimes k}$, we have that 
$$\tau_{\lambda_2,c_{-\lambda_2}}^{-1}(c_1)\cdot  (v_{\lambda_2} \otimes
  v_{\lambda_0} \otimes \dots \otimes v_{\lambda_0}) = 0. $$ 
Since the total charge of $\tau_{\lambda_2,c_{-\lambda_2}}^{-1}(c_1) $
is less than the total charge of $a$, we have that
$$\tau_{\lambda_2,c_{-\lambda_2}}^{-1}(c_1) \in I_{(k-1)\Lambda_0 +
  \Lambda_2}.$$ Applying
$\tau_{\lambda_2,c_{-\lambda_2}}^{(k-1)\Lambda_0 + \Lambda_2}$, we
have 
\begin{eqnarray*}\tau_{\lambda_2,c_{-\lambda_2}}^{(k-1)\Lambda_0 +
  \Lambda_2}(\tau_{\lambda_2,c_{-\lambda_2}}^{-1}(c_1)) &=&
c_1x_{\alpha_2}(-1)\\ &\in& I_{\Lambda_1 + (k-1)\Lambda_2}
\end{eqnarray*}
by Lemma \ref{taulemma}. So we have that $$c_1x_{\alpha_2}(-1) \in I_{k\Lambda_0} +
 U(\bar{\goth{n}}) x_{\alpha_2}(-1)^2 + U(\bar{\goth{n}})x_{\alpha_1}(-1)^k+
U(\bar{\goth{n}})x_{\alpha_1+\alpha_2}(-1).$$ Since the $\lambda_1$-charge of $a$ is $0$, we may
write $$c_1x_{\alpha_2}(-1) = b' + c'x_{\alpha_2}(-1)^2$$ for some $b'
\in I_{k\Lambda_0}$ and $c' \in U(\bar{\goth{n}})$. Thus, we have
that $$a =b+b'+c'x_{\alpha_2}(-1)^2 = b_2 + c_2x_{\alpha_2}(-1)^2,$$
where we set $ b_2 = b+b' \in I_{k\Lambda_0}$ and $c_2 = c'$.
Continuing in this way, and applying the operator
$\mathcal{Y}_c(e^{\lambda_1},x)^{\otimes j} \otimes
\mathcal{Y}_c(e^{\lambda_2},x)^{\otimes (k-j)}$ at each step, we
eventually obtain that
$$ a = b_k + c_kx_{\alpha_2}(-1)^k$$ for some $b_k \in I_{k\Lambda_0}$
and $c_k \in U(\bar{\goth{n}})$. 
We have that 
$$c_kx_{\alpha_2}(-1)^k\cdot (v_{\lambda_0} \otimes \dots \otimes
v_{\lambda_0}) = (a-b)\cdot (v_{\lambda_0} \otimes \dots \otimes
v_{\lambda_0}) = 0,$$
which implies 
\begin{eqnarray*}
e_{\omega_2}^{\otimes k}e_{\lambda_2}^{\otimes k}(\tau_{\alpha_2, c_{-\alpha_2}}^{-1}(c_k)\cdot 
(v_{\lambda_0} \otimes \dots \otimes v_{\lambda_0})) =
\epsilon(\omega_2,\lambda_2)^{k}c_kx_{\alpha_2}(-1)^k\cdot (v_{\lambda_0}
\otimes \dots \otimes v_{\lambda_0}) = 0.
\end{eqnarray*}
By injectivity of both $e_{\lambda_2}^{\otimes k}$ and $e_{\omega_2}^{\otimes k}$, we have that $\tau_{\alpha_2, c_{-\alpha_2}}^{-1}(c_k)
(v_{\lambda_0} \otimes \dots \otimes v_{\lambda_0})=0$. Since
$\tau_{\alpha_2, c_{-\alpha_2}}^{-1}(c_k)$ is of lower total charge
than $a$, we have that $\tau_{\alpha_2, c_{-\alpha_2}}^{-1}(c_k) \in
I_{k\Lambda_0}$. Thus, by Lemmas \ref{taulemma} and \ref{sigmalemma},  we have that $$ (\sigma_{\omega_2,c_{-\omega_2}} \circ
\tau_{\lambda_2,c_{-\lambda_2}})(\tau_{\alpha_2,
  c_{-\alpha_2}}^{-1}(c_k)) = c_kx_{\alpha_2}(-1)^k \in
I_{k\Lambda_0}.$$ Hence $$a = b_k + c_k x_{\alpha_2}(-1)^k \in
I_{k\Lambda_0} \subset I_{k\Lambda_1},$$ a contradiction. So we have
that $\Lambda \neq k\Lambda_1$. Similarly, we have that $\Lambda \neq
k\Lambda_2$.

We now show that $\Lambda \neq k_0 \Lambda_0 + k_1 \Lambda_1$ for some
$k_0,k_1 \in \mathbb{N}$ with $k_0 + k_1 = k$. We've already shown
this for $k_0 = 0$, so we proceed by induction on $k_0$ as $k_0$
ranges from $0$ to $k$. Suppose we've shown that $\Lambda \neq k_0
\Lambda_0 + k_1 \Lambda_1$ for some $k_0,k_1 \in \mathbb{N}$ with $k_0
+ k_1 = k$ and $k_1>1$, and suppose $\Lambda = (k_0+1)\Lambda_0 +
(k_1-1)\Lambda_1$. In this case, we have that
\begin{eqnarray*}
a \cdot (\underbrace{v_{\Lambda_0} \otimes \dots \otimes
  v_{\Lambda_0}}_{k_0+1 \; \; \mbox{times}} \otimes \underbrace{v_{\Lambda_1} \otimes \dots \otimes
  v_{\Lambda_1}}_{k_1-1 \; \; \mbox{times}}) = 0.
\end{eqnarray*}
Applying the operator $1^{\otimes k_0} \otimes
\mathcal{Y}_c(e^{\lambda_1},x) \otimes 1^{\otimes k_1-1}$, we have that
\begin{eqnarray*}
a \cdot (\underbrace{v_{\Lambda_0} \otimes \dots \otimes
  v_{\Lambda_0}}_{k_0 \; \; \mbox{times}} \otimes \underbrace{v_{\Lambda_1} \otimes \dots \otimes
  v_{\Lambda_1}}_{k_1 \; \; \mbox{times}}) = 0.
\end{eqnarray*}
By our inductive hypothesis, we have that $$a \in I_{k_0\Lambda_0 + k_1
  \Lambda_1} = I_{k\Lambda_0} +
U(\bar{\goth{n}})x_{\alpha_1}(-1)^{k_0+1},$$ so we may write $$a =
b+cx_{\alpha_1}(-1)^{k_0+1}$$ for some $b \in I_{k\Lambda_0}$ and $c
\in U(\bar{\goth{n}})x_{\alpha_1}(-1)^{k_0+1}$. Here, $c \neq
0$ (otherwise, $a = b \in I_{k\Lambda_0}$, a contradiction). So we have that
\begin{eqnarray*} 
\lefteqn{cx_{\alpha_1}(-1)^{k_0+1}\cdot (\underbrace{v_{\Lambda_0} \otimes \dots \otimes
  v_{\Lambda_0}}_{k_0+1 \; \; \mbox{times}} \otimes \underbrace{v_{\Lambda_1} \otimes \dots \otimes
  v_{\Lambda_1}}_{k_1-1 \; \; \mbox{times}})}\\ &=& (a-b) \cdot (\underbrace{v_{\Lambda_0} \otimes \dots \otimes
  v_{\Lambda_0}}_{k_0+1 \; \; \mbox{times}} \otimes \underbrace{v_{\Lambda_1} \otimes \dots \otimes
  v_{\Lambda_1}}_{k_1-1 \; \; \mbox{times}})\\ &=& 0
\end{eqnarray*}
and this, using the diagonal action (\ref{comultiplication}) of $x_{\alpha_1}(-1)^{k_0+1}$, implies that
\begin{eqnarray*}
c\cdot  (\underbrace{x_{\alpha_1}(-1)v_{\Lambda_0} \otimes \dots \otimes
  x_{\alpha_1}(-1)v_{\Lambda_0}}_{k_0+1 \; \; \mbox{times}} \otimes \underbrace{v_{\Lambda_1} \otimes \dots \otimes
  v_{\Lambda_1}}_{k_1-1 \; \; \mbox{times}}) = 0.
\end{eqnarray*}
We may rewrite this as 
\begin{eqnarray*}
\lefteqn{c\cdot  (\underbrace{x_{\alpha_1}(-1)v_{\Lambda_0} \otimes \dots
    \otimes x_{\alpha_1}(-1)v_{\Lambda_0}}_{k_0+1 \; \; \mbox{times}}
  \otimes \underbrace{v_{\Lambda_1} \otimes \dots \otimes
    v_{\Lambda_1}}_{k_1-1 \; \; \mbox{times}})}\\ &=&
\epsilon(\omega_1,\lambda_1)^{-(k_0+1)}\epsilon(\omega_1,\lambda_2)^{-(k_1-1)}
e_{\omega_1}^{\otimes k}(\tau_{\omega_1,c_{-\omega_1}}^{-1}(c)\cdot 
(\underbrace{v_{\Lambda_1} \otimes \dots \otimes v_{\Lambda_1}}_{k_0+1
  \; \; \mbox{times}} \otimes \underbrace{v_{\Lambda_2} \otimes \dots
  \otimes v_{\Lambda_2}}_{k_1-1 \; \; \mbox{times}}))\\ &=& 0
\end{eqnarray*}
and so, by the injectivity of $e_{\omega_1}^{\otimes k}$, we have that
\begin{eqnarray*}
\tau_{\omega_1,c_{-\omega_1}}^{-1}(c)\cdot  (\underbrace{v_{\Lambda_1} \otimes \dots \otimes
  v_{\Lambda_1}}_{k_0+1 \; \; \mbox{times}} \otimes \underbrace{v_{\Lambda_2} \otimes \dots \otimes
  v_{\Lambda_2}}_{k_1-1 \; \; \mbox{times}}) = 0.
\end{eqnarray*}
Now, since the total charge of $\tau_{\omega_1,c_{-\omega_1}}^{-1}(c)$
is less than the total charge of $a$, we have that
$\tau_{\omega_1,c_{-\omega_1}}^{-1}(c) \in I_{(k_0+1)\Lambda_1 +
  (k_1-1)\Lambda_2}$. Using Lemma \ref{sigmalemma}, we obtain
\begin{eqnarray*}
  \sigma_{\omega_1,c_{-\omega_1}}^{(k_0+1)\Lambda_1 +
    (k_1-1)\Lambda_2}(\tau_{\omega_1,c_{-\omega_1}}^{-1}(c)) =
  cx_{\alpha_1}(-1)^{k_0+1} \in I_{(k_0+1)\Lambda_0 +
    (k_1-1)\Lambda_1} 
\end{eqnarray*}
and so 
\begin{eqnarray*}
a = b+cx_{\alpha_1}(-1)^{k_0+1} \in  I_{(k_0+1)\Lambda_0 +
    (k_1-1)\Lambda_1},
\end{eqnarray*}
a contradiction. Hence, $\Lambda \neq (k_0+1)\Lambda_0 +
(k_1-1)\Lambda_1$, completing our induction.

This gives us that $\Lambda \neq k_0\Lambda_0+k_1\Lambda_1$ for any
choice of $k_0,k_1 \in \mathbb{N}$ with $k_0 + k_1 = k$. A similar
argument shows that $\Lambda \neq k_0\Lambda_0+k_2\Lambda_2$ for any
choice of $k_0,k_2 \in \mathbb{N}$ with $k_0 + k_2 = k$.  We now
proceed to show that $\Lambda \neq k_1\Lambda_1 + k_2\Lambda_2$ for
any choice of $k_1,k_2 \in \mathbb{N}$ with $k_1 + k_2 = k$.

Suppose, for contradiction, that $\Lambda = k_1\Lambda_1 + k_2\Lambda_2$ for
some $k_1,k_2 \in \mathbb{N}$ with $k_1+k_2 = k$. We show by induction
that, given
$1 \le j \le k_2$, $a$ can
be written in the form
\begin{eqnarray*}
a = b + cx_{\alpha_1}(-1)^j
\end{eqnarray*}
for some $b \in I_{k_1\Lambda_1 + k_2\Lambda_2}$ and $c \in
U(\bar{\goth{n}})$. First, we prove the claim for $j=1$.  We have that 
\begin{eqnarray*}
a\cdot  (\underbrace{v_{\Lambda_1} \otimes \dots \otimes
  v_{\Lambda_1}}_{k_1 \; \; \mbox{times}} \otimes \underbrace{v_{\Lambda_2} \otimes \dots \otimes
  v_{\Lambda_2}}_{k_2 \; \; \mbox{times}}) = 0.
\end{eqnarray*}
Applying the operator $1^{\otimes k_1} \otimes
\mathcal{Y}_c(e^{\lambda_1},x)^{\otimes k_2}$, we have that 
\begin{eqnarray*}
a \cdot (\underbrace{v_{\Lambda_1} \otimes \dots \otimes
  v_{\Lambda_1}}_{k_1 \; \; \mbox{times}} \otimes \underbrace{e_{\lambda_1}v_{\Lambda_2} \otimes \dots \otimes
  e_{\lambda_1}v_{\Lambda_2}}_{k_2 \; \; \mbox{times}}) = 0
\end{eqnarray*}
so that 
\begin{eqnarray*}
\lefteqn{a \cdot (\underbrace{v_{\Lambda_1} \otimes \dots \otimes
  v_{\Lambda_1}}_{k_1 \; \; \mbox{times}} \otimes \underbrace{e_{\lambda_1}v_{\Lambda_2} \otimes \dots \otimes
  e_{\lambda_1}v_{\Lambda_2}}_{k_2 \; \; \mbox{times}})} \\ &=& 
e_{\lambda_1}^{\otimes k}(\tau_{\lambda_1,c_{-\lambda_1}}^{-1}(a) \cdot (\underbrace{v_{\Lambda_0} \otimes \dots \otimes
  v_{\Lambda_0}}_{k_1 \; \; \mbox{times}} \otimes \underbrace{v_{\Lambda_2} \otimes \dots \otimes
  v_{\Lambda_2}}_{k_2 \; \; \mbox{times}}))\\
&=& 0.
\end{eqnarray*}
Since $e_{\lambda_1}^{\otimes k}$ is injective, we have that 
\begin{eqnarray*}
\tau_{\lambda_1,c_{-\lambda_1}}^{-1}(a)\cdot  (\underbrace{v_{\Lambda_0} \otimes \dots \otimes
  v_{\Lambda_0}}_{k_1 \; \; \mbox{times}} \otimes \underbrace{v_{\Lambda_2} \otimes \dots \otimes
  v_{\Lambda_2}}_{k_2 \; \; \mbox{times}}) = 0.
\end{eqnarray*}
Now, $\tau_{\lambda_1,c_{-\lambda_1}}^{-1}(a)$ has the same total
charge as $a$, and satisfies
$\mathrm{wt}(\tau_{\lambda_1,c_{-\lambda_1}}^{-1}(a)) \le \mathrm{wt}(a)$. Since we've
shown that $\Lambda = k_1\Lambda_0 + k_2 \Lambda_2$ does not give the
element smallest weight among those of smallest total charge in (\ref{mcset}), we
have that
\begin{eqnarray*}
\tau_{\lambda_1,c_{-\lambda_1}}^{-1}(a)  \in
I_{k_1\Lambda_0 + k_2 \Lambda_2} = I_{k\Lambda_0} + U(\bar{\goth{n}})x_{\alpha_2}(-1)^{k_1+1}.
\end{eqnarray*}
Applying $\tau_{\lambda_1,c_{-\lambda_1}}$ and using Lemma
\ref{taulemma}, we see that 
\begin{eqnarray*}
a \in I_{k\Lambda_1} + U(\bar{\goth{n}})x_{\alpha_2}(-1)^{k_1+1} =
I_{k\Lambda_0} + U(\bar{\goth{n}})x_{\alpha_1}(-1) + U(\bar{\goth{n}})x_{\alpha_2}(-1)^{k_1+1}.
\end{eqnarray*}
So we may write
\begin{eqnarray*}
a = b + cx_{\alpha_1}(-1)
\end{eqnarray*}
for some $b \in I_{k\Lambda_0} +
U(\bar{\goth{n}})x_{\alpha_2}(-1)^{k_1+1} \subset I_{k_1\Lambda_1 +
  k_2\Lambda_2}$ and $c \in U(\bar{\goth{n}})$ and our claim holds
for $j=1$.  Now, suppose for induction that we may write
\begin{eqnarray*}
a = b+cx_{\alpha_1}(-1)^j
\end{eqnarray*}
for some $b\in I_{k_1\Lambda_1 + k_2 \Lambda_2}, c \in U(\bar{\goth{n}}), 1
\le j \le k_2-1$. We have
\begin{eqnarray*}
\lefteqn{cx_{\alpha_1}(-1)^{j}\cdot (\underbrace{v_{\Lambda_1} \otimes \dots \otimes
  v_{\Lambda_1}}_{k_1 \; \; \mbox{times}} \otimes \underbrace{v_{\Lambda_2} \otimes \dots \otimes
  v_{\Lambda_2}}_{k_2\; \; \mbox{times}})}\\ &=& (a-b)\cdot  (\underbrace{v_{\Lambda_1} \otimes \dots \otimes
  v_{\Lambda_1}}_{k_1 \; \; \mbox{times}} \otimes \underbrace{v_{\Lambda_2} \otimes \dots \otimes
  v_{\Lambda_2}}_{k_2 \; \; \mbox{times}})\\ &=& 0.
\end{eqnarray*}
Applying the operator $1^{\otimes k_1} \otimes
\mathcal{Y}_c(e^{\lambda_1},x)^{k_2-j} \otimes 1^{\otimes j}$, we have 
\begin{eqnarray*}
cx_{\alpha_1}(-1)^{j}\cdot (\underbrace{v_{\Lambda_1} \otimes \dots \otimes
  v_{\Lambda_1}}_{k_1 \; \; \mbox{times}} \otimes
\underbrace{e_{\lambda_1}v_{\Lambda_2} \otimes \dots \otimes
  e_{\lambda_1}v_{\Lambda_2}}_{k_2-j\; \; \mbox{times}} \otimes
\underbrace{v_{\Lambda_2} \otimes \dots \otimes v_{\Lambda_2}}_{j\; \;
  \mbox{times}}) = 0
\end{eqnarray*} which, using the diagonal action (\ref{comultiplication}), implies
\begin{eqnarray*}
c\cdot (\underbrace{v_{\Lambda_1} \otimes \dots \otimes v_{\Lambda_1}}_{k_1
  \; \; \mbox{times}} \otimes \underbrace{e_{\lambda_1}v_{\Lambda_2}
  \otimes \dots \otimes e_{\lambda_1}v_{\Lambda_2}}_{k_2-j\; \;
  \mbox{times}} \otimes \underbrace{e_{\lambda_1}v_{\Lambda_1} \otimes
  \dots \otimes e_{\lambda_1}v_{\Lambda_1}}_{j\; \; \mbox{times}}) = 0
\end{eqnarray*} so that
\begin{eqnarray*}
\lefteqn{c\cdot (\underbrace{v_{\Lambda_1} \otimes \dots \otimes
    v_{\Lambda_1}}_{k_1 \; \; \mbox{times}} \otimes
  \underbrace{e_{\lambda_1}v_{\Lambda_2} \otimes \dots \otimes
    e_{\lambda_1}v_{\Lambda_2}}_{k_2-j\; \; \mbox{times}} \otimes
  \underbrace{e_{\lambda_1}v_{\Lambda_1} \otimes \dots \otimes
    e_{\lambda_1}v_{\Lambda_1}}_{j \; \; \mbox{times}})}\\ &=&
e_{\lambda_1}^{\otimes
  k}(\tau_{\lambda_1,c_{-\lambda_1}}^{-1}(c)\cdot (\underbrace{v_{\Lambda_0}
  \otimes \dots \otimes v_{\Lambda_0}}_{k_1 \; \; \mbox{times}}
\otimes \underbrace{v_{\Lambda_2} \otimes \dots \otimes
  v_{\Lambda_2}}_{k_2-j\; \; \mbox{times}} \otimes
\underbrace{v_{\Lambda_1} \otimes \dots \otimes v_{\Lambda_1}}_{j\; \;
  \mbox{times}}))\\ &=&0.
\end{eqnarray*}
Since $e_{\lambda_1}^{\otimes k}$ is injective, we have
\begin{eqnarray*}
\tau_{\lambda_1,c_{-\lambda_1}}^{-1}(c)\cdot (\underbrace{v_{\Lambda_0}
  \otimes \dots \otimes v_{\Lambda_0}}_{k_1 \; \; \mbox{times}}
\otimes \underbrace{v_{\Lambda_2} \otimes \dots \otimes
  v_{\Lambda_2}}_{k_2-j\; \; \mbox{times}} \otimes
\underbrace{v_{\Lambda_1} \otimes \dots \otimes v_{\Lambda_1}}_{j\; \;
  \mbox{times}})=0.
\end{eqnarray*}
Now, since the total charge of
$\tau_{\lambda_1,c_{-\lambda_1}}^{-1}(c)$ is less than the total
charge of $a$, we have that $$\tau_{\lambda_1,c_{-\lambda_1}}^{-1}(c)
\in I_{k_1\Lambda_0 + j\Lambda_1 + (k_2-j)\Lambda_2}.$$ Recall
\begin{eqnarray*}
I_{k_1\Lambda_0 + j\Lambda_1 + (k_2-j)\Lambda_2} &=& I_{k\Lambda_0} +
U(\bar{\goth{n}})x_{\alpha_1}(-1)^{k_1+k_2-j+1} +
U(\bar{\goth{n}})x_{\alpha_2}(-1)^{k_1+j+1} \\
&& + U(\bar{\goth{n}})x_{\alpha_1 +
  \alpha_2}(-1)^{k_1 + 1}
\end{eqnarray*}
so that, applying $\tau_{\lambda_1,c_{-\lambda_1}}$, Lemma
\ref{taulemma} gives us 
\begin{eqnarray*}
c \in I_{k\Lambda_1} +
U(\bar{\goth{n}})x_{\alpha_1}(-2)^{k_1+k_2-j+1} +
U(\bar{\goth{n}})x_{\alpha_2}(-1)^{k_1+j+1} + U(\bar{\goth{n}})x_{\alpha_1 +
  \alpha_2}(-2)^{k_1 + 1}.
\end{eqnarray*}
So we may write
\begin{eqnarray*}
c = c_1 + c_2x_{\alpha_1}(-1) + c_3x_{\alpha_1}(-2)^{k_1+k_2-j+1} + c_4x_{\alpha_2}(-1)^{k_1+j+1}
+ c_5 x_{\alpha_1 + \alpha_2}(-2)^{k_1+1}
\end{eqnarray*} 
so that 
\begin{eqnarray*}
cx_{\alpha_1}(-1)^j &=& c_1x_{\alpha_1}(-1)^j +
c_2x_{\alpha_1}(-1)^{j+1} +
c_3x_{\alpha_1}(-2)^{k_1+k_2-j+1}x_{\alpha_1}(-1)^j\\ && +
c_4x_{\alpha_2}(-1)^{k_1+j+1}x_{\alpha_1}(-1)^j + c_5 x_{\alpha_1 +
  \alpha_2}(-2)^{k_1+1}x_{\alpha_1}(-1)^j
\end{eqnarray*} 
for some $c_1 \in I_{k\Lambda_0}$ and $c_2,c_3,c_4,c_5 \in
U(\bar{\goth{n}})$. We now analyze each of these terms. By Corollary
\ref{Rcommute}, we have that $c_1x_{\alpha_1}(-1)^j$ is of the form
$c_1'+ c_1''x_{\alpha_1 + \alpha_2}(-1)$ for some $c_1' \in
I_{k\Lambda_0},c_1'' \in U(\bar{\goth{n}})$, so that in particular
$c_1x_{\alpha_1}(-1)^j \in I_{k_1\Lambda_1 + k_2 \Lambda_2}$. Clearly
$c_2x_{\alpha_1}(-1)^{j+1}$ is of the desired form. For
$c_3x_{\alpha_1}(-2)^{k_1+k_2-j+1}x_{\alpha_1}(-1)^j$, we have that
\begin{eqnarray*}
c_3x_{\alpha_1}(-2)^{k_1+k_2-j+1}x_{\alpha_1}(-1)^j =
c_3'R_{-1,2(k_1+k_2-j+1)+j}^1 + c_3''x_{\alpha_1}(-1)^{j+1}
\end{eqnarray*}
for some $c_3',c_3'' \in U(\bar{\goth{n}})$ which means that
$c_3x_{\alpha_1}(-2)^{k_1+k_2-j+1}x_{\alpha_1}(-1)^j$ is of the
desired form.
For $c_4x_{\alpha_2}(-1)^{k_1+j+1}$, we use Lemma
\ref{rootcommute} to obtain 
\begin{eqnarray*}
\lefteqn{c_4x_{\alpha_2}(-1)^{k_1+j+1}x_{\alpha_1}(-1)^j} \\
&=& c_4(x_{\alpha_1}(-1)^jx_{\alpha_2}(-1)^{k_1+j+1} +
n_1x_{\alpha_1}(-1)^{j-1}x_{\alpha_1+\alpha_2}(-2)x_{\alpha_2}(-1)^{k_1+j} +\\
& & \dots + n_jx_{\alpha_1+\alpha_2}(-1)^jx_{\alpha_2}(-1)^{k_1+1})
\end{eqnarray*}
for some constants $n_1,\dots ,n_j \in \mathbb{C}$, which is clearly
an element of $I_{k_1\Lambda_1 + k_2 \Lambda_2}$. Finally, for $ c_5
x_{\alpha_1 + \alpha_2}(-2)^{k_1+1}x_{\alpha_1}(-1)^j $, again using
Lemma \ref{rootcommute}, we have that 
\begin{eqnarray*}
\lefteqn{ c_5 x_{\alpha_1 + \alpha_2}(-2)^{k_1+1}x_{\alpha_1}(-1)^j }
\\ &=& c_5(x_{\alpha_1}(-1)^{k_1+1+j}x_{\alpha_2}(-1)^{k_1+1} -
x_{\alpha_2}(-1)^{k_1+1}x_{\alpha_1}(-1)^{k_1+1+j} +\\ & & n_1
x_{\alpha_2}(-1)^{k_1}x_{\alpha_1+\alpha_2}(-2)
x_{\alpha_1}(-1)^{k_1+j} + \dots +
n_{k_1}x_{\alpha_2}(-1)x_{\alpha_1+\alpha_2}(-2)^{k_1}x_{\alpha_1}(-1)^{j+1})
\end{eqnarray*}
for some constants $n_1, \dots, n_{k_1}\in \mathbb{C}$,
which is clearly of the desired form. This completes our induction. In
particular, we have that $a$ has the form  $b +
cx_{\alpha_1}(-1)^{k_2}$ for some $b \in I_{k_1\Lambda_1 +
  k_2\Lambda_2}$ and  $c \in U(\bar{\goth{n}})$. This gives us 
\begin{eqnarray*}
\lefteqn{cx_{\alpha_1}(-1)^{k_2}\cdot (\underbrace{v_{\Lambda_1} \otimes
    \dots \otimes v_{\Lambda_1}}_{k_1 \; \; \mbox{times}} \otimes
  \underbrace{v_{\Lambda_2} \otimes \dots \otimes
    v_{\Lambda_2}}_{k_2\; \; \mbox{times}})}\\ &=& (a-b)\cdot 
(\underbrace{v_{\Lambda_1} \otimes \dots \otimes v_{\Lambda_1}}_{k_1
  \; \; \mbox{times}} \otimes \underbrace{v_{\Lambda_2} \otimes \dots
  \otimes v_{\Lambda_2}}_{k_2 \; \; \mbox{times}})\\ &=& 0.
\end{eqnarray*}
which, using the diagonal action (\ref{comultiplication}), implies that 
\begin{eqnarray*}
\lefteqn{c\cdot (\underbrace{v_{\Lambda_1} \otimes \dots \otimes
    v_{\Lambda_1}}_{k_1 \; \; \mbox{times}} \otimes
  \underbrace{e_{\lambda_1} v_{\Lambda_1} \otimes \dots \otimes
    e_{\lambda_1}v_{\Lambda_1}}_{k_2\; \; \mbox{times}})}\\ &=&
e_{\lambda_1}^{\otimes
  k}(\tau_{\lambda_1,c_{-\lambda_1}}^{-1}(c)\cdot (\underbrace{v_{\Lambda_0}
  \otimes \dots \otimes v_{\Lambda_0}}_{k_1 \; \; \mbox{times}}
\otimes \underbrace{ v_{\Lambda_1} \otimes \dots \otimes
  v_{\Lambda_1}}_{k_2\; \; \mbox{times}}))\\ &=& 0.
\end{eqnarray*}
By the injectivity of $e_{\lambda_1}^{\otimes k}$, we have that
$$\tau_{\lambda_1,c_{-\lambda_1}}^{-1}(c)\cdot (\underbrace{v_{\Lambda_0}
  \otimes \dots \otimes v_{\Lambda_0}}_{k_1 \; \; \mbox{times}}
\otimes \underbrace{ v_{\Lambda_1} \otimes \dots \otimes
  v_{\Lambda_1}}_{k_2\; \; \mbox{times}}) = 0$$ Now, since
$\tau_{\lambda_1,c_{-\lambda_1}}^{-1}(c)$ has lower total charge than
$a$, we have that $$\tau_{\lambda_1,c_{-\lambda_1}}^{-1}(c) \in
I_{k_1\Lambda_0 + k_2\Lambda_1} = I_{k\Lambda_0} +
U(\bar{\goth{n}})x_{\alpha_1}(-1)^{k_1+1}.$$ So, by Lemma
\ref{taulemma}, we have that $$c \in I_{k\Lambda_1}
+U(\bar{\goth{n}})x_{\alpha_1}(-2)^{k_1+1}.$$ We may thus write
\begin{eqnarray*}
c = c_1 + c_2x_{\alpha_1}(-1) + c_3x_{\alpha_1}(-2)^{k_1+1}
\end{eqnarray*}
for some $c_1 \in I_{k\Lambda_0}$ and $c_2,c_3 \in U(\bar{\goth{n}})$ and so
\begin{eqnarray*}
cx_{\alpha_1}(-1)^{k_2} = c_1x_{\alpha_1}(-1)^{k_2} +
c_2x_{\alpha_1}(-1)^{k_2+1} +
c_3x_{\alpha_1}(-2)^{k_1+1}x_{\alpha_1}(-1)^{k_2}.
\end{eqnarray*}
By Lemma \ref{Rcommute}, we have that $$c_1x_{\alpha_1}(-1)^{k_2} \in
I_{k\Lambda_0} + U(\bar{\goth{n}})x_{\alpha_1+\alpha_2}(-1) \subset
I_{k_1\Lambda_1 + k_2\Lambda_2}.$$ Clearly,
$c_2x_{\alpha_1}(-1)^{k_2+1} \in I_{k_1\Lambda_1 + k_2
  \Lambda_2}$. Finally,
$$c_3x_{\alpha_1}(-2)^{k_1+1}x_{\alpha_1}(-1)^{k_2} =
c_3'R_{-1,2(k_1+1)+k_2}^1 + c_3''x_{\alpha_1}(-1)^{k_2+1} \in
I_{k_1\Lambda_1+k_2\Lambda_2}$$ for some $c_3',c_3'' \in
U(\bar{\goth{n}})$. So we have that 
\begin{eqnarray*}
a = b + cx_{\alpha_1}(-1)^{k_2} \in I_{k_1\Lambda_1 + k_2\Lambda_2},
\end{eqnarray*}
which is a contradiction. Hence, we have that $\Lambda \neq
k_1\Lambda_1 + k_2\Lambda_2$ for some $k_1,k_2 \in \mathbb{N}$ such
that $k_1 + k_2 =k$. 

We now proceed to show, via induction, that $\Lambda \neq k_0\Lambda_0
+ k_1\Lambda_1 + k_2\Lambda_2$ for all $k_0,k_1,k_2 \in \mathbb{N}$
such that $k_0+k_1+k_2 = k$. We've already shown this in the case when
$k_0 = 0$. Suppose now for induction that we've shown $\Lambda \neq
k_0\Lambda_0 + k_1 \Lambda_1 + k_2\Lambda_2$ for some $0 \le k_0 < k$
and for all $k_1,k_2 \in \mathbb{N}$ such that $k_0+k_1+k_2 =
k$. Consider the weight $\Lambda = (k_0+1)\Lambda_0 + k_1\Lambda_1 +
k_2\Lambda_2$ such that $k_0+k_1+k_2+1 = k$. In this case, we have
that
\begin{eqnarray*}
a\cdot (\underbrace{v_{\Lambda_0} \otimes \dots \otimes
  v_{\Lambda_0}}_{k_0+1 \; \; \mbox{times}} \otimes
\underbrace{v_{\Lambda_1} \otimes \dots \otimes v_{\Lambda_1}}_{k_1\;
  \; \mbox{times}} \otimes \underbrace{v_{\Lambda_2} \otimes \dots
  \otimes v_{\Lambda_2}}_{k_2\; \; \mbox{times}}) = 0\\
\end{eqnarray*}
We now claim that, for any $m,n \in \mathbb{N}$ satisfying $1 \le m
\le k_0+1$ and $0 \le n \le k_0+k_2+2-m$, there exist
$b \in I_{(k_0+1)\Lambda_0 + k_1\Lambda_1 + k_2\Lambda_2}$ and $c,d \in
U(\bar{\goth{n}})$ such that
\begin{eqnarray*}
a = b + cx_{\alpha_1}(-1)^nx_{\alpha_1+\alpha_2}(-1)^m + dx_{\alpha_1+\alpha_2}(-1)^{m+1}.
\end{eqnarray*}
To show this, we use a nested induction. First, we show the claim for
$m=1$ and $n=0$, and then proceed to show it is true for $m=1$ and all
$0\le n \le k_0+k_2+1$.  

Applying the operators $1^{\otimes k_0}
\otimes \mathcal{Y}_c(e^{\lambda_1},x) \otimes 1^{\otimes (k_1 +
  k_2)}$, we obtain
\begin{eqnarray}\label{KerContain}
a\cdot (\underbrace{v_{\Lambda_0} \otimes \dots \otimes
  v_{\Lambda_0}}_{k_0 \; \; \mbox{times}} \otimes \underbrace{v_{\Lambda_1} \otimes \dots \otimes
  v_{\Lambda_1}}_{k_1+1\; \; \mbox{times}} \otimes \underbrace{v_{\Lambda_2} \otimes \dots \otimes
  v_{\Lambda_2}}_{k_2\; \; \mbox{times}}) = 0.\\ \nonumber
\end{eqnarray}
By our inductive hypothesis, we have that 
\begin{eqnarray*}a \in I_{k_0\Lambda_0 +
  (k_1+1)\Lambda_1 + k_2 \Lambda_2} &=& I_{k\Lambda_0} +
U(\bar{\goth{n}})x_{\alpha_1}(-1)^{k_0+k_2+1} +
U(\bar{\goth{n}})x_{\alpha_2}(-1)^{k_0 + k_1 + 2}\\
&& + U(\bar{\goth{n}})x_{\alpha_1 +
  \alpha_2}(-1)^{k_0+1}.
\end{eqnarray*}
So we may write
\begin{eqnarray*}
a = a_1 + a_2x_{\alpha_1}(-1)^{k_0+k_2+1} +
a_3x_{\alpha_2}(-1)^{k_0+k_1+2} + a_4x_{\alpha_1 + \alpha_2}(-1)^{k_0+1}
\end{eqnarray*}
for some $a_1 \in I_{k\Lambda_0}$ and $a_2,a_3,a_4 \in U(\bar{\goth{n}})$. 
Clearly, $$a_1, a_3x_{\alpha_2}(-1)^{k_0+k_1+2} \in I_{(k_0+1)\Lambda_0
  + k_1 \Lambda_1 + k_2 \Lambda_2}.$$  We thus have that 
\begin{eqnarray*}
\lefteqn{(a_2x_{\alpha_1}(-1)^{k_0+k_2+1} + a_4x_{\alpha_1 +
    \alpha_2}(-1)^{k_0+1})}\\ && \hspace{1.5in} \cdot
(\underbrace{v_{\Lambda_0} \otimes \dots \otimes v_{\Lambda_0}}_{k_0+1
  \; \; \mbox{times}} \otimes \underbrace{v_{\Lambda_1} \otimes \dots
  \otimes v_{\Lambda_1}}_{k_1\; \; \mbox{times}} \otimes
\underbrace{v_{\Lambda_2} \otimes \dots \otimes v_{\Lambda_2}}_{k_2\;
  \; \mbox{times}}) \\ &=& (a - a_1 - a_3x_{\alpha_2}(-1)^{k_0+k_1+2})
\\ && \hspace{1.5in} \cdot (\underbrace{v_{\Lambda_0} \otimes
  \dots \otimes v_{\Lambda_0}}_{k_0+1 \; \; \mbox{times}} \otimes
\underbrace{v_{\Lambda_1} \otimes \dots \otimes v_{\Lambda_1}}_{k_1\;
  \; \mbox{times}} \otimes \underbrace{v_{\Lambda_2} \otimes \dots
  \otimes v_{\Lambda_2}}_{k_2\; \; \mbox{times}})\\ &=& 0.
\end{eqnarray*}
Applying the operator $1^{\otimes k_0}
\otimes \mathcal{Y}_c(e^{\lambda_2},x) \otimes 1^{\otimes k_1 +
  k_2}$, we obtain 
\begin{eqnarray*}
\lefteqn{(a_2x_{\alpha_1}(-1)^{k_0+k_2+1} + a_4x_{\alpha_1 +
    \alpha_2}(-1)^{k_0+1})}\\ && \cdot (\underbrace{v_{\Lambda_0}
  \otimes \dots \otimes v_{\Lambda_0}}_{k_0 \; \; \mbox{times}}
\otimes v_{\Lambda_2}\otimes \underbrace{v_{\Lambda_1} \otimes \dots
  \otimes v_{\Lambda_1}}_{k_1\; \; \mbox{times}} \otimes
\underbrace{v_{\Lambda_2} \otimes \dots \otimes v_{\Lambda_2}}_{k_2\;
  \; \mbox{times}})\\ &=& a_2x_{\alpha_1}(-1)^{k_0+k_2+1}\cdot 
(\underbrace{v_{\Lambda_0} \otimes \dots \otimes v_{\Lambda_0}}_{k_0
  \; \; \mbox{times}} \otimes v_{\Lambda_2}\otimes
\underbrace{v_{\Lambda_1} \otimes \dots \otimes v_{\Lambda_1}}_{k_1\;
  \; \mbox{times}} \otimes \underbrace{v_{\Lambda_2} \otimes \dots
  \otimes v_{\Lambda_2}}_{k_2\; \; \mbox{times}}) \\ &=& 0,
\end{eqnarray*}
since $a_4x_{\alpha_1 +
  \alpha_2}(-1)^{k_0+1} \in I_{k_0\Lambda_0 + k_1\Lambda_1 +
  (k_2+1)\Lambda_2}$.
In particular, we may write
\begin{eqnarray*}
a_2x_{\alpha_1}(-1)^{k_0+k_2+1} \cdot (\underbrace{v_{\Lambda_0} \otimes
  \dots \otimes v_{\Lambda_0}}_{k_0 \; \; \mbox{times}} \otimes
\underbrace{v_{\Lambda_1} \otimes \dots \otimes v_{\Lambda_1}}_{k_1\;
  \; \mbox{times}} \otimes \underbrace{v_{\Lambda_2} \otimes \dots
  \otimes v_{\Lambda_2}}_{k_2+1\; \; \mbox{times}}) = 0.
\end{eqnarray*}
Applying the operator $\mathcal{Y}_c(e^{\lambda_2},x)^{\otimes k_0}
\otimes 1^{\otimes k_1+k_2+1}$, we have that 
\begin{eqnarray*}
a_2x_{\alpha_1}(-1)^{k_0+k_2+1} \cdot (\underbrace{v_{\Lambda_2} \otimes
  \dots \otimes v_{\Lambda_2}}_{k_0 \; \; \mbox{times}} \otimes
\underbrace{v_{\Lambda_1} \otimes \dots \otimes v_{\Lambda_1}}_{k_1\;
  \; \mbox{times}} \otimes \underbrace{v_{\Lambda_2} \otimes \dots
  \otimes v_{\Lambda_2}}_{k_2+1\; \; \mbox{times}}) = 0
\end{eqnarray*}
which, using the diagonal action (\ref{comultiplication}), implies
\begin{eqnarray*}
\lefteqn{a_2\cdot (\underbrace{ x_{\alpha_1}(-1)v_{\Lambda_2} \otimes \dots
    \otimes x_{\alpha_1}(-1)v_{\Lambda_2}}_{k_0 \; \;
    \mbox{times}}\otimes}\\ && \hspace{0.5in}
\underbrace{v_{\Lambda_1} \otimes \dots \otimes v_{\Lambda_1}}_{k_1\;
  \; \mbox{times}} \otimes \underbrace{ x_{\alpha_1}(-1)v_{\Lambda_2}
  \otimes \dots \otimes x_{\alpha_1}(-1)v_{\Lambda_2}}_{k_2+1\; \;
  \mbox{times}}) = 0.
\end{eqnarray*}
From this, we have that 
\begin{eqnarray*}
a_2 \cdot (\underbrace{e_{\lambda_1}v_{\Lambda_1} \otimes \dots \otimes
  e_{\lambda_1}v_{\Lambda_1}}_{k_0 \; \; \mbox{times}} \otimes
\underbrace{v_{\Lambda_1} \otimes \dots \otimes v_{\Lambda_1}}_{k_1\;
  \; \mbox{times}} \otimes \underbrace{e_{\lambda_1}v_{\Lambda_1}
  \otimes \dots \otimes e_{\lambda_1}v_{\Lambda_1}}_{k_2+1\; \;
  \mbox{times}}) = 0
\end{eqnarray*}
which implies that 
\begin{eqnarray*}
e_{\lambda_1}^{\otimes k} (\tau_{\lambda_1,c_{-\lambda_1}}^{-1}(a_2)\cdot 
(\underbrace{v_{\Lambda_1} \otimes \dots \otimes v_{\Lambda_1}}_{k_0
  \; \; \mbox{times}} \otimes \underbrace{v_{\Lambda_0} \otimes \dots
  \otimes v_{\Lambda_0}}_{k_1\; \; \mbox{times}} \otimes
\underbrace{v_{\Lambda_1} \otimes \dots \otimes
  v_{\Lambda_1}}_{k_2+1\; \; \mbox{times}})) = 0.
\end{eqnarray*}
Now, by injectivity of $e_{\lambda_1}^{\otimes k}$, we have that
\begin{eqnarray*}
\tau_{\lambda_1,c_{-\lambda_1}}^{-1}(a_2)\cdot  (\underbrace{v_{\Lambda_1}
  \otimes \dots \otimes v_{\Lambda_1}}_{k_0 \; \; \mbox{times}}
\otimes \underbrace{v_{\Lambda_0} \otimes \dots \otimes
  v_{\Lambda_0}}_{k_1\; \; \mbox{times}} \otimes
\underbrace{v_{\Lambda_1} \otimes \dots \otimes
  v_{\Lambda_1}}_{k_2+1\; \; \mbox{times}}) = 0.
\end{eqnarray*}
Since the total charge of $\tau_{\lambda_1,c_{-\lambda_1}}^{-1}(a_2)$
is less than the total charge of $a$, we have that
$$\tau_{\lambda_1,c_{-\lambda_1}}^{-1}(a_2) \in I_{k_1\Lambda_0 +
  (k_0+k_2+1)\Lambda_1} = I_{k\Lambda_0} + U(\bar{\goth{n}})
x_{\alpha_1}(-1)^{k_1+1}.$$ So, in particular,
applying $\tau_{\lambda_1,c_{-\lambda_1}}$, we obtain 
\begin{eqnarray*}
a_2 \in  I_{k\Lambda_1} + U(\bar{\goth{n}})
x_{\alpha_1}(-2)^{k_1+1}.
\end{eqnarray*}
So we may write $a_2 = a_{2,1} + a_{2,2}x_{\alpha_1}(-1) +
a_{2,3}x_{\alpha_1}(-2)^{k_1+1}$ for some $a_{2,1} \in I_{k\Lambda_0}$ and
$a_{2,2},a_{2,3} \in U(\bar{\goth{n}})$. So we have that 
\begin{eqnarray*}
\lefteqn{a_2x_{\alpha_1}(-1)^{k_0+k_2+1}}\\ &=& a_{2,1}
x_{\alpha_1}(-1)^{k_0+k_2+1}+ a_{2,2}x_{\alpha_1}(-1)^{k_0+k_2+2} +
a_{2,3}x_{\alpha_1}(-2)^{k_1+1}x_{\alpha_1}(-1)^{k_0+k_2+1}.
\end{eqnarray*}
By Lemma \ref{Rcommute}, we have that $$a_{2,1}
x_{\alpha_1}(-1)^{k_0+k_2+1} \in I_{k\Lambda_0} +
U(\bar{\goth{n}})x_{\alpha_1 + \alpha_2}(-1).$$  Clearly,
$$a_{2,2}x_{\alpha_1}(-1)^{k_0+k_2+2} \in I_{(k_0+1)\Lambda_0 +
  k_1\Lambda_1 + k_2\Lambda_2}.$$ Finally,
\begin{eqnarray*}
a_{2,3}x_{\alpha_1}(-2)^{k_1+1}x_{\alpha_1}(-1)^{k_0+k_2+1} &=&
rR_{-1,2(k_1+1) + k_0+k_2+1} + a_{2,3}'x_{\alpha_1}(-1)^{k_0 + k_2 +1}\\
&\in& I_{(k_0+1)\Lambda_0 + k_1\Lambda_1 + k_2\Lambda_2}
\end{eqnarray*}
for some constant $r \in \mathbb{C}$. Thus, we have that
$$a_2x_{\alpha_1}(-1)^{k_0+k_2+1} \in I_{(k_0+1)\Lambda_0 +
  k_1\Lambda_1 + k_2\Lambda_2} + U(\bar{\goth{n}})x_{\alpha_1 +
  \alpha_2}(-1).$$ Clearly
\begin{eqnarray*}
a_1, a_3x_{\alpha_1}(-1)^{k_0+k_2 +1},
a_4x_{\alpha_1+\alpha_2}(-1)^{k_0+1} \in I_{(k_0+1)\Lambda_0 +
  k_1\Lambda_1 + k_2\Lambda_2} + U(\bar{\goth{n}})x_{\alpha_1 +
  \alpha_2}(-1)
\end{eqnarray*}
and so we may write
\begin{eqnarray*}
a = b+cx_{\alpha_1+\alpha_2}(-1)
\end{eqnarray*}
for some $b \in I_{(k_0+1)\Lambda_0 + k_1\Lambda_1 + k_2\Lambda_2}$
and $c \in U(\bar{\goth{n}})$, and so our claim holds for $m=1,n=0$
(notice here the $d$ term is $0$).

Now, we assume our claim holds for $m=1$ and for some $n\in \mathbb{N}$ satisfying $0 \le n \le
k_0+k_2$. We show that it holds for $n+1$.  Suppose
\begin{eqnarray*}
a = b + cx_{\alpha_1}(-1)^nx_{\alpha_1+\alpha_2}(-1) + dx_{\alpha_1+\alpha_2}(-1)^2
\end{eqnarray*}
for some $b \in I_{(k_0+1)\Lambda_0 + k_1\Lambda_1 + k_2\Lambda_2}$
and $c,d \in U(\bar{\goth{n}})$.  We must now consider the cases when $0\le n <
k_2$ and $k_2 \le n \le k_0+k_2$ separately. First assume $0\le n <
k_2$. In this case, we have 
\begin{eqnarray*}
a\cdot (\underbrace{v_{\Lambda_0} \otimes \dots \otimes
  v_{\Lambda_0}}_{k_0+1 \; \; \mbox{times}} \otimes \underbrace{v_{\Lambda_1} \otimes \dots \otimes
  v_{\Lambda_1}}_{k_1\; \; \mbox{times}} \otimes \underbrace{v_{\Lambda_2} \otimes \dots \otimes
  v_{\Lambda_2}}_{k_2\; \; \mbox{times}}) =0\\
\end{eqnarray*}
which implies
\begin{eqnarray*}
\lefteqn{(cx_{\alpha_1}(-1)^nx_{\alpha_1+\alpha_2}(-1) + dx_{\alpha_1+\alpha_2}(-1)^2 )} \\
&& \hspace{1in} \cdot (\underbrace{v_{\Lambda_0} \otimes \dots \otimes
  v_{\Lambda_0}}_{k_0+1 \; \; \mbox{times}} \otimes \underbrace{v_{\Lambda_1} \otimes \dots \otimes
  v_{\Lambda_1}}_{k_1\; \; \mbox{times}} \otimes \underbrace{v_{\Lambda_2} \otimes \dots \otimes
  v_{\Lambda_2}}_{k_2\; \; \mbox{times}})\\
&=& (a-b)\cdot  (\underbrace{v_{\Lambda_0} \otimes \dots \otimes
  v_{\Lambda_0}}_{k_0+1 \; \; \mbox{times}} \otimes \underbrace{v_{\Lambda_1} \otimes \dots \otimes
  v_{\Lambda_1}}_{k_1\; \; \mbox{times}} \otimes \underbrace{v_{\Lambda_2} \otimes \dots \otimes
  v_{\Lambda_2}}_{k_2\; \; \mbox{times}})\\
&=& 0.
\end{eqnarray*}
Applying the operator $1 \otimes
\mathcal{Y}_c(e^{\lambda_1},x)^{\otimes k_0} \otimes 1^{\otimes k_1}
\otimes \mathcal{Y}_c(e^{\lambda_1},x)^{\otimes k_2-n} \otimes
1^{\otimes n}$, we obtain
\begin{eqnarray*} 
\lefteqn{(cx_{\alpha_1}(-1)^nx_{\alpha_1+\alpha_2}(-1) +
  dx_{\alpha_1+\alpha_2}(-1)^2 )} \\ && \cdot (v_{\Lambda_0} \otimes
\underbrace{v_{\Lambda_1} \otimes \dots \otimes
  v_{\Lambda_1}}_{k_0+k_1 \; \; \mbox{times}} \otimes
\underbrace{e_{\lambda_1}v_{\Lambda_2} \otimes \dots \otimes
  e_{\lambda_1}v_{\Lambda_2}}_{k_2-n\; \; \mbox{times}} \otimes
\underbrace{v_{\Lambda_2} \otimes \dots \otimes v_{\Lambda_2}}_{n\; \;
  \mbox{times}}) \\ &=&
cx_{\alpha_1}(-1)^nx_{\alpha_1+\alpha_2}(-1)\\ && \cdot (v_{\Lambda_0}
\otimes \underbrace{v_{\Lambda_1} \otimes \dots \otimes
  v_{\Lambda_1}}_{k_0+k_1 \; \; \mbox{times}} \otimes
\underbrace{e_{\lambda_1}v_{\Lambda_2} \otimes \dots \otimes
  e_{\lambda_1}v_{\Lambda_2}}_{k_2-n\; \; \mbox{times}} \otimes
\underbrace{v_{\Lambda_2} \otimes \dots \otimes v_{\Lambda_2}}_{n\; \;
  \mbox{times}})\\ &=& 0.
\end{eqnarray*}
This implies that 
\begin{eqnarray*}
c \cdot (e_{\lambda_1}v_{\Lambda_2} \otimes \underbrace{v_{\Lambda_1}
  \otimes \dots \otimes v_{\Lambda_1}}_{k_0+k_1 \; \; \mbox{times}}
\otimes \underbrace{e_{\lambda_1}v_{\Lambda_2} \otimes \dots \otimes
  e_{\lambda_1}v_{\Lambda_2}}_{k_2-n\; \; \mbox{times}} \otimes
\underbrace{e_{\lambda_1}v_{\Lambda_1} \otimes \dots \otimes
  e_{\lambda_1}v_{\Lambda_1}}_{n\; \; \mbox{times}}) = 0
\end{eqnarray*}
which gives us 
\begin{eqnarray*}
e_{\lambda_1}^{\otimes k} (\tau_{\lambda_1,c_{-\lambda_1}}^{-1}(c)\cdot 
(v_{\Lambda_2} \otimes \underbrace{v_{\Lambda_0} \otimes \dots \otimes
  v_{\Lambda_0}}_{k_0+k_1 \; \; \mbox{times}} \otimes
\underbrace{v_{\Lambda_2} \otimes \dots \otimes
  v_{\Lambda_2}}_{k_2-n\; \; \mbox{times}} \otimes
\underbrace{v_{\Lambda_1} \otimes \dots \otimes v_{\Lambda_1}}_{n\; \;
  \mbox{times}})) = 0.
\end{eqnarray*}
By the injectivity of $e_{\lambda_1}^{\otimes k}$, we have that 
\begin{eqnarray*}
\tau_{\lambda_1,c_{-\lambda_1}}^{-1}(c) \cdot (v_{\Lambda_2} \otimes
\underbrace{v_{\Lambda_0} \otimes \dots \otimes
  v_{\Lambda_0}}_{k_0+k_1 \; \; \mbox{times}} \otimes
\underbrace{v_{\Lambda_2} \otimes \dots \otimes
  v_{\Lambda_2}}_{k_2-n\; \; \mbox{times}} \otimes
\underbrace{v_{\Lambda_1} \otimes \dots \otimes v_{\Lambda_1}}_{n\; \;
  \mbox{times}})= 0.
\end{eqnarray*}
Since $\tau_{\lambda_1,c_{-\lambda_1}}^{-1}(c)$ is of lower total
charge than $a$, we have that 
\begin{eqnarray*}
 \tau_{\lambda_1,c_{-\lambda_1}}^{-1}(c)
&\in& I_{(k_0+k_1)\Lambda_0 +n\Lambda_1 + (k_2-n+1)\Lambda_2}\\ &=&
I_{k\Lambda_0} + U(\bar{\goth{n}})x_{\alpha_1}(-1)^{k_0+k_1+k_2-n+2} +
U(\bar{\goth{n}})x_{\alpha_2}(-1)^{k_0+k_1+n+1} \\
&& + U(\bar{\goth{n}})x_{\alpha_1+\alpha_2}(-1)^{k_0+k_1+1}.
\end{eqnarray*}
Applying the map $\tau_{\lambda_1,c_{-\lambda_1}}$, Lemma
\ref{taulemma} gives us 
\begin{eqnarray*}
c &\in& I_{k\Lambda_1} + U(\bar{\goth{n}})x_{\alpha_1}(-2)^{k_0+k_1+k_2-n+2} +
U(\bar{\goth{n}})x_{\alpha_2}(-1)^{k_0+k_1+n+1}\\
&& + U(\bar{\goth{n}})x_{\alpha_1+\alpha_2}(-2)^{k_0+k_1+1}.
\end{eqnarray*}
So we may write
\begin{eqnarray*}
c = c_1 +c_2x_{\alpha_1}(-1) + c_3x_{\alpha_1}(-2)^{k_0+k_1+k_2-n+2} +
c_4x_{\alpha_2}(-1)^{k_0 + k_1 + n +1} + c_5x_{\alpha_1 +
  \alpha_2}(-2)^{k_0+k_1+1}
\end{eqnarray*}
for some $c_1 \in I_{k\Lambda_0}$ and $c_2,c_3,c_4,c_5 \in U(\bar{\goth{n}})$.
This gives us 
\begin{eqnarray*}
\lefteqn{cx_{\alpha_1}(-1)^nx_{\alpha_1+\alpha_2}(-1)}\\ &=&
c_1x_{\alpha_1}(-1)^nx_{\alpha_1+\alpha_2}(-1) +
c_2x_{\alpha_1}(-1)^{n+1}x_{\alpha_1+\alpha_2}(-1)\\
&& +c_3x_{\alpha_1}(-2)^{k_0+k_1+k_2-n+2}x_{\alpha_1}(-1)^nx_{\alpha_1+\alpha_2}(-1)
\\ &&+
c_4x_{\alpha_2}(-1)^{k_0 + k_1 + n +1}x_{\alpha_1}(-1)^nx_{\alpha_1+\alpha_2}(-1)\\
&& + c_5x_{\alpha_1 + \alpha_2}(-2)^{k_0+k_1+1}x_{\alpha_1}(-1)^nx_{\alpha_1+\alpha_2}(-1).
\end{eqnarray*}
By Lemma \ref{Rcommute}, we may write
\begin{eqnarray*}
c_1x_{\alpha_1}(-1)^nx_{\alpha_1+\alpha_2}(-1) = c_1' +
c_1''x_{\alpha_1 + \alpha_2}(-1)^2
\end{eqnarray*}
for some $c_1' \in I_{k\Lambda_0}$ and $c_1'' \in U(\bar{\goth{n}})$. Clearly
$c_2x_{\alpha_1}(-1)^{n+1}x_{\alpha_1+\alpha_2}(-1)$ is already of the
desired form. We may write 
\begin{eqnarray*}
\lefteqn{c_3x_{\alpha_1}(-2)^{k_0+k_1+k_2-n+2}x_{\alpha_1}(-1)^nx_{\alpha_1+\alpha_2}(-1)}\\
&=& c_3'R_{-1,2(k_0+k_1+k_2-n+2) + n}^1 +
c_3''x_{\alpha_1}(-1)^{n+1}x_{\alpha_1}(-1)^{n+1}x_{\alpha_1 + \alpha_2}(-1)
\end{eqnarray*}
for some $c_3',c_3'' \in U(\bar{\goth{n}})$, which gives us an element of the
desired form. Using Lemma \ref{rootcommute}, we may write
\begin{eqnarray*}
\lefteqn{c_4x_{\alpha_2}(-1)^{k_0 + k_1 + n
    +1}x_{\alpha_1}(-1)^nx_{\alpha_1+\alpha_2}(-1)}\\
&=&
c_4(x_{\alpha_1}(-1)^nx_{\alpha_2}(-1)^{k_0+k_1+n+1}x_{\alpha_1+\alpha_2}(-1)\\
  && + m_1 x_{\alpha_1}(-1)^{n-1}x_{\alpha_1
    +\alpha_2}(-2)x_{\alpha_2}(-1)^{k_0+k_1+n}x_{\alpha_1 +
    \alpha_2}(-1) +\\ && \dots 
+m_nx_{\alpha_1 +
  \alpha_2}(-2)^nx_{\alpha_2}(-1)^{k_0+k_1+1}x_{\alpha_1 +
  \alpha_2}(-1))\\
&=& c_4(m_1'
x_{\alpha_1}(-1)^n[x_{\alpha_1}(0),\dots[x_{\alpha_1}(0),x_{\alpha_2}(-1)^{k_0+k_1+j+2}]\dots ]
\\
&& + \dots + m_n'x_{\alpha_1 +
  \alpha_2}(-2)^n[x_{\alpha_1}(0),x_{\alpha_2}(-1)^{k_0+k_1+2}]\\
&\in& I_{(k_0+1)\Lambda_0 + k_1\Lambda_1 + k_2\Lambda_2}
\end{eqnarray*}
for some constants $m_1, \dots m_n,m_1' \dots m_n' \in
\mathbb{C}$,which is an element of the desired form.  Finally, we have
that
\begin{eqnarray*}
 \lefteqn{c_5x_{\alpha_1 +
     \alpha_2}(-2)^{k_0+k_1+1}x_{\alpha_1}(-1)^nx_{\alpha_1+\alpha_2}(-1)}\\
&=& c_5'[x_{\alpha_1}(0) ,\dots
[x_{\alpha_1}(0),x_{\alpha_2}(-1)^{k_0+k_1+2+n}] \dots ]\\
&\in& I_{(k_0+1)\Lambda_0 + k_1 \Lambda_1 + k_2\Lambda_2}
\end{eqnarray*}
for some $c_5' \in U(\bar{\goth{n}})$, which gives us an element of
the desired form. So we have that
\begin{eqnarray*}
cx_{\alpha_1}(-1)^nx_{\alpha_1 + \alpha_2}(-1) &\in&
I_{(k_0+1)\Lambda_0 + k_1 \Lambda_1 + k_2 \Lambda_2} + U(\bar{\goth{n}})x_{\alpha_1}(-1)^{n+1}x_{\alpha_1+\alpha_2}(-1) \\
&&+U(\bar{\goth{n}})x_{\alpha_1 + \alpha_2}(-1)^2.
\end{eqnarray*}
So we may conclude that $$a
\in I_{(k_0+1)\Lambda_0 + k_1 \Lambda_1 + k_2 \Lambda_2} + U(\bar{\goth{n}})x_{\alpha_1}(-1)^{n+1}x_{\alpha_1+\alpha_2}(-1) +
U(\bar{\goth{n}})x_{\alpha_1 + \alpha_2}(-1)^2$$ proving our claim for $n+1$
when $0 \le n < k_2$. 

We now need to consider the case that $k_2 \le n \le k_0 + k_2$. In
this case, we again have
\begin{eqnarray*}
a\cdot (\underbrace{v_{\Lambda_0} \otimes \dots \otimes
  v_{\Lambda_0}}_{k_0+1 \; \; \mbox{times}} \otimes
\underbrace{v_{\Lambda_1} \otimes \dots \otimes v_{\Lambda_1}}_{k_1\;
  \; \mbox{times}} \otimes \underbrace{v_{\Lambda_2} \otimes \dots
  \otimes v_{\Lambda_2}}_{k_2\; \; \mbox{times}}) =0\\
\end{eqnarray*}
which implies
\begin{eqnarray*}
\lefteqn{(cx_{\alpha_1}(-1)^nx_{\alpha_1+\alpha_2}(-1) +
  dx_{\alpha_1+\alpha_2}(-1)^2 )}\\ && \hspace{1in}\cdot
(\underbrace{v_{\Lambda_0} \otimes \dots \otimes v_{\Lambda_0}}_{k_0+1
  \; \; \mbox{times}} \otimes \underbrace{v_{\Lambda_1} \otimes \dots
  \otimes v_{\Lambda_1}}_{k_1\; \; \mbox{times}} \otimes
\underbrace{v_{\Lambda_2} \otimes \dots \otimes v_{\Lambda_2}}_{k_2\;
  \; \mbox{times}})\\ &=& (a-b) \cdot (\underbrace{v_{\Lambda_0} \otimes
  \dots \otimes v_{\Lambda_0}}_{k_0+1 \; \; \mbox{times}} \otimes
\underbrace{v_{\Lambda_1} \otimes \dots \otimes v_{\Lambda_1}}_{k_1\;
  \; \mbox{times}} \otimes \underbrace{v_{\Lambda_2} \otimes \dots
  \otimes v_{\Lambda_2}}_{k_2\; \; \mbox{times}}))\\ &=& 0.
\end{eqnarray*}
Applying the operator $1 \otimes
\mathcal{Y}_c(e^{\lambda_1},x)^{\otimes k_0+k_2-n}\otimes
\mathcal{Y}_c(e^{\lambda_2},x)^{\otimes n-k_2} \otimes 1^{\otimes
  k_1+k_2}$, we obtain
\begin{eqnarray*}
\lefteqn{(cx_{\alpha_1}(-1)^nx_{\alpha_1+\alpha_2}(-1) +
  dx_{\alpha_1+\alpha_2}(-1)^2 )}\\ && \cdot(v_{\Lambda_0} \otimes
\underbrace{v_{\Lambda_1} \otimes \dots \otimes
  v_{\Lambda_1}}_{k_0+k_2-n \; \; \mbox{times}} \otimes
\underbrace{v_{\Lambda_2} \otimes \dots \otimes
  v_{\Lambda_2}}_{n-k_2\; \; \mbox{times}} \otimes
\underbrace{v_{\Lambda_1} \otimes \dots \otimes v_{\Lambda_1}}_{k_1\;
  \; \mbox{times}}\otimes \underbrace{v_{\Lambda_2} \otimes \dots
  \otimes v_{\Lambda_2}}_{k_2\; \; \mbox{times}})\\ &=&
cx_{\alpha_1}(-1)^nx_{\alpha_1+\alpha_2}(-1)\\ && \cdot(v_{\Lambda_0}
\otimes \underbrace{v_{\Lambda_1} \otimes \dots \otimes
  v_{\Lambda_1}}_{k_0+k_2-n \; \; \mbox{times}} \otimes
\underbrace{v_{\Lambda_2} \otimes \dots \otimes
  v_{\Lambda_2}}_{n-k_2\; \; \mbox{times}} \otimes
\underbrace{v_{\Lambda_1} \otimes \dots \otimes v_{\Lambda_1}}_{k_1\;
  \; \mbox{times}}\otimes \underbrace{v_{\Lambda_2} \otimes \dots
  \otimes v_{\Lambda_2}}_{k_2\; \; \mbox{times}})\\ &=& 0.
\end{eqnarray*}
This implies that 
\begin{align*}
c\cdot (e_{\lambda_1}v_{\Lambda_2} \otimes \underbrace{v_{\Lambda_1} \otimes
  \dots \otimes v_{\Lambda_1}}_{k_0+k_2-n \; \; \mbox{times}} \otimes
\underbrace{e_{\lambda_1}v_{\Lambda_1} \otimes \dots \otimes
  e_{\lambda_1}v_{\Lambda_1}}_{n-k_2\; \; \mbox{times}} \otimes
\underbrace{v_{\Lambda_1} \otimes \dots \otimes v_{\Lambda_1}}_{k_1\;
  \; \mbox{times}}\otimes \underbrace{e_{\lambda_1}v_{\Lambda_1}
  \otimes \dots \otimes e_{\lambda_1}v_{\Lambda_1}}_{k_2\; \;
  \mbox{times}})\\ =0
\end{align*}
which gives
\begin{align*}
e_{\lambda_1}^{\otimes
  k}(\tau_{\lambda_1,c_{-\lambda_1}}^{-1}(c)\cdot (v_{\Lambda_2} \otimes
\underbrace{v_{\Lambda_0} \otimes \dots \otimes
  v_{\Lambda_0}}_{k_0+k_2-n \; \; \mbox{times}} \otimes
\underbrace{v_{\Lambda_1} \otimes \dots \otimes
  v_{\Lambda_1}}_{n-k_2\; \; \mbox{times}} \otimes
\underbrace{v_{\Lambda_0} \otimes \dots \otimes v_{\Lambda_0}}_{k_1\;
  \; \mbox{times}}\otimes \underbrace{v_{\Lambda_1} \otimes \dots
  \otimes v_{\Lambda_1}}_{k_2\; \; \mbox{times}})) \\ =0.
\end{align*}
Since $e_{\lambda_1}^{\otimes k}$ is injective, we have that 
\begin{eqnarray*}
\tau_{\lambda_1,c_{-\lambda_1}}^{-1}(c)\cdot (v_{\Lambda_2} \otimes
\underbrace{v_{\Lambda_0} \otimes \dots \otimes
  v_{\Lambda_0}}_{k_0+k_2-n \; \; \mbox{times}} \otimes
\underbrace{v_{\Lambda_1} \otimes \dots \otimes
  v_{\Lambda_1}}_{n-k_2\; \; \mbox{times}} \otimes
\underbrace{v_{\Lambda_0} \otimes \dots \otimes v_{\Lambda_0}}_{k_1\;
  \; \mbox{times}}\otimes \underbrace{v_{\Lambda_1} \otimes \dots
  \otimes v_{\Lambda_1}}_{k_2\; \; \mbox{times}}) =0.
\end{eqnarray*}
Now, since the total charge of
$\tau_{\lambda_1,c_{-\lambda_1}}^{-1}(c)$ is less than the total
charge of $a$, we have that 
\begin{eqnarray*}
\tau_{\lambda_1,c_{-\lambda_1}}^{-1}(c) &\in&
I_{(k_0+k_1+k_2-n)\Lambda_0 + n\Lambda_1 + \Lambda_2} \\
&=& I_{k\Lambda_0}
+ U(\bar{\goth{n}})x_{\alpha_1}(-1)^{k_0+k_1+k_2-n+2} +
U(\bar{\goth{n}})x_{\alpha_2}(-1)^{k} \\
&& + U(\bar{\goth{n}})x_{\alpha_1 + \alpha_2}(-1)^{k_0+k_1+k_2-n+1}.
\end{eqnarray*}
Applying $\tau_{\lambda_1,c_{-\lambda_1}}$, Lemma \ref{taulemma} gives us
\begin{eqnarray*}
c \in I_{k\Lambda_1} +
U(\bar{\goth{n}})x_{\alpha_1}(-2)^{k_0+k_1+k_2-n+2} +
U(\bar{\goth{n}})x_{\alpha_2}(-1)^{k} + U(\bar{\goth{n}})x_{\alpha_1 +
  \alpha_2}(-2)^{k_0+k_1+k_2-n+1}.
\end{eqnarray*}
So we may write
\begin{eqnarray*}
c=c_1+c_2x_{\alpha_1}(-1) + c_3 x_{\alpha_1}(-2)^{k_0+k_1+k_2-n+2} +
c_4 x_{\alpha_2}(-1)^{k}  + c_5 x_{\alpha_1 + \alpha_2}(-2)^{k_0+k_1+k_2-n+1}
\end{eqnarray*}
for some $c_1 \in I_{k\Lambda_0}$ and $c_2,c_3,c_4,c_5 \in
U(\bar{\goth{n}})$, so that 
\begin{eqnarray*}
\lefteqn{cx_{\alpha_1}(-1)^nx_{\alpha_1+\alpha_2}(-1)}\\
&=& c_1x_{\alpha_1}(-1)^nx_{\alpha_1+\alpha_2}(-1)+c_2x_{\alpha_1}(-1)^{n+1}x_{\alpha_1
  + \alpha_2}(-1)\\
&& + c_3 x_{\alpha_1}(-2)^{k_0+k_1+k_2-n+2}x_{\alpha_1}(-1)^nx_{\alpha_1+\alpha_2}(-1) +
c_4 x_{\alpha_2}(-1)^{k}x_{\alpha_1}(-1)^nx_{\alpha_1+\alpha_2}(-1) \\
 && + c_5 x_{\alpha_1 + \alpha_2}(-2)^{k_0+k_1+k_2-n+1}x_{\alpha_1}(-1)^nx_{\alpha_1+\alpha_2}(-1).
\end{eqnarray*}
As before, we have that
\begin{eqnarray*}
\lefteqn{c_1x_{\alpha_1}(-1)^nx_{\alpha_1+\alpha_2}(-1),
  c_2x_{\alpha_1}(-1)^{n+1}x_{\alpha_1+\alpha_2}(-1)}\\ &\in&
I_{(k_0+1)\Lambda_0 + k_1\Lambda_1 + k_2\Lambda_2} +
U(\bar{\goth{n}})x_{\alpha_1}(-1)^{n+1}x_{\alpha_1+\alpha_2}(-1) +
U(\bar{\goth{n}})x_{\alpha_1 + \alpha_2}(-1)^2
\end{eqnarray*}
and so have the desired form.  We also have that
\begin{eqnarray*}
\lefteqn{c_3
  x_{\alpha_1}(-2)^{k_0+k_1+k_2-n+2}x_{\alpha_1}(-1)^nx_{\alpha_1+\alpha_2}(-1)}\\
&=& c_3'R_{-1,2(k_0+k_1+k_2-n+2) + n}^1 +
c_3''x_{\alpha_1}(-1)^{n+1}x_{\alpha_1 + \alpha_2}(-1)
\end{eqnarray*}
for some $c_3',c_3'' \in U(\bar{\goth{n}})$, which is of the desired
form. Using Lemma \ref{rootcommute}, we
have 
\begin{eqnarray*}
\lefteqn{c_4
  x_{\alpha_2}(-1)^{k}x_{\alpha_1}(-1)^nx_{\alpha_1+\alpha_2}(-1)}\\
&=& c_4(x_{\alpha_1}(-1)^n
x_{\alpha_2}(-1)^{k}x_{\alpha_1+\alpha_2}(-1) + m_1
x_{\alpha_1}(-1)^{n-1}x_{\alpha_1 +
  \alpha_2}(-2)x_{\alpha_2}(-1)^{k-1}x_{\alpha_1 + \alpha_2}(-1)\\
&&+ \dots +m_nx_{\alpha_1 + \alpha_2}(-2)^n
x_{\alpha_2}(-1)^{k-n}x_{\alpha_1+\alpha_2}(-1))\\
&=& c_4(m_0' x_{\alpha_1}(-1)^n
[x_{\alpha_1}(0),R_{-1,k+1}^2] + m_1'
x_{\alpha_1}(-1)^{n-1}[x_{\alpha_1}(0),[x_{\alpha_1}(0),R_{-1,k+2}^2]]
\\
&& + \dots + m_n'[x_{\alpha_1}(0),\dots
[x_{\alpha_1}(0),R_{-1,k+n+1}^2]\dots ]) + c_4'x_{\alpha_1+\alpha_2}(-1)^2
\end{eqnarray*}
for some $m_1,\dots m_n,m_0', \dots m_n' \in \mathbb{C}$ and $c_4' \in
U(\bar{\goth{n}})$, which is of the desired form.   Finally,
\begin{eqnarray*}
\lefteqn{c_5 x_{\alpha_1 +
  \alpha_2}(-2)^{k_0+k_1+k_2-n+1}x_{\alpha_1}(-1)^nx_{\alpha_1+\alpha_2}(-1)}\\
&=&  c_5 [x_{\alpha_2}(0),\dots [x_{\alpha_2}(0),R_{-1,2(k-n)+n+1}^1]
\dots ] + c_5' x_{\alpha_1}(-1)^{n+1}x_{\alpha_1+\alpha_2}(-1)
\end{eqnarray*}
for some $c_5' \in U(\bar{\goth{n}})$. So we have that 
\begin{eqnarray*}
cx_{\alpha_1}(-1)^nx_{\alpha_1+\alpha_2}(-1)  &\in& I_{(k_0+1)\Lambda_0 +
  k_1\Lambda_1 + k_2\Lambda_2} +
U(\bar{\goth{n}})x_{\alpha_1+\alpha_2}(-1)^{n+1}x_{\alpha_1 + \alpha_2}(-1)\\
&& + U(\bar{\goth{n}})x_{\alpha_1 + \alpha_2}(-1)^2,
\end{eqnarray*}
which gives us that $$a \in I_{(k_0+1)\Lambda_0 + k_1\Lambda_1 +
  k_2\Lambda_2} +
U(\bar{\goth{n}})x_{\alpha_1+\alpha_2}(-1)^{n+1}x_{\alpha_1 +
  \alpha_2}(-1) + U(\bar{\goth{n}})x_{\alpha_1 + \alpha_2}(-1)^2,$$
completing our induction on $n$.  We've thus shown that
\begin{eqnarray*}
a = b + cx_{\alpha_1}(-1)^nx_{\alpha_1 + \alpha_2}(-1)^m +
dx_{\alpha_1 + \alpha_2}(-1)^{m+1}
\end{eqnarray*}
for some $b \in  I_{(k_0+1)\Lambda_0 +
  k_1\Lambda_1 + k_2\Lambda_2}$ and $c,d \in U(\bar{\goth{n}})$ holds when
$m=1$ and $0 \le n \le k_0+k_2+2-m$.

Now, we induct on $m$. Assume that we've shown, for some $m \in \mathbb{N}$ satisfying $1 \le m \le
k_0$ and all $0 \le n \le k_0+k_2+2-m$ that there exist $b \in  I_{(k_0+1)\Lambda_0 +
  k_1\Lambda_1 + k_2\Lambda_2}$ and $c,d \in U(\bar{\goth{n}})$ such that
\begin{eqnarray} \label{mninduction}
a = b + cx_{\alpha_1}(-1)^nx_{\alpha_1 + \alpha_2}(-1)^m +
dx_{\alpha_1 + \alpha_2}(-1)^{m+1}.
\end{eqnarray}
 We show that
this holds for $m+1$ as well.  
We have by (\ref{mninduction}) that
\begin{eqnarray*}
a = b + cx_{\alpha_1}(-1)^{k_0+k_2+2-m}x_{\alpha_1 + \alpha_2}(-1)^m +
dx_{\alpha_1+\alpha_2}(-1)^{m+1}
\end{eqnarray*}
for some $b \in  I_{(k_0+1)\Lambda_0 +
  k_1\Lambda_1 + k_2\Lambda_2}$ and $c,d \in U(\bar{\goth{n}})$.  Here, we
have that 
\begin{eqnarray*}
cx_{\alpha_1}(-1)^{k_0+k_2+2-m}x_{\alpha_1 + \alpha_2}(-1)^m &=&
c'[x_{\alpha_2}(0),\dots
[x_{\alpha_2}(0),x_{\alpha_1}(-1)^{k_0+k_2+1}] \dots ]\\
 &\in& I_{(k_0+1)\Lambda_0 + k_1 \Lambda_1 + k_2\Lambda_2},
\end{eqnarray*}
for some $c' \in U(\bar{\goth{n}})$, so $a \in I_{(k_0+1)\Lambda_0 +
  k_1\Lambda_1 + k_2\Lambda_2} +
U(\bar{\goth{n}})x_{\alpha_1+\alpha_2}(-1)^{m+1}$. In this case, we
have that
\begin{eqnarray*}
a = b' + c'x_{\alpha_1}(-1)^0x_{\alpha_1+\alpha_2}(-1)^{m+1} + d'x_{\alpha_1+\alpha_2}(-1)^{m+2}
\end{eqnarray*}
for some $b' \in I_{(k_0+1)\Lambda_0 + k_1 \Lambda_1 + k_2\Lambda_2}$
and $c',d' \in U(\bar{\goth{n}})$ (here, $d=0$). So our claim holds for
$m+1$ and $n=0$. Now, assume that we have shown 
\begin{eqnarray*}
a = b + cx_{\alpha_1}(-1)^nx_{\alpha_1+\alpha_2}(-1)^{m+1} +
dx_{\alpha_1+\alpha_2}(-1)^{m+2}
\end{eqnarray*}
for some $b \in I_{(k_0+1)\Lambda_0 + k_1 \Lambda_1 + k_2\Lambda_2}$,
$c,d \in U(\bar{\goth{n}})$ and $n$ satisfying $0 \le n \le k_0+k_2-m$.  As with
$m=1$, we need to consider the cases $0 \le n < k_2$ and $k_2 \le n
\le k_0+k_2-m$ separately. First, suppose $0 \le n < k_2$. Here, we
have that
\begin{eqnarray*}
\lefteqn{(cx_{\alpha_1}(-1)^nx_{\alpha_1+\alpha_2}(-1)^{m+1} +
  dx_{\alpha_1+\alpha_2}(-1)^{m+2})}\\ && \hspace{1in}
\cdot(\underbrace{v_{\Lambda_0} \otimes \dots \otimes
  v_{\Lambda_0}}_{k_0+1 \; \; \mbox{times}} \otimes
\underbrace{v_{\Lambda_1} \otimes \dots \otimes v_{\Lambda_1}}_{k_1\;
  \; \mbox{times}} \otimes \underbrace{v_{\Lambda_2} \otimes \dots
  \otimes v_{\Lambda_2}}_{k_2\; \; \mbox{times}}) \\ &=& (a-b)\cdot 
(\underbrace{v_{\Lambda_0} \otimes \dots \otimes v_{\Lambda_0}}_{k_0+1
  \; \; \mbox{times}} \otimes \underbrace{v_{\Lambda_1} \otimes \dots
  \otimes v_{\Lambda_1}}_{k_1\; \; \mbox{times}} \otimes
\underbrace{v_{\Lambda_2} \otimes \dots \otimes v_{\Lambda_2}}_{k_2\;
  \; \mbox{times}})\\ &=& 0.
\end{eqnarray*}
Applying the operator $1^{\otimes m+1} \otimes
\mathcal{Y}_c(e^{\lambda_1},x)^{\otimes k_0-m} \otimes 1^{\otimes k_1}
\otimes \mathcal{Y}_c(e^{\lambda_1},x)^{\otimes k_2-n} \otimes
1^{\otimes n}$, we have 
\begin{align*}
cx_{\alpha_1}(-1)^nx_{\alpha_1+\alpha_2}(-1)^{m+1}\cdot 
(&\underbrace{v_{\Lambda_0} \otimes \dots \otimes v_{\Lambda_0}}_{m+1
  \; \; \mbox{times}} \otimes \underbrace{v_{\Lambda_1} \otimes \dots
  \otimes v_{\Lambda_1}}_{k_0+k_1-m\; \; \mbox{times}} \otimes
\\ & \hspace{0.5in}\underbrace{e_{\lambda_1}v_{\Lambda_2} \otimes
  \dots \otimes e_{\lambda_1}v_{\Lambda_2}}_{k_2-n \; \;
  \mbox{times}}\otimes \underbrace{v_{\Lambda_2} \otimes \dots \otimes
  v_{\Lambda_2}}_{n \; \; \mbox{times}}) = 0,
\end{align*}
which implies that
\begin{align*}
c\cdot  (&\underbrace{e_{\lambda_1}v_{\Lambda_2} \otimes \dots \otimes
  e_{\lambda_1}v_{\Lambda_2}}_{m+1 \; \; \mbox{times}} \otimes
\underbrace{v_{\Lambda_1} \otimes \dots \otimes
  v_{\Lambda_1}}_{k_0+k_1-m\; \; \mbox{times}} \otimes
\\ &\hspace{0.5in} \underbrace{e_{\lambda_1}v_{\Lambda_2} \otimes
  \dots \otimes e_{\lambda_1}v_{\Lambda_2}}_{k_2-n \; \;
  \mbox{times}}\otimes \underbrace{e_{\lambda_1}v_{\Lambda_1} \otimes
  \dots \otimes e_{\lambda_1}v_{\Lambda_1}}_{n \; \; \mbox{times}}) =
0
\end{align*}
and so we have
\begin{eqnarray*}
e_{\lambda_1}^{\otimes k}(\tau_{\lambda_1,c_{-\lambda_1}}^{-1}(c)\cdot 
(\underbrace{v_{\Lambda_2} \otimes \dots \otimes v_{\Lambda_2}}_{m+1
  \; \; \mbox{times}} \otimes \underbrace{v_{\Lambda_0} \otimes \dots
  \otimes v_{\Lambda_0}}_{k_0+k_1-m\; \; \mbox{times}} \otimes
\underbrace{v_{\Lambda_2} \otimes \dots \otimes v_{\Lambda_2}}_{k_2-n
  \; \; \mbox{times}}\otimes \underbrace{v_{\Lambda_1} \otimes \dots
  \otimes v_{\Lambda_1}}_{n \; \; \mbox{times}})) = 0.
\end{eqnarray*}
Since $e_{\lambda_1}^{\otimes k}$ is injective, we have that 
\begin{eqnarray*}
\tau_{\lambda_1,c_{-\lambda_1}}^{-1}(c) \cdot (\underbrace{v_{\Lambda_2}
  \otimes \dots \otimes v_{\Lambda_2}}_{m+1 \; \; \mbox{times}}
\otimes \underbrace{v_{\Lambda_0} \otimes \dots \otimes
  v_{\Lambda_0}}_{k_0+k_1-m\; \; \mbox{times}} \otimes
\underbrace{v_{\Lambda_2} \otimes \dots \otimes v_{\Lambda_2}}_{k_2-n
  \; \; \mbox{times}}\otimes \underbrace{v_{\Lambda_1} \otimes \dots
  \otimes v_{\Lambda_1}}_{n \; \; \mbox{times}}) = 0.
\end{eqnarray*}
Now, since $\tau_{\lambda_1,c_{-\lambda_1}}^{-1}(c)$ has lower total
charge than $a$, we have that 
\begin{eqnarray*}
\tau_{\lambda_1,c_{-\lambda_1}}^{-1}(c) &\in& I_{(k_0+k_1-m)\Lambda_0 +
  n\Lambda_1 + (k_2-n+m+1)\Lambda_2}\\
 &=& I_{k\Lambda_0} +
U(\bar{\goth{n}})x_{\alpha_1}(-1)^{k-n+1} +
U(\bar{\goth{n}})x_{\alpha_2}(-1)^{k_0+k_1-m+n+1}\\
&& + U(\bar{\goth{n}})x_{\alpha_1 + \alpha_2}(-1)^{k_0+k_1-m+1}.
\end{eqnarray*}
Applying the map $\tau_{\lambda_1,c_{-\lambda_1}}$, by Lemma
\ref{taulemma} we have that 
\begin{eqnarray*}
c \in I_{k\Lambda_1} + U(\bar{\goth{n}})x_{\alpha_1}(-2)^{k-n+1} +
U(\bar{\goth{n}})x_{\alpha_2}(-1)^{k_0+k_1-m+n+1} +
U(\bar{\goth{n}})x_{\alpha_1 + \alpha_2}(-2)^{k_0+k_1-m+1}.
\end{eqnarray*}
So we may write
\begin{eqnarray*}
c = c_1 + c_2 x_{\alpha_1}(-1) + c_3x_{\alpha_1}(-2)^{k-n+1} + c_4
x_{\alpha_2}(-1)^{k_0+k_1-m+n+1} + c_5 x_{\alpha_1 + \alpha_2}(-2)^{k_0+k_1-m+1}
\end{eqnarray*}
for some $c_1 \in I_{k\Lambda_0}$ and $c_2,c_3,c_4,c_5 \in
U(\bar{\goth{n}})$, and we have
\begin{eqnarray*}
\lefteqn{cx_{\alpha_1}(-1)^nx_{\alpha_1+\alpha_2}(-1)^{m+1} =
c_1x_{\alpha_1}(-1)^nx_{\alpha_1+\alpha_2}(-1)^{m+1}}\\
&& + c_2 x_{\alpha_1}(-1)^{n+1}x_{\alpha_1+\alpha_2}(-1)^{m+1} +
c_3x_{\alpha_1}(-2)^{k-n+1}x_{\alpha_1}(-1)^nx_{\alpha_1+\alpha_2}(-1)^{m+1}\\
&& + c_4
x_{\alpha_2}(-1)^{k_0+k_1-m+n+1}x_{\alpha_1}(-1)^nx_{\alpha_1+\alpha_2}(-1)^{m+1}\\
&& + c_5 x_{\alpha_1 + \alpha_2}(-2)^{k_0+k_1-m+1}x_{\alpha_1}(-1)^nx_{\alpha_1+\alpha_2}(-1)^{m+1}.
\end{eqnarray*}
By Lemma \ref{Rcommute}, we have that
$$c_1x_{\alpha_1}(-1)^nx_{\alpha_1+\alpha_2}(-1)^{m+1} \in
I_{k\Lambda_0} + U(\bar{\goth{n}})x_{\alpha_1 + \alpha_2}(-1)^{m+2}.$$ The summand $c_2
x_{\alpha_1}(-1)^{n+1}x_{\alpha_1+\alpha_2}(-1)^{m+1}$ is clearly of
the desired form. We have that 
\begin{eqnarray*}
\lefteqn{c_3x_{\alpha_1}(-2)^{k-n+1}x_{\alpha_1}(-1)^nx_{\alpha_1+\alpha_2}(-1)^{m+1}}\\
&=& c_3'x_{\alpha_1+\alpha_2}(-1)^{m+1}R_{-1,2(k-n+1) + n}^1 + c_3''x_{\alpha_1}(-1)^{n+1}x_{\alpha_1+\alpha_2}(-1)^{m+1}
\end{eqnarray*}
for some $c_3',c_3'' \in U(\bar{\goth{n}})$, which is clearly of the
desired form.  For the next summand, by Lemma \ref{rootcommute} we
have that
\begin{eqnarray*}
 \lefteqn{c_4x_{\alpha_2}(-1)^{k_0+k_1-m+n+1}x_{\alpha_1}(-1)^nx_{\alpha_1+\alpha_2}(-1)^{m+1}}\\
&=&
c_4(x_{\alpha_1}(-1)^nx_{\alpha_2}(-1)^{k_0+k_1-m+n+1}x_{\alpha_1+\alpha_2}(-1)^{m+1}\\
&&+ m_1 x_{\alpha_1}(-1)^{n-1}x_{\alpha_1 +
  \alpha_2}(-2)x_{\alpha_2}(-1)^{k_0+k_1-m+n}x_{\alpha_1+\alpha_2}(-1)^{m+1}\\
&& + \dots +
m_nx_{\alpha_1+\alpha_2}(-2)^nx_{\alpha_2}(-1)^{k_0+k_1-m+1}x_{\alpha_1+\alpha_2}(-1)^{m+1})\\
&=& c_4(m_0' x_{\alpha_1}(-1)^n[x_{\alpha_1}(0), \dots
[x_{\alpha_1}(0),x_{\alpha_2}(-1)^{k_0+k_1+2+n}] \dots ]\\
&& +m_1'x_{\alpha_1}(-1)^{n-1}x_{\alpha_1+\alpha_2}(-2)[x_{\alpha_1}(0),
\dots [x_{\alpha_1}(0),x_{\alpha_2}(-1)^{k_0+k_1+n+1}] \dots ]\\
&& + \dots + m_n'x_{\alpha_1+\alpha_2}(-1)^n[x_{\alpha_1}, \dots
[x_{\alpha_1}(0),x_{\alpha_2}(-1)^{k_0+k_1+2}] \dots ])
\end{eqnarray*}
for some $m_1,\dots m_n, m_0', \dots m_n' \in \mathbb{C}$, which is an
element of $I_{(k_0+1)\Lambda_0 + k_1\Lambda_1 + k_2\Lambda_2}$.
Finally, applying Lemma \ref{rootcommute}, we have that
\begin{eqnarray*}
\lefteqn{c_5 x_{\alpha_1 + \alpha_2}(-2)^{k_0+k_1-m+1}x_{\alpha_1}(-1)^nx_{\alpha_1+\alpha_2}(-1)^{m+1}}\\
&=& c_5( m_0
[x_{\alpha_1}(-1)^{k_0+k_1-m+1+n},x_{\alpha_2}(-1)^{k_0+k_1-m+1}] +\\
&& m_1x_{\alpha_2}(-1)x_{\alpha_1 +
  \alpha_2}(-2)^{k_0+k_1-m}x_{\alpha_1}(-1)^{n+1} + \\
&& \dots+m_{k_0+k_1-m}x_{\alpha_2}(-1)^{k_0+k_1-m}x_{\alpha_1+\alpha_2}(-2)x_{\alpha_1}(-1)^{k_0+k_1-m+n})x_{\alpha_1
  + \alpha_2}(-1)
\end{eqnarray*}
for some $m_0, \dots , m_{k_0+k_1-m} \in \mathbb{C}$.  Notice, all the
terms in the right hand side are elements of \\
$$U(\bar{\goth{n}})x_{\alpha_1}(-1)^{n+1}x_{\alpha_1 + \alpha_2}(-1)^{m+1}$$
except for the term\\
$$m_0c_5x_{\alpha_1}(-1)^{k_0+k_1-m+1+n}x_{\alpha_2}(-1)^{k_0+k_1-m+1}x_{\alpha_1+\alpha_2}(-1)^{m+1}.$$ This
term, however, may be written as
\begin{eqnarray*}
\lefteqn{m_0c_5x_{\alpha_1}(-1)^{k_0+k_1-m+1+n}x_{\alpha_2}(-1)^{k_0+k_1-m+1}x_{\alpha_1+\alpha_2}(-1)^{m+1}}\\
&=& c_5'x_{\alpha_1}(-1)^{k_0+k_1-m+1+n}[x_{\alpha_1}(0), \dots
[x_{\alpha_1}(0),x_{\alpha_2}(-1)^{k_0+k_1+2}] \dots ]\\
&\in& I_{(k_0+1)\Lambda_0 + k_1\Lambda_1 + k_2\Lambda_2}
\end{eqnarray*}
for some $c_5' \in U(\bar{\goth{n}})$ Hence, we have
that
\begin{eqnarray*} 
cx_{\alpha_1}(-1)^n
x_{\alpha_1+\alpha_2}(-1)^{m+1} &\in& I_{(k_0+1)\Lambda_0 + k_1
  \Lambda_1 + k_2 \Lambda_2} +
U(\bar{\goth{n}})x_{\alpha_1}(-1)^{n+1}x_{\alpha_1+\alpha_2}(-1)^{m+1}\\
&& +U(\bar{\goth{n}})x_{\alpha_1+\alpha_2}(-1)^{m+2}
\end{eqnarray*}
and is of the desired form. From this, we conclude that we may write
\begin{eqnarray*}
a = b' + c'x_{\alpha_1}(-1)^{n+1}x_{\alpha_1+\alpha_2}(-1)^{m+1} +
d'x_{\alpha_1+\alpha_2}(-1)^{m+2}
\end{eqnarray*}
for some $b' \in I_{(k_0+1)\Lambda_0 + k_1\Lambda_1 + k_2\Lambda_2}$
and $c',d' \in U(\bar{\goth{n}})$, completing our induction for $0 \le n <
k_2$.

We now assume that $k_2 \le n \le k_0+k_2-m$. As before, we have that
\begin{eqnarray*}
a = b+cx_{\alpha_1}(-1)^nx_{\alpha_1+\alpha_2}(-1)^{m+1} + dx_{\alpha_1+\alpha_2}(-1)^{m+2}
\end{eqnarray*}
for some $b \in I_{(k_0+1)\Lambda_0 + k_1\Lambda_1 + k_2\Lambda_2}$
and  $c,d \in U(\bar{\goth{n}})$. So, as before, we have that
\begin{eqnarray*}
\lefteqn{(cx_{\alpha_1}(-1)^nx_{\alpha_1+\alpha_2}(-1)^{m+1} +
  dx_{\alpha_1+\alpha_2}(-1)^{m+2})}\\ && \hspace{1in}\cdot
(\underbrace{v_{\Lambda_0} \otimes \dots \otimes v_{\Lambda_0}}_{k_0+1
  \; \; \mbox{times}} \otimes \underbrace{v_{\Lambda_1} \otimes \dots
  \otimes v_{\Lambda_1}}_{k_1\; \; \mbox{times}} \otimes
\underbrace{v_{\Lambda_2} \otimes \dots \otimes v_{\Lambda_2}}_{k_2\;
  \; \mbox{times}}) \\ &=& (a-b)\cdot  (\underbrace{v_{\Lambda_0} \otimes
  \dots \otimes v_{\Lambda_0}}_{k_0+1 \; \; \mbox{times}} \otimes
\underbrace{v_{\Lambda_1} \otimes \dots \otimes v_{\Lambda_1}}_{k_1\;
  \; \mbox{times}} \otimes \underbrace{v_{\Lambda_2} \otimes \dots
  \otimes v_{\Lambda_2}}_{k_2\; \; \mbox{times}})\; \; \; \; \; \; \;
\; \; \;\; \; \; \; \; \; \; \; \; \; \; \; \; \\ &=& 0.
\end{eqnarray*}
Applying the operator $1^{\otimes (m+1)} \otimes
\mathcal{Y}_c(e^{\lambda_1},x)^{\otimes (k_0+k_2-m-n)} \otimes
\mathcal{Y}_c(e^{\lambda_2},x)^{\otimes (n-k_2)} \otimes 1^{\otimes
  (k_1+k_2)}$, we have that
\begin{eqnarray*}
\lefteqn{(cx_{\alpha_1}(-1)^nx_{\alpha_1+\alpha_2}(-1)^{m+1} +
  dx_{\alpha_1+\alpha_2}(-1)^{m+2})}\\ &&\cdot
(\underbrace{v_{\Lambda_0} \otimes \dots \otimes v_{\Lambda_0}}_{m+1
  \; \; \mbox{times}} \otimes \underbrace{v_{\Lambda_1} \otimes \dots
  \otimes v_{\Lambda_1}}_{k_0+k_2-m-n \; \; \mbox{times}} \otimes
\underbrace{v_{\Lambda_2} \otimes \dots \otimes
  v_{\Lambda_2}}_{n-k_2\; \; \mbox{times}} \otimes
\underbrace{v_{\Lambda_1} \otimes \dots \otimes v_{\Lambda_1}}_{k_1\;
  \; \mbox{times}} \otimes \underbrace{v_{\Lambda_2} \otimes \dots
  \otimes v_{\Lambda_2}}_{k_2\; \; \mbox{times}}) \\ &=&
(cx_{\alpha_1}(-1)^nx_{\alpha_1+\alpha_2}(-1)^{m+1})\\ &&\cdot
(\underbrace{v_{\Lambda_0} \otimes \dots \otimes v_{\Lambda_0}}_{m+1
  \; \; \mbox{times}} \otimes \underbrace{v_{\Lambda_1} \otimes \dots
  \otimes v_{\Lambda_1}}_{k_0+k_2-m-n \; \; \mbox{times}} \otimes
\underbrace{v_{\Lambda_2} \otimes \dots \otimes
  v_{\Lambda_2}}_{n-k_2\; \; \mbox{times}} \otimes
\underbrace{v_{\Lambda_1} \otimes \dots \otimes v_{\Lambda_1}}_{k_1\;
  \; \mbox{times}} \otimes \underbrace{v_{\Lambda_2} \otimes \dots
  \otimes v_{\Lambda_2}}_{k_2\; \; \mbox{times}}) \\ &=& 0.
\end{eqnarray*}
This implies that 
\begin{align*}
cx_{\alpha_1}(-1)^{n}\cdot (&\underbrace{e_{\lambda_1}v_{\Lambda_2} \otimes
  \dots \otimes e_{\lambda_1}v_{\Lambda_2}}_{m+1 \; \; \mbox{times}}
\otimes \underbrace{v_{\Lambda_1} \otimes \dots \otimes
  v_{\Lambda_1}}_{k_0+k_2-m-n \; \; \mbox{times}}
\otimes\\ &\hspace{0.5in} \underbrace{v_{\Lambda_2} \otimes \dots
  \otimes v_{\Lambda_2}}_{n-k_2\; \; \mbox{times}} \otimes
\underbrace{v_{\Lambda_1} \otimes \dots \otimes v_{\Lambda_1}}_{k_1\;
  \; \mbox{times}} \otimes \underbrace{v_{\Lambda_2} \otimes \dots
  \otimes v_{\Lambda_2}}_{k_2\; \; \mbox{times}}) = 0
\end{align*}
which implies
\begin{align*}
c\cdot (&\underbrace{e_{\lambda_1}v_{\Lambda_2} \otimes \dots \otimes
  e_{\lambda_1}v_{\Lambda_2}}_{m+1 \; \; \mbox{times}} \otimes
\underbrace{v_{\Lambda_1} \otimes \dots \otimes
  v_{\Lambda_1}}_{k_0+k_2-m-n \; \; \mbox{times}}
\otimes\\ & \hspace{0.3in} \underbrace{e_{\lambda_1}v_{\Lambda_1}
  \otimes \dots \otimes e_{\lambda_1}v_{\Lambda_1}}_{n-k_2\; \;
  \mbox{times}} \otimes \underbrace{v_{\Lambda_1} \otimes \dots
  \otimes v_{\Lambda_1}}_{k_1\; \; \mbox{times}} \otimes
\underbrace{e_{\lambda_1}v_{\Lambda_1} \otimes \dots \otimes
  e_{\lambda_1}v_{\Lambda_1}}_{k_2\; \; \mbox{times}}) = 0.
\end{align*}
From this, we have that 
\begin{align*}
e_{\lambda_1}^{\otimes
  k}(&\tau_{\lambda_1,c_{-\lambda_1}}^{-1}(c)\cdot (\underbrace{v_{\Lambda_2}
  \otimes \dots \otimes v_{\Lambda_2}}_{m+1 \; \; \mbox{times}}
\otimes \underbrace{v_{\Lambda_0} \otimes \dots \otimes
  v_{\Lambda_0}}_{k_0+k_2-m-n \; \; \mbox{times}}
\otimes\\ & \hspace{0.3in} \underbrace{v_{\Lambda_1} \otimes \dots
  \otimes v_{\Lambda_1}}_{n-k_2\; \; \mbox{times}} \otimes
\underbrace{v_{\Lambda_0} \otimes \dots \otimes v_{\Lambda_0}}_{k_1\;
  \; \mbox{times}} \otimes \underbrace{v_{\Lambda_1} \otimes \dots
  \otimes v_{\Lambda_1}}_{k_2\; \; \mbox{times}})) = 0.
\end{align*}
Since $e_{\lambda_1}^{\otimes k}$ is injective, we have that
\begin{align*}
\tau_{\lambda_1,c_{-\lambda_1}}^{-1}(c)\cdot (&\underbrace{v_{\Lambda_2}
  \otimes \dots \otimes v_{\Lambda_2}}_{m+1 \; \; \mbox{times}}
\otimes \underbrace{v_{\Lambda_0} \otimes \dots \otimes
  v_{\Lambda_0}}_{k_0+k_2-m-n \; \; \mbox{times}} \otimes
\\ & \hspace{0.3in}\underbrace{v_{\Lambda_1} \otimes \dots \otimes
  v_{\Lambda_1}}_{n-k_2\; \; \mbox{times}} \otimes
\underbrace{v_{\Lambda_0} \otimes \dots \otimes v_{\Lambda_0}}_{k_1\;
  \; \mbox{times}} \otimes \underbrace{v_{\Lambda_1} \otimes \dots
  \otimes v_{\Lambda_1}}_{k_2\; \; \mbox{times}}) = 0.
\end{align*}
Since the total charge of $\tau_{\lambda_1,c_{-\lambda_1}}^{-1}(c)$ is
less than the total charge of $a$, we have that
\begin{eqnarray*}
\tau_{\lambda_1,c_{-\lambda_1}}^{-1}(c) &\in&
I_{(k_0+k_1+k_2-m-n)\Lambda_0 + n\Lambda_1 + (m+1)\Lambda_2}\\
&=& I_{k\Lambda_0} + U(\bar{\goth{n}})x_{\alpha_1}(-1)^{k_0+k_1+k_2-n+2}
+U(\bar{\goth{n}})x_{\alpha_2}(-1)^{k_0+k_1+k_2-m+1} \\
&&+ U(\bar{\goth{n}})x_{\alpha_1
  + \alpha_2}(-1)^{k_0+k_1+k_2-m-n+1}.
\end{eqnarray*}
Applying $\tau_{\lambda_1,c_{-\lambda_1}}$ to both sides, Lemma
\ref{taulemma} gives us
\begin{eqnarray*}
c &\in& I_{k\Lambda_1} + U(\bar{\goth{n}})x_{\alpha_1}(-2)^{k_0+k_1+k_2-n+2}
+U(\bar{\goth{n}})x_{\alpha_2}(-1)^{k_0+k_1+k_2-m+1}\\
&& + U(\bar{\goth{n}})x_{\alpha_1
  + \alpha_2}(-2)^{k_0+k_1+k_2-m-n+1}.
\end{eqnarray*}
So we may write
\begin{eqnarray*}
c &=& c_1 + c_2x_{\alpha_1}(-1) + c_3 x_{\alpha_1}(-2)^{k_0+k_1+k_2-n+2}
+ c_4 x_{\alpha_2}(-1)^{k_0+k_1+k_2-m+1}\\
&& + c_5 x_{\alpha_1
  + \alpha_2}(-2)^{k_0+k_1+k_2-m-n+1}
\end{eqnarray*}
for some $c_1 \in I_{k\Lambda_0}$ and $c_2,c_3,c_4,c_5 \in
U(\bar{\goth{n}})$, which gives
\begin{eqnarray*}
cx_{\alpha_1}(-1)^nx_{\alpha_1+\alpha_2}(-1)^{m+1} &=&
c_1x_{\alpha_1}(-1)^nx_{\alpha_1+\alpha_2}(-1)^{m+1} +
c_2x_{\alpha_1}(-1)^{n+1}x_{\alpha_1+\alpha_2}(-1)^{m+1}\\
&& + c_3
x_{\alpha_1}(-2)^{k_0+k_1+k_2-n+2}x_{\alpha_1}(-1)^nx_{\alpha_1+\alpha_2}(-1)^{m+1}\\
&&
+ c_4
x_{\alpha_2}(-1)^{k_0+k_1+k_2-m+1}x_{\alpha_1}(-1)^nx_{\alpha_1+\alpha_2}(-1)^{m+1}\\
&& + c_5 x_{\alpha_1
  + \alpha_2}(-2)^{k_0+k_1+k_2-m-n+1}x_{\alpha_1}(-1)^nx_{\alpha_1+\alpha_2}(-1)^{m+1}.
\end{eqnarray*}
As before, we analyze each summand and show that it is of the desired
form.  By Corollary \ref{Rcommute}, we have that
\begin{eqnarray*}
c_1x_{\alpha_1}(-1)^nx_{\alpha_1+\alpha_2}(-1)^{m+1} \in
I_{k\Lambda_0} + U(\bar{\goth{n}})x_{\alpha_1 + \alpha_2}(-1)^{m+2}.
\end{eqnarray*}
The summand $c_2x_{\alpha_1}(-1)^{n+1}x_{\alpha_1+\alpha_2}(-1)^{m+1}$
is clearly of the desired form.  For $$ c_3
x_{\alpha_1}(-2)^{k_0+k_1+k_2-n+2}x_{\alpha_1}(-1)^nx_{\alpha_1+\alpha_2}(-1)^{m+1}$$
we have that
\begin{eqnarray*}
\lefteqn{c_3x_{\alpha_1}(-2)^{k_0+k_1+k_2-n+2}x_{\alpha_1}(-1)^nx_{\alpha_1+\alpha_2}(-1)^{m+1}}\\
&=& c_3'R_{-1,2(k_0+k_1+k_2+2-n)+n}^1 + c_3''x_{\alpha_1}(-1)^{n+1}x_{\alpha_1+\alpha_2}(-1)^{m+1}
\end{eqnarray*}
for some $c_3',c_3'' \in U(\bar{\goth{n}})$, which is of the desired
form. We also have, by Lemma \ref{rootcommute}, that
\begin{eqnarray*}
\lefteqn{c_4x_{\alpha_2}(-1)^{k_0+k_1+k_2-m+1}x_{\alpha_1}(-1)^nx_{\alpha_1+\alpha_2}(-1)^{m+1}}
\\ &=&
c_4(x_{\alpha_1}(-1)^nx_{\alpha_2}(-1)^{k_0+k_1+k_2-m+1}x_{\alpha_1+\alpha_2}(-1)^{m+1}\\ &&
+
m_1x_{\alpha_1}(-1)^{n-1}x_{\alpha_1+\alpha_2}(-2)x_{\alpha_2}(-1)^{k_0+k_1+k_2-m}x_{\alpha_1+\alpha_2}(-1)^{m+1}\\ &&
+ \dots +
m_nx_{\alpha_1+\alpha_2}(-2)^nx_{\alpha_2}(-1)^{k_0+k_1+k_2-m-n+1}x_{\alpha_1+\alpha_2}(-1)^{m+1})\\ &=&
c_4(m_0'x_{\alpha_1}(-1)^n[x_{\alpha_1}(0), \dots
  [x_{\alpha_1}(0),R_{-1,k+1}^2] \dots ] +\\ &&
m_1'x_{\alpha_1}(-1)^{n-1}[x_{\alpha_1}(0), \dots [x_{\alpha_1}(0),
    R_{-1,k+2}^2] \dots ] \\ && + \dots + m_n'[x_{\alpha_1}(0), \dots
  [x_{\alpha_1}(0),R_{-1,k + n+1}^2 ] \dots ]) +
c_4'x_{\alpha_1+\alpha_2}(-1)^{m+2}
\end{eqnarray*}
for some constants $m_1,\dots m_n,m_0', \dots m_n' \in \mathbb{C}$ and
$c_4' \in U(\bar{\goth{n}})$, which is of the desired form. Finally,
we have that
\begin{eqnarray*}
\lefteqn{c_5 x_{\alpha_1+
    \alpha_2}(-2)^{k_0+k_1+k_2-m-n+1}x_{\alpha_1}(-1)^nx_{\alpha_1+\alpha_2}(-1)^{m+1}}\\
&=& c_5' [x_{\alpha_2}(0), \dots
[x_{\alpha_2}(0),R_{2(k_0+k_1+k_2-m-n+1) + n + m+1}^1] \dots ] +
c_5''x_{\alpha_1}(-1)^{n+1}x_{\alpha_1 + \alpha_2}(-1)^{m+1}
\end{eqnarray*}
for some $c_5',c_5'' \in U(\bar{\goth{n}})$, which is of the desired
form.  Hence, we may write
\begin{eqnarray*}
a = b' + c'x_{\alpha_1}(-1)^{n+1}x_{\alpha_1+\alpha_2}(-1)^{m+1} +
d'x_{\alpha_1 + \alpha_2}(-1)^{m+2}
\end{eqnarray*}
for some $b' \in I_{(k_0+1)\Lambda_0 + k_1\Lambda_1 + k_2\Lambda_2}$
and $c',d' \in U(\bar{\goth{n}})$, completing our induction.

So, in particular, we may find $b \in I_{(k_0+1)\Lambda_0 +
  k_1\Lambda_1 + k_2\Lambda_2}$ and $c,d \in U(\bar{\goth{n}})$ such
that
\begin{eqnarray*}
a = b + cx_{\alpha_1}(-1)^{k_2 + 1}x_{\alpha_1 + \alpha_2}(-1)^{k_0+1}
+ dx_{\alpha_1+\alpha_2}(-1)^{k_0+2}.
\end{eqnarray*}
Clearly we have that $b,dx_{\alpha_1+\alpha_2}(-1)^{k_0+2} \in
I_{(k_0+1)\Lambda_0 + k_1\Lambda_1 + k_2\Lambda_2}$. We also have that
\begin{eqnarray*}
\lefteqn{cx_{\alpha_1}(-1)^{k_2 + 1}x_{\alpha_1 +
    \alpha_2}(-1)^{k_0+1}} \\
&=& c'[x_{\alpha_2}(0), \dots
[x_{\alpha_2}(0),x_{\alpha_1}(-1)^{k_0+k_2+2}] \dots ]
\end{eqnarray*}
which is an element of $I_{(k_0+1)\Lambda_0 + k_1\Lambda_1 +
  k_2\Lambda_2}$. So we have that $a \in I_{(k_0+1)\Lambda_0 +
  k_1\Lambda_1 + k_2\Lambda_2}$, a contradiction, and so
$\Lambda \ne (k_0+1)\Lambda_0 + k_1\Lambda_1 + k_2\Lambda_2$, completing our proof.

\begin{remark} \em
The proof of the main theorem is similar in structure to the proof found in
\cite{CalLM2}, in that we show that our ``minimal counterexample'' element $a$ cannot
be in $(\mathrm{Ker}f_\Lambda) \setminus I_\Lambda$ for each $\Lambda$ by eliminating each
$\Lambda$ in a certain order. In \cite{CalLM2}, it was shown that $\Lambda \neq
k\Lambda_1$, $\Lambda \neq \Lambda_0 + (k-1)\Lambda_1$, $\Lambda \neq 2\Lambda_0 +
(k-2)\Lambda_1, \dots, \Lambda \neq k\Lambda_0$, in that order. This choice was incredibly
important in the proof. In the proof above, the order in which the possible $\Lambda$
are eliminated is equally important, but more choices have to be made. We show that 
$\Lambda \neq
k\Lambda_i$, $\Lambda \neq \Lambda_0 + (k-1)\Lambda_i$, $\Lambda \neq 2\Lambda_0 +
(k-2)\Lambda_i, \dots, \Lambda \neq k\Lambda_0$ for $i=1,2$. Then, we proceed to show that 
$\Lambda \neq k_1\Lambda_1 + k_2\Lambda_2$ for each $k_1,k_2 \in \mathbb{N}$ satisfying
 $k_1 + k_2 = k$, and this can only be shown once all 
$\Lambda = k_0\Lambda_0 + k_i\Lambda_i$ with $k_0, k_i \in \mathbb{N}$ satisfying
 $k_0 + k_i = k$, $i=1,2$ have been ruled out. Once we have
 $\Lambda \neq k_1\Lambda_1 + k_2\Lambda_2$, we
 then proceed to show that $\Lambda \neq k_0\Lambda_0 + k_1\Lambda_1 + k_2\Lambda_2$
for all remaining choices of $k_0,k_1,k_2 \in \mathbb{N}$ with $k_0+k_1+k_2=k$ 
by taking $k_0=1, k_0=2, \dots$ until all remaining choices of $\Lambda$ 
have been eliminated.
\end{remark}

\begin{remark} \em
The part of the proof which considers weights of the form $k\Lambda_i$
and $k_0\Lambda_0 + k_i\Lambda_i$ for $i=0,1,2$ uses generalizations
of ideas of \cite{CalLM1}-\cite{CalLM3}. These ideas no longer work in
the general case, so a new method was developed to handle the
remaining cases. This method ``rebuilds'' the element the ``minimal counterexample''
element $a$ in order to reach certain desired contradictions.
 This new method works equally well for weights of the
form $k\Lambda_i$ and $k_0\Lambda_0 + k_i\Lambda_i$. \em
\end{remark}

\begin{remark}\em
The method developed in the above proof,
which ``rebuilds'' our ``minimal counterexample'' element $a$ to show that it is actually
in the ideal $I_\Lambda$ can be used to show all the presentations considered in
\cite{CalLM1}--\cite{CalLM3} in the type $A$ case. In this sense, it is a unifying method,
and should be able to be generalized to prove presentations for the principal subspaces
of all the standard $\widehat{\goth{sl}(n+1)}$-modules, but at this stage it is 
not completely clear how this can be done in the case where $n>2$ and $k>1$. \em
\end{remark}

\begin{remark}\em
Unlike in \cite{CalLM2} and \cite{C3}, we use only intertwining operators for level $1$ standard
modules to prove certain inclusions of kernels inside other kernels.
For example, when we we wanted to show that 
$$
\mbox{Ker}f_{(k_0+1)\Lambda_0 + k_1\Lambda_1 + k_2\Lambda_2} \subset 
\mbox{Ker}f_{k_0\Lambda_0 + (k_1+1)\Lambda_1 + k_2\Lambda_2}
$$
in (\ref{KerContain}), we simply applied the operator 
$1^{\otimes k_0} \otimes \mathcal{Y}_c(e^{\lambda_1},x)\otimes 1^{\otimes k_1+k_2}$ 
to 
$$
a\cdot v_{(k_0+1)\Lambda_0 + k_1\Lambda_1 + k_2\Lambda_2} = 0
$$
to obtain
$$
a \cdot v_{k_0\Lambda_0 + (k_1+1)\Lambda_1 + k_2\Lambda_2} = 0,
$$
so that if $a \in \mbox{Ker}f_{(k_0+1)\Lambda_0 + k_1\Lambda_1 + k_2\Lambda_2}$
then $a \in \mbox{Ker}f_{k_0\Lambda_0 + (k_1+1)\Lambda_1 + k_2\Lambda_2}$. Such methods work
equally well in the cases considered in \cite{CalLM2} and \cite{C3} in showing
similar inclusions. Using these types of maps gives an alternate technique for proving such inclusions
which does not require use of intertwining operators for higher level standard modules
as in \cite{CalLM2} and \cite{C3}.
\em
\end{remark}

\section{The higher rank case}
\setcounter{equation}{0}
Much of the proof above can be generalized to the case of $\widehat{\goth{sl}(n+1)}$. We
now formulate the presentations of the principal subspace of the standard $\widehat{\goth{sl}(n+1)}$-modules
and present them as a conjecture.

Consider the $\widehat{\goth{sl}(n+1)}$-subalgebra 
\begin{equation} 
\bar{\goth{n}}= \goth{n} \otimes \mathbb{C}[t, t^{-1}].
\end{equation}
Given a $\widehat{\goth{sl}(n+1)}$-module $L(\Lambda)$ of positive integral level $k$ with highest weight vector $v_\Lambda$,
the {\it principal subspace} of $L(\Lambda)$ is defined in the
obvious way:
$$W(\Lambda) = U(\bar{\goth{n}})\cdot v_\Lambda.$$
As in the case of $n=2$ above, we define the formal sums
\begin{equation} 
R_{t}^i= \sum_{m_1+\dots+m_{k+1}=-t}x_{\alpha_i}(m_1)\cdots x_{\alpha_i}(m_{k+1}), \ \ t
\in \mathbb{Z}, \ \ i=1,\dots, n.
\end{equation}
We also define their truncations
\begin{equation} 
R_{-1, t}^i= \sum_{m_1+\dots +m_{k+1}=-t, \ m_i \le -1}x_{\alpha_i}(m_1)x_{\alpha_i}(m_2)\cdots x_{\alpha_i}(m_{k+1})
\end{equation}
for $t \in \mathbb{Z}$, $t \ge k+1$, and $i=1, \dots , n$.
Let $J$ be the left
ideal of $U(\bar{\goth{n}})$ generated by the elements $R_{-1, t}^i$
for $t \geq k+1$ and $i=1,\dots, n$:
\begin{equation}
J = \sum_{i=1}^n \sum_{t \geq k+1} U(\bar{\goth{n}}) R_{-1,t}^i.
\end{equation}
 We define certain left ideals of $U(\bar{\goth{n}})$ by:
\begin{equation*} 
I_{k\Lambda_{0}}= J+U(\bar{\goth{n}})\bar{\goth{n}}_{+}
\end{equation*}
and
\begin{eqnarray*} 
I_\Lambda = I_{k\Lambda_0} + \sum_{\alpha \in \Delta_+} U(\bar{\goth{n}})x_{\alpha}(-1)^{k+1 - \langle \alpha, \Lambda \rangle}
\end{eqnarray*}
for any $\Lambda = k_0 \Lambda_0 + \dots + k_n \Lambda_n$, where $k_0, \dots ,k_n \in \mathbb{N}$,
$k_0 + \dots + k_n = k$, and $\Lambda_0, \dots ,\Lambda_n$ are the fundamental weights of
$\widehat{\goth{sl}(n+1)}$.
For each such $\Lambda$, we have a surjective map
\begin{eqnarray} 
F_{\Lambda}: U(\widehat{\goth{g}}) & \longrightarrow &
L(\Lambda) \\ a &\mapsto& a \cdot v_{\Lambda}. \nonumber
\end{eqnarray}
and its surjective restriction $f_{\Lambda}$:
\begin{eqnarray}
f_{\Lambda}: U(\bar{\goth{n}}) & \longrightarrow &
W(\Lambda)\\ a & \mapsto & a \cdot v_{\Lambda}. \nonumber
\end{eqnarray}
We are now ready to give our conjectured presentations:
\begin{conjecture}
For each $\Lambda = k_0 \Lambda_0 + \dots + k_n \Lambda_n$ with $k_0, \dots, k_n \in \mathbb{N}$ and
$k_0 + \dots + k_n = k$, we have that
$$\mathrm{Ker}f_\Lambda = I_\Lambda.$$
\end{conjecture}

\begin{remark}\em
Many of the terms in the above conjecture may be redundant, depending on the weight $\Lambda$.
For example, in the case when $n=2$ and $\Lambda = \Lambda_0 + \Lambda_1$, we have that
$$I_\Lambda = I_{2\Lambda_0} + U(\bar{\goth{n}})x_{\alpha_1}(-1)^{2} +
U(\bar{\goth{n}})x_{\alpha_2}(-1)^3 + U(\bar{\goth{n}})x_{\alpha_1+\alpha_2}(-1)^2.$$
In this case,
$$x_{\alpha_2}(-1)^3 = cR_{-1,3}^2 \in I_{2\Lambda_0}$$
for some constant $c \in \mathbb{C}$ and
$$x_{\alpha_1+\alpha_2}(-1)^2 = d [x_{\alpha_2}(0), [x_{\alpha_2}(0),x_{\alpha_1}(-1)^2]]
 \in I_{2\Lambda_0} + U(\bar{\goth{n}})x_{\alpha_1}(-1)^2$$
for some constant $d \in \mathbb{C}$.
So we have that $$I_\Lambda \subset I_{2\Lambda_0} + U(\bar{\goth{n}})x_{\alpha_1}(-1)^2.$$
In some cases, however, this conjecture does not produce redundant terms. For $\Lambda = \Lambda_1 + \Lambda_2$,
we have that
$$I_\Lambda = I_{2\Lambda_0} + U(\bar{\goth{n}})x_{\alpha_1}(-1)^2
+U(\bar{\goth{n}})x_{\alpha_2}(-1)^2 + U(\bar{\goth{n}})x_{\alpha_1+\alpha_2}(-1),$$
which does not contain redundant terms. \em
\end{remark}

\vspace{.4in}
  
\noindent {\small \sc Department of Mathematics, Rutgers University,
Piscataway, NJ 08854}\\ 
{\em E--mail address}: sadowski@math.rutgers.edu

 \end{document}